\documentclass[a4paper,12pt]{report}
\usepackage{amssymb}
\usepackage[english,francais]{babel}

\newcommand{\theoremname}{Theorem}
\addto\captionsfrench{\renewcommand{\theoremname}{Th\'eor\`eme}}
\addto\captionsenglish{\renewcommand{\theoremname}{Theorem}}

\usepackage{amsthm,amssymb}
\newtheorem{thm}{\theoremname}[chapter]
\newtheorem{proposition}{Proposition}[chapter]{\bf}{\it}
\newtheorem{remark}{Remarque}[chapter]{\it}{\rm}

\catcode`@=11 \catcode`!=11

\expandafter\ifx\csname fiverm\endcsname\relax
  \let\fiverm\fivrm
\fi
  
\let\!latexendpicture=\endpicture 
\let\!latexframe=\frame
\let\!latexlinethickness=\linethickness
\let\!latexmultiput=\multiput
\let\!latexput=\put
 
\def\@picture(#1,#2)(#3,#4){%
  \@picht #2\unitlength
  \setbox\@picbox\hbox to #1\unitlength\bgroup 
  \let\endpicture=\!latexendpicture
  \let\frame=\!latexframe
  \let\linethickness=\!latexlinethickness
  \let\multiput=\!latexmultiput
  \let\put=\!latexput
  \hskip -#3\unitlength \lower #4\unitlength \hbox\bgroup}

\catcode`@=12 \catcode`!=12
\font\fiverm=cmr5

\catcode`!=11 
 
\def\PiC{P\kern-.12em\lower.5ex\hbox{I}\kern-.075emC}
\def\PiCTeX{\PiC\kern-.11em\TeX}

\def\!ifnextchar#1#2#3{%
  \let\!testchar=#1%
  \def\!first{#2}%
  \def\!second{#3}%
  \futurelet\!nextchar\!testnext}
\def\!testnext{%
  \ifx \!nextchar \!spacetoken 
    \let\!next=\!skipspacetestagain
  \else
    \ifx \!nextchar \!testchar
      \let\!next=\!first
    \else 
      \let\!next=\!second 
    \fi 
  \fi
  \!next}
\def\\{\!skipspacetestagain} 
  \expandafter\def\\ {\futurelet\!nextchar\!testnext} 
\def\\{\let\!spacetoken= } \\  

\def\!tfor#1:=#2\do#3{%
  \edef\!fortemp{#2}%
  \ifx\!fortemp\!empty 
    \else
    \!tforloop#2\!nil\!nil\!!#1{#3}%
  \fi}
\def\!tforloop#1#2\!!#3#4{%
  \def#3{#1}%
  \ifx #3\!nnil
    \let\!nextwhile=\!fornoop
  \else
    #4\relax
    \let\!nextwhile=\!tforloop
  \fi 
  \!nextwhile#2\!!#3{#4}}

\def\!etfor#1:=#2\do#3{%
  \def\!!tfor{\!tfor#1:=}%
  \edef\!!!tfor{#2}%
  \expandafter\!!tfor\!!!tfor\do{#3}}

\def\!cfor#1:=#2\do#3{%
  \edef\!fortemp{#2}%
  \ifx\!fortemp\!empty 
  \else
    \!cforloop#2,\!nil,\!nil\!!#1{#3}%
  \fi}
\def\!cforloop#1,#2\!!#3#4{%
  \def#3{#1}%
  \ifx #3\!nnil
    \let\!nextwhile=\!fornoop 
  \else
    #4\relax
    \let\!nextwhile=\!cforloop
  \fi
  \!nextwhile#2\!!#3{#4}}

\def\!ecfor#1:=#2\do#3{%
  \def\!!cfor{\!cfor#1:=}%
  \edef\!!!cfor{#2}%
  \expandafter\!!cfor\!!!cfor\do{#3}}

\def\!empty{}
\def\!nnil{\!nil}
\def\!fornoop#1\!!#2#3{}

\def\!ifempty#1#2#3{%
  \edef\!emptyarg{#1}%
  \ifx\!emptyarg\!empty
    #2%
  \else
    #3%
  \fi}
 
\def\!getnext#1\from#2{%
  \expandafter\!gnext#2\!#1#2}%
\def\!gnext\\#1#2\!#3#4{%
  \def#3{#1}%
  \def#4{#2\\{#1}}%
  \ignorespaces}

\def\!getnextvalueof#1\from#2{%
  \expandafter\!gnextv#2\!#1#2}%
\def\!gnextv\\#1#2\!#3#4{%
  #3=#1%
  \def#4{#2\\{#1}}%
  \ignorespaces}

\def\!copylist#1\to#2{%
  \expandafter\!!copylist#1\!#2}
\def\!!copylist#1\!#2{%
  \def#2{#1}\ignorespaces}

\def\!wlet#1=#2{%
  \let#1=#2 
  \wlog{\string#1=\string#2}}
 
\def\!listaddon#1#2{%
  \expandafter\!!listaddon#2\!{#1}#2}
\def\!!listaddon#1\!#2#3{%
  \def#3{#1\\#2}}
 
\def\!rightappend#1\withCS#2\to#3{\expandafter\!!rightappend#3\!#2{#1}#3}
\def\!!rightappend#1\!#2#3#4{\def#4{#1#2{#3}}}

\def\!leftappend#1\withCS#2\to#3{\expandafter\!!leftappend#3\!#2{#1}#3}
\def\!!leftappend#1\!#2#3#4{\def#4{#2{#3}#1}}

\def\!lop#1\to#2{\expandafter\!!lop#1\!#1#2}
\def\!!lop\\#1#2\!#3#4{\def#4{#1}\def#3{#2}}

\def\!loop#1\repeat{\def\!body{#1}\!iterate}
\def\!iterate{\!body\let\!next=\!iterate\else\let\!next=\relax\fi\!next}
 
\def\!!loop#1\repeat{\def\!!body{#1}\!!iterate}
\def\!!iterate{\!!body\let\!!next=\!!iterate\else\let\!!next=\relax\fi\!!next}
 
\def\!removept#1#2{\edef#2{\expandafter\!!removePT\the#1}}
{\catcode`p=12 \catcode`t=12 \gdef\!!removePT#1pt{#1}}

\def\placevalueinpts of <#1> in #2 {%
  \!removept{#1}{#2}}
 
\def\!mlap#1{\hbox to 0pt{\hss#1\hss}}
\def\!vmlap#1{\vbox to 0pt{\vss#1\vss}}
 
\def\!not#1{%
  #1\relax
    \!switchfalse
  \else
    \!switchtrue
  \fi
  \if!switch
  \ignorespaces}

\let\!!!wlog=\wlog              
\def\wlog#1{}    

\newdimen\headingtoplotskip     
\newdimen\linethickness         
\newdimen\longticklength        
\newdimen\plotsymbolspacing     
\newdimen\shortticklength       
\newdimen\stackleading          
\newdimen\tickstovaluesleading  
\newdimen\totalarclength        
\newdimen\valuestolabelleading  

\newbox\!boxA                   
\newbox\!boxB                   
\newbox\!picbox                 
\newbox\!plotsymbol             
\newbox\!putobject              
\newbox\!shadesymbol            

\newcount\!countA               
\newcount\!countB               
\newcount\!countC               
\newcount\!countD               
\newcount\!countE               
\newcount\!countF               
\newcount\!countG               
\newcount\!fiftypt              
\newcount\!intervalno           
\newcount\!npoints              
\newcount\!nsegments            
\newcount\!ntemp                
\newcount\!parity               
\newcount\!scalefactor          
\newcount\!tfs                  
\newcount\!tickcase             

\newdimen\!Xleft                
\newdimen\!Xright               
\newdimen\!Xsave                
\newdimen\!Ybot                 
\newdimen\!Ysave                
\newdimen\!Ytop                 
\newdimen\!angle                
\newdimen\!arclength            
\newdimen\!areabloc             
\newdimen\!arealloc             
\newdimen\!arearloc             
\newdimen\!areatloc             
\newdimen\!bshrinkage           
\newdimen\!checkbot             
\newdimen\!checkleft            
\newdimen\!checkright           
\newdimen\!checktop             
\newdimen\!dimenA               
\newdimen\!dimenB               
\newdimen\!dimenC               
\newdimen\!dimenD               
\newdimen\!dimenE               
\newdimen\!dimenF               
\newdimen\!dimenG               
\newdimen\!dimenH               
\newdimen\!dimenI               
\newdimen\!distacross           
\newdimen\!downlength           
\newdimen\!dp                   
\newdimen\!dshade               
\newdimen\!dxpos                
\newdimen\!dxprime              
\newdimen\!dypos                
\newdimen\!dyprime              
\newdimen\!ht                   
\newdimen\!leaderlength         
\newdimen\!lshrinkage           
\newdimen\!midarclength         
\newdimen\!offset               
\newdimen\!plotheadingoffset    
\newdimen\!plotsymbolxshift     
\newdimen\!plotsymbolyshift     
\newdimen\!plotxorigin          
\newdimen\!plotyorigin          
\newdimen\!rootten              
\newdimen\!rshrinkage           
\newdimen\!shadesymbolxshift    
\newdimen\!shadesymbolyshift    
\newdimen\!tenAa                
\newdimen\!tenAc                
\newdimen\!tenAe                
\newdimen\!tshrinkage           
\newdimen\!uplength             
\newdimen\!wd                   
\newdimen\!wmax                 
\newdimen\!wmin                 
\newdimen\!xB                   
\newdimen\!xC                   
\newdimen\!xE                   
\newdimen\!xM                   
\newdimen\!xS                   
\newdimen\!xaxislength          
\newdimen\!xdiff                
\newdimen\!xleft                
\newdimen\!xloc                 
\newdimen\!xorigin              
\newdimen\!xpivot               
\newdimen\!xpos                 
\newdimen\!xprime               
\newdimen\!xright               
\newdimen\!xshade               
\newdimen\!xshift               
\newdimen\!xtemp                
\newdimen\!xunit                
\newdimen\!xxE                  
\newdimen\!xxM                  
\newdimen\!xxS                  
\newdimen\!xxloc                
\newdimen\!yB                   
\newdimen\!yC                   
\newdimen\!yE                   
\newdimen\!yM                   
\newdimen\!yS                   
\newdimen\!yaxislength          
\newdimen\!ybot                 
\newdimen\!ydiff                
\newdimen\!yloc                 
\newdimen\!yorigin              
\newdimen\!ypivot               
\newdimen\!ypos                 
\newdimen\!yprime               
\newdimen\!yshade               
\newdimen\!yshift               
\newdimen\!ytemp                
\newdimen\!ytop                 
\newdimen\!yunit                
\newdimen\!yyE                  
\newdimen\!yyM                  
\newdimen\!yyS                  
\newdimen\!yyloc                
\newdimen\!zpt                  

\newif\if!axisvisible           
\newif\if!gridlinestoo          
\newif\if!keepPO                
\newif\if!placeaxislabel        
\newif\if!switch                
\newif\if!xswitch               

\newtoks\!axisLaBeL             
\newtoks\!keywordtoks           

\newwrite\!replotfile           

\newhelp\!keywordhelp{The keyword mentioned in the error message in unknown. 
Replace NEW KEYWORD in the indicated response by the keyword that 
should have been specified.}    

\!wlet\!!origin=\!xM                   
\!wlet\!!unit=\!uplength               
\!wlet\!Lresiduallength=\!dimenG       
\!wlet\!Rresiduallength=\!dimenF       
\!wlet\!axisLength=\!distacross        
\!wlet\!axisend=\!ydiff                
\!wlet\!axisstart=\!xdiff              
\!wlet\!axisxlevel=\!arclength         
\!wlet\!axisylevel=\!downlength        
\!wlet\!beta=\!dimenE                  
\!wlet\!gamma=\!dimenF                 
\!wlet\!shadexorigin=\!plotxorigin     
\!wlet\!shadeyorigin=\!plotyorigin     
\!wlet\!ticklength=\!xS                
\!wlet\!ticklocation=\!xE              
\!wlet\!ticklocationincr=\!yE          
\!wlet\!tickwidth=\!yS                 
\!wlet\!totalleaderlength=\!dimenE     
\!wlet\!xone=\!xprime                  
\!wlet\!xtwo=\!dxprime                 
\!wlet\!ySsave=\!yM                    
\!wlet\!ybB=\!yB                       
\!wlet\!ybC=\!yC                       
\!wlet\!ybE=\!yE                       
\!wlet\!ybM=\!yM                       
\!wlet\!ybS=\!yS                       
\!wlet\!ybpos=\!yyloc                  
\!wlet\!yone=\!yprime                  
\!wlet\!ytB=\!xB                       
\!wlet\!ytC=\!xC                       
\!wlet\!ytE=\!downlength               
\!wlet\!ytM=\!arclength                
\!wlet\!ytS=\!distacross               
\!wlet\!ytpos=\!xxloc                  
\!wlet\!ytwo=\!dyprime                 

\!zpt=0pt                              
\!xunit=1pt
\!yunit=1pt
\!arearloc=\!xunit
\!areatloc=\!yunit
\!dshade=5pt
\!leaderlength=24in
\!tfs=256                              
\!wmax=5.3pt                           
\!wmin=2.7pt                           
\!xaxislength=\!xunit
\!xpivot=\!zpt
\!yaxislength=\!yunit 
\!ypivot=\!zpt
  \!dimenA=50pt \!fiftypt=\!dimenA     

\!rootten=3.162278pt                   
\!tenAa=8.690286pt                     
\!tenAc=2.773839pt                     
\!tenAe=2.543275pt                     

\def\!cosrotationangle{1}      
\def\!sinrotationangle{0}      
\def\!xpivotcoord{0}           
\def\!xref{0}                  
\def\!xshadesave{0}            
\def\!ypivotcoord{0}           
\def\!yref{0}                  
\def\!yshadesave{0}            
\def\!zero{0}                  

\let\wlog=\!!!wlog
\def\normalgraphs{%
  \longticklength=.4\baselineskip
  \shortticklength=.25\baselineskip
  \tickstovaluesleading=.25\baselineskip
  \valuestolabelleading=.8\baselineskip
  \linethickness=.4pt
  \stackleading=.17\baselineskip
  \headingtoplotskip=1.5\baselineskip
  \visibleaxes
  \ticksout
  \nogridlines
  \unloggedticks}
\def\setplotarea x from #1 to #2, y from #3 to #4 {%
  \!arealloc=\!M{#1}\!xunit \advance \!arealloc -\!xorigin
  \!areabloc=\!M{#3}\!yunit \advance \!areabloc -\!yorigin
  \!arearloc=\!M{#2}\!xunit \advance \!arearloc -\!xorigin
  \!areatloc=\!M{#4}\!yunit \advance \!areatloc -\!yorigin
  \!initinboundscheck
  \!xaxislength=\!arearloc  \advance\!xaxislength -\!arealloc
  \!yaxislength=\!areatloc  \advance\!yaxislength -\!areabloc
  \!plotheadingoffset=\!zpt
  \!dimenput {{\setbox0=\hbox{}\wd0=\!xaxislength\ht0=\!yaxislength\box0}}
     [bl] (\!arealloc,\!areabloc)}
\def\visibleaxes{%
  \def\!axisvisibility{\!axisvisibletrue}}

\def\!fixkeyword#1{%
  \errhelp=\!keywordhelp
  \errmessage{Unrecognized keyword `#1': \the\!keywordtoks{NEW KEYWORD}'}}

\!keywordtoks={enter `i\fixkeyword}

\def\fixkeyword#1{%
  \!nextkeyword#1 }

\def\axis {%
  \def\!nextkeyword##1 {%
    \expandafter\ifx\csname !axis##1\endcsname \relax
      \def\!next{\!fixkeyword{##1}}%
    \else
      \def\!next{\csname !axis##1\endcsname}%
    \fi
    \!next}%
  \!offset=\!zpt
  \!axisvisibility
  \!placeaxislabelfalse
  \!nextkeyword}

\def\!axisbottom{%
  \!axisylevel=\!areabloc
  \def\!tickxsign{0}%
  \def\!tickysign{-}%
  \def\!axissetup{\!axisxsetup}%
  \def\!axislabeltbrl{t}%
  \!nextkeyword}

\def\!axistop{%
  \!axisylevel=\!areatloc
  \def\!tickxsign{0}%
  \def\!tickysign{+}%
  \def\!axissetup{\!axisxsetup}%
  \def\!axislabeltbrl{b}%
  \!nextkeyword}

\def\!axisleft{%
  \!axisxlevel=\!arealloc
  \def\!tickxsign{-}%
  \def\!tickysign{0}%
  \def\!axissetup{\!axisysetup}%
  \def\!axislabeltbrl{r}%
  \!nextkeyword}

\def\!axisright{%
  \!axisxlevel=\!arearloc
  \def\!tickxsign{+}%
  \def\!tickysign{0}%
  \def\!axissetup{\!axisysetup}%
  \def\!axislabeltbrl{l}%
  \!nextkeyword}

\def\!axisshiftedto#1=#2 {%
  \if 0\!tickxsign
    \!axisylevel=\!M{#2}\!yunit
    \advance\!axisylevel -\!yorigin
  \else
    \!axisxlevel=\!M{#2}\!xunit
    \advance\!axisxlevel -\!xorigin
  \fi
  \!nextkeyword}

\def\!axisvisible{%
  \!axisvisibletrue  
  \!nextkeyword}

\def\!axisinvisible{%
  \!axisvisiblefalse
  \!nextkeyword}

\def\!axislabel#1 {%
  \!axisLaBeL={#1}%
  \!placeaxislabeltrue
  \!nextkeyword}

\expandafter\def\csname !axis/\endcsname{%
  \!axissetup 
  \if!placeaxislabel
    \!placeaxislabel
  \fi
  \if +\!tickysign 
    \!dimenA=\!axisylevel
    \advance\!dimenA \!offset 
    \advance\!dimenA -\!areatloc 
    \ifdim \!dimenA>\!plotheadingoffset
      \!plotheadingoffset=\!dimenA 
    \fi
  \fi}

\def\grid #1 #2 {%
  \!countA=#1\advance\!countA 1
  \axis bottom invisible ticks length <\!zpt> andacross quantity {\!countA} /
  \!countA=#2\advance\!countA 1
  \axis left   invisible ticks length <\!zpt> andacross quantity {\!countA} / }

\def\plotheading#1 {%
  \advance\!plotheadingoffset \headingtoplotskip
  \!dimenput {#1} [B] <.5\!xaxislength,\!plotheadingoffset>
    (\!arealloc,\!areatloc)}

\def\!axisxsetup{%
  \!axisxlevel=\!arealloc
  \!axisstart=\!arealloc
  \!axisend=\!arearloc
  \!axisLength=\!xaxislength
  \!!origin=\!xorigin
  \!!unit=\!xunit
  \!xswitchtrue
  \if!axisvisible 
    \!makeaxis
  \fi}

\def\!axisysetup{%
  \!axisylevel=\!areabloc
  \!axisstart=\!areabloc
  \!axisend=\!areatloc
  \!axisLength=\!yaxislength
  \!!origin=\!yorigin
  \!!unit=\!yunit
  \!xswitchfalse
  \if!axisvisible
    \!makeaxis
  \fi}

\def\!makeaxis{%
  \setbox\!boxA=\hbox{
    \beginpicture
      \!setdimenmode
      \setcoordinatesystem point at {\!zpt} {\!zpt}   
      \putrule from {\!zpt} {\!zpt} to
        {\!tickysign\!tickysign\!axisLength} 
        {\!tickxsign\!tickxsign\!axisLength}
    \endpicturesave <\!Xsave,\!Ysave>}%
    \wd\!boxA=\!zpt
    \!placetick\!axisstart}

\def\!placeaxislabel{%
  \advance\!offset \valuestolabelleading
  \if!xswitch
    \!dimenput {\the\!axisLaBeL} [\!axislabeltbrl]
      <.5\!axisLength,\!tickysign\!offset> (\!axisxlevel,\!axisylevel)
    \advance\!offset \!dp  
    \advance\!offset \!ht  
  \else
    \!dimenput {\the\!axisLaBeL} [\!axislabeltbrl]
      <\!tickxsign\!offset,.5\!axisLength> (\!axisxlevel,\!axisylevel)
  \fi
  \!axisLaBeL={}}

\def\arrow <#1> [#2,#3]{%
  \!ifnextchar<{\!arrow{#1}{#2}{#3}}{\!arrow{#1}{#2}{#3}<\!zpt,\!zpt> }}

\def\!arrow#1#2#3<#4,#5> from #6 #7 to #8 #9 {%
  \!xloc=\!M{#8}\!xunit   
  \!yloc=\!M{#9}\!yunit
  \!dxpos=\!xloc  \!dimenA=\!M{#6}\!xunit  \advance \!dxpos -\!dimenA
  \!dypos=\!yloc  \!dimenA=\!M{#7}\!yunit  \advance \!dypos -\!dimenA
  \let\!MAH=\!M
  \!setdimenmode
  \!xshift=#4\relax  \!yshift=#5\relax
  \!reverserotateonly\!xshift\!yshift
  \advance\!xshift\!xloc  \advance\!yshift\!yloc
  \!xS=-\!dxpos  \advance\!xS\!xshift
  \!yS=-\!dypos  \advance\!yS\!yshift
  \!start (\!xS,\!yS)
  \!ljoin (\!xshift,\!yshift)
  \!Pythag\!dxpos\!dypos\!arclength
  \!divide\!dxpos\!arclength\!dxpos  
  \!dxpos=32\!dxpos  \!removept\!dxpos\!!cos
  \!divide\!dypos\!arclength\!dypos  
  \!dypos=32\!dypos  \!removept\!dypos\!!sin
  \!halfhead{#1}{#2}{#3}
  \!halfhead{#1}{-#2}{-#3}
  \let\!M=\!MAH
  \ignorespaces}
  \def\!halfhead#1#2#3{%
    \!dimenC=-#1%
    \divide \!dimenC 2 
    \!dimenD=#2\!dimenC
    \!rotate(\!dimenC,\!dimenD)by(\!!cos,\!!sin)to(\!xM,\!yM)
    \!dimenC=-#1
    \!dimenD=#3\!dimenC
    \!dimenD=.5\!dimenD
    \!rotate(\!dimenC,\!dimenD)by(\!!cos,\!!sin)to(\!xE,\!yE)
    \!start (\!xshift,\!yshift)
    \advance\!xM\!xshift  \advance\!yM\!yshift
    \advance\!xE\!xshift  \advance\!yE\!yshift
    \!qjoin (\!xM,\!yM) (\!xE,\!yE) 
    \ignorespaces}

\def\betweenarrows #1#2 from #3 #4 to #5 #6 {%
  \!xloc=\!M{#3}\!xunit  \!xxloc=\!M{#5}\!xunit%
  \!yloc=\!M{#4}\!yunit  \!yyloc=\!M{#6}\!yunit%
  \!dxpos=\!xxloc  \advance\!dxpos by -\!xloc
  \!dypos=\!yyloc  \advance\!dypos by -\!yloc
  \advance\!xloc .5\!dxpos
  \advance\!yloc .5\!dypos
  \let\!MBA=\!M
  \!setdimenmode
  \ifdim\!dypos=\!zpt
    \ifdim\!dxpos<\!zpt \!dxpos=-\!dxpos \fi
    \put {\!lrarrows{\!dxpos}{#1}}#2{} at {\!xloc} {\!yloc}
  \else
    \ifdim\!dxpos=\!zpt
      \ifdim\!dypos<\!zpt \!dypos=-\!zpt \fi
      \put {\!udarrows{\!dypos}{#1}}#2{} at {\!xloc} {\!yloc}
    \fi
  \fi
  \let\!M=\!MBA
  \ignorespaces}

\def\!lrarrows#1#2{
  {\setbox\!boxA=\hbox{$\mkern-2mu\mathord-\mkern-2mu$}%
   \setbox\!boxB=\hbox{$\leftarrow$}\!dimenE=\ht\!boxB
   \setbox\!boxB=\hbox{}\ht\!boxB=2\!dimenE
   \hbox to #1{$\mathord\leftarrow\mkern-6mu
     \cleaders\copy\!boxA\hfil
     \mkern-6mu\mathord-$%
     \kern.4em $\vcenter{\box\!boxB}$$\vcenter{\hbox{#2}}$\kern.4em
     $\mathord-\mkern-6mu
     \cleaders\copy\!boxA\hfil
     \mkern-6mu\mathord\rightarrow$}}}

\def\!udarrows#1#2{
  {\setbox\!boxB=\hbox{#2}%
   \setbox\!boxA=\hbox to \wd\!boxB{\hss$\vert$\hss}%
   \!dimenE=\ht\!boxA \advance\!dimenE \dp\!boxA \divide\!dimenE 2
   \vbox to #1{\offinterlineskip
      \vskip .05556\!dimenE
      \hbox to \wd\!boxB{\hss$\mkern.4mu\uparrow$\hss}\vskip-\!dimenE
      \cleaders\copy\!boxA\vfil
      \vskip-\!dimenE\copy\!boxA
      \vskip\!dimenE\copy\!boxB\vskip.4em
      \copy\!boxA\vskip-\!dimenE
      \cleaders\copy\!boxA\vfil
      \vskip-\!dimenE \hbox to \wd\!boxB{\hss$\mkern.4mu\downarrow$\hss}
      \vskip .05556\!dimenE}}}

\def\putbar#1breadth <#2> from #3 #4 to #5 #6 {%
  \!xloc=\!M{#3}\!xunit  \!xxloc=\!M{#5}\!xunit%
  \!yloc=\!M{#4}\!yunit  \!yyloc=\!M{#6}\!yunit%
  \!dypos=\!yyloc  \advance\!dypos by -\!yloc
  \!dimenI=#2  
  \ifdim \!dimenI=\!zpt 
    \putrule#1from {#3} {#4} to {#5} {#6} 
  \else 
    \let\!MBar=\!M
    \!setdimenmode 
    \divide\!dimenI 2
    \ifdim \!dypos=\!zpt             
      \advance \!yloc -\!dimenI 
      \advance \!yyloc \!dimenI
    \else
      \advance \!xloc -\!dimenI 
      \advance \!xxloc \!dimenI
    \fi
    \putrectangle#1corners at {\!xloc} {\!yloc} and {\!xxloc} {\!yyloc}
    \let\!M=\!MBar 
  \fi
  \ignorespaces}

\def\setbars#1breadth <#2> baseline at #3 = #4 {%
  \edef\!barshift{#1}%
  \edef\!barbreadth{#2}%
  \edef\!barorientation{#3}%
  \edef\!barbaseline{#4}%
  \def\!bardobaselabel{\!bardoendlabel}%
  \def\!bardoendlabel{\!barfinish}%
  \let\!drawcurve=\!barcurve
  \!setbars}
\def\!setbars{%
  \futurelet\!nextchar\!!setbars}
\def\!!setbars{%
  \if b\!nextchar
    \def\!!!setbars{\!setbarsbget}%
  \else 
    \if e\!nextchar
      \def\!!!setbars{\!setbarseget}%
    \else
      \def\!!!setbars{\relax}%
    \fi
  \fi
  \!!!setbars}
\def\!setbarsbget baselabels (#1) {%
  \def\!barbaselabelorientation{#1}%
  \def\!bardobaselabel{\!!bardobaselabel}%
  \!setbars}
\def\!setbarseget endlabels (#1) {%
  \edef\!barendlabelorientation{#1}%
  \def\!bardoendlabel{\!!bardoendlabel}%
  \!setbars}

\def\!barcurve #1 #2 {%
  \if y\!barorientation
    \def\!basexarg{#1}%
    \def\!baseyarg{\!barbaseline}%
  \else
    \def\!basexarg{\!barbaseline}%
    \def\!baseyarg{#2}%
  \fi
  \expandafter\putbar\!barshift breadth <\!barbreadth> from {\!basexarg}
    {\!baseyarg} to {#1} {#2}
  \def\!endxarg{#1}%
  \def\!endyarg{#2}%
  \!bardobaselabel}

\def\!!bardobaselabel "#1" {%
  \put {#1}\!barbaselabelorientation{} at {\!basexarg} {\!baseyarg}
  \!bardoendlabel}
 
\def\!!bardoendlabel "#1" {%
  \put {#1}\!barendlabelorientation{} at {\!endxarg} {\!endyarg}
  \!barfinish}

\def\!barfinish{%
  \!ifnextchar/{\!finish}{\!barcurve}}

\def\putrectangle{%
  \!ifnextchar<{\!putrectangle}{\!putrectangle<\!zpt,\!zpt> }}
\def\!putrectangle<#1,#2> corners at #3 #4 and #5 #6 {%
  \!xone=\!M{#3}\!xunit  \!xtwo=\!M{#5}\!xunit%
  \!yone=\!M{#4}\!yunit  \!ytwo=\!M{#6}\!yunit%
  \ifdim \!xtwo<\!xone
    \!dimenI=\!xone  \!xone=\!xtwo  \!xtwo=\!dimenI
  \fi
  \ifdim \!ytwo<\!yone
    \!dimenI=\!yone  \!yone=\!ytwo  \!ytwo=\!dimenI
  \fi
  \!dimenI=#1\relax  \advance\!xone\!dimenI  \advance\!xtwo\!dimenI
  \!dimenI=#2\relax  \advance\!yone\!dimenI  \advance\!ytwo\!dimenI
  \let\!MRect=\!M
  \!setdimenmode
  \!shaderectangle
  \!dimenI=.5\linethickness
  \advance \!xone  -\!dimenI
  \advance \!xtwo   \!dimenI
  \putrule from {\!xone} {\!yone} to {\!xtwo} {\!yone} 
  \putrule from {\!xone} {\!ytwo} to {\!xtwo} {\!ytwo} 
  \advance \!xone   \!dimenI
  \advance \!xtwo  -\!dimenI%
  \advance \!yone  -\!dimenI
  \advance \!ytwo   \!dimenI
  \putrule from {\!xone} {\!yone} to {\!xone} {\!ytwo} 
  \putrule from {\!xtwo} {\!yone} to {\!xtwo} {\!ytwo} 
  \let\!M=\!MRect
  \ignorespaces}
 
\def\shaderectanglesoff{%
  \def\!shaderectangle{}%
  \ignorespaces}

\shaderectanglesoff
 
\def\!!shaderectangle{%
  \!dimenA=\!xtwo  \advance \!dimenA -\!xone
  \!dimenB=\!ytwo  \advance \!dimenB -\!yone
  \ifdim \!dimenA<\!dimenB
    \!startvshade (\!xone,\!yone,\!ytwo)
    \!lshade      (\!xtwo,\!yone,\!ytwo)
  \else
    \!starthshade (\!yone,\!xone,\!xtwo)
    \!lshade      (\!ytwo,\!xone,\!xtwo)
  \fi
  \ignorespaces}
  
\def\frame{%
  \!ifnextchar<{\!frame}{\!frame<\!zpt> }}
\long\def\!frame<#1> #2{%
  \beginpicture
    \setcoordinatesystem units <1pt,1pt> point at 0 0 
    \put {#2} [Bl] at 0 0 
    \!dimenA=#1\relax
    \!dimenB=\!wd \advance \!dimenB \!dimenA
    \!dimenC=\!ht \advance \!dimenC \!dimenA
    \!dimenD=\!dp \advance \!dimenD \!dimenA
    \let\!MFr=\!M
    \!setdimenmode
    \putrectangle corners at {-\!dimenA} {-\!dimenD} and {\!dimenB} {\!dimenC}
    \!setcoordmode
    \let\!M=\!MFr
  \endpicture
  \ignorespaces}
 
\def\rectangle <#1> <#2> {%
  \setbox0=\hbox{}\wd0=#1\ht0=#2\frame {\box0}}

\def\plot{%
  \!ifnextchar"{\!plotfromfile}{\!drawcurve}}
\def\!plotfromfile"#1"{%
  \expandafter\!drawcurve \input #1 /}

\def\setquadratic{%
  \let\!drawcurve=\!qcurve
  \let\!!Shade=\!!qShade
  \let\!!!Shade=\!!!qShade}

\def\setlinear{%
  \let\!drawcurve=\!lcurve
  \let\!!Shade=\!!lShade
  \let\!!!Shade=\!!!lShade}

\def\sethistograms{%
  \let\!drawcurve=\!hcurve}

\def\!qcurve #1 #2 {%
  \!start (#1,#2)
  \!Qjoin}
\def\!Qjoin#1 #2 #3 #4 {%
  \!qjoin (#1,#2) (#3,#4)             
  \!ifnextchar/{\!finish}{\!Qjoin}}

\def\!lcurve #1 #2 {%
  \!start (#1,#2)
  \!Ljoin}
\def\!Ljoin#1 #2 {%
  \!ljoin (#1,#2)                    
  \!ifnextchar/{\!finish}{\!Ljoin}}

\def\!finish/{\ignorespaces}

\def\!hcurve #1 #2 {%
  \edef\!hxS{#1}%
  \edef\!hyS{#2}%
  \!hjoin}
\def\!hjoin#1 #2 {%
  \putrectangle corners at {\!hxS} {\!hyS} and {#1} {#2}
  \edef\!hxS{#1}%
  \!ifnextchar/{\!finish}{\!hjoin}}

\def\vshade #1 #2 #3 {%
  \!startvshade (#1,#2,#3)
  \!Shadewhat}

\def\hshade #1 #2 #3 {%
  \!starthshade (#1,#2,#3)
  \!Shadewhat}

\def\!Shadewhat{%
  \futurelet\!nextchar\!Shade}
\def\!Shade{%
  \if <\!nextchar
    \def\!nextShade{\!!Shade}%
  \else
    \if /\!nextchar
      \def\!nextShade{\!finish}%
    \else
      \def\!nextShade{\!!!Shade}%
    \fi
  \fi
  \!nextShade}
\def\!!lShade<#1> #2 #3 #4 {%
  \!lshade <#1> (#2,#3,#4)                 
  \!Shadewhat}
\def\!!!lShade#1 #2 #3 {%
  \!lshade (#1,#2,#3)
  \!Shadewhat} 
\def\!!qShade<#1> #2 #3 #4 #5 #6 #7 {%
  \!qshade <#1> (#2,#3,#4) (#5,#6,#7)      
  \!Shadewhat}
\def\!!!qShade#1 #2 #3 #4 #5 #6 {%
  \!qshade (#1,#2,#3) (#4,#5,#6)
  \!Shadewhat} 

\setlinear

\def\setdashpattern <#1>{%
  \def\!Flist{}\def\!Blist{}\def\!UDlist{}%
  \!countA=0
  \!ecfor\!item:=#1\do{%
    \!dimenA=\!item\relax
    \expandafter\!rightappend\the\!dimenA\withCS{\\}\to\!UDlist%
    \advance\!countA  1
    \ifodd\!countA
      \expandafter\!rightappend\the\!dimenA\withCS{\!Rule}\to\!Flist%
      \expandafter\!leftappend\the\!dimenA\withCS{\!Rule}\to\!Blist%
    \else 
      \expandafter\!rightappend\the\!dimenA\withCS{\!Skip}\to\!Flist%
      \expandafter\!leftappend\the\!dimenA\withCS{\!Skip}\to\!Blist%
    \fi}%
  \!leaderlength=\!zpt
  \def\!Rule##1{\advance\!leaderlength  ##1}%
  \def\!Skip##1{\advance\!leaderlength  ##1}%
  \!Flist%
  \ifdim\!leaderlength>\!zpt 
  \else
    \def\!Flist{\!Skip{24in}}\def\!Blist{\!Skip{24in}}\ignorespaces
    \def\!UDlist{\\{\!zpt}\\{24in}}\ignorespaces
    \!leaderlength=24in
  \fi
  \!dashingon}

\def\!dashingon{%
  \def\!advancedashing{\!!advancedashing}%
  \def\!drawlinearsegment{\!lineardashed}%
  \def\!puthline{\!putdashedhline}%
  \def\!putvline{\!putdashedvline}%
  \ignorespaces}%
\def\!dashingoff{%
  \def\!advancedashing{\relax}%
  \def\!drawlinearsegment{\!linearsolid}%
  \def\!puthline{\!putsolidhline}%
  \def\!putvline{\!putsolidvline}%
  \ignorespaces}

\def\setdots{%
  \!ifnextchar<{\!setdots}{\!setdots<5pt>}}
\def\!setdots<#1>{%
  \!dimenB=#1\advance\!dimenB -\plotsymbolspacing
  \ifdim\!dimenB<\!zpt
    \!dimenB=\!zpt
  \fi
\setdashpattern <\plotsymbolspacing,\!dimenB>}
 
\def\setdotsnear <#1> for <#2>{%
  \!dimenB=#2\relax  \advance\!dimenB -.05pt  
  \!dimenC=#1\relax  \!countA=\!dimenC 
  \!dimenD=\!dimenB  \advance\!dimenD .5\!dimenC  \!countB=\!dimenD
  \divide \!countB  \!countA
  \ifnum 1>\!countB 
    \!countB=1
  \fi
  \divide\!dimenB  \!countB
  \setdots <\!dimenB>}
 
\def\setdashes{%
  \!ifnextchar<{\!setdashes}{\!setdashes<5pt>}}
\def\!setdashes<#1>{\setdashpattern <#1,#1>}
 
\def\setdashesnear <#1> for <#2>{%
  \!dimenB=#2\relax  
  \!dimenC=#1\relax  \!countA=\!dimenC 
  \!dimenD=\!dimenB  \advance\!dimenD .5\!dimenC  \!countB=\!dimenD
  \divide \!countB  \!countA
  \ifodd \!countB 
  \else 
    \advance \!countB  1
  \fi
  \divide\!dimenB  \!countB
  \setdashes <\!dimenB>}
 
\def\setsolid{%
  \def\!Flist{\!Rule{24in}}\def\!Blist{\!Rule{24in}}%
  \def\!UDlist{\\{24in}\\{\!zpt}}%
  \!dashingoff}  
\setsolid

\def\!divide#1#2#3{%
  \!dimenB=#1
  \!dimenC=#2
  \!dimenD=\!dimenB
  \divide \!dimenD \!dimenC
  \!dimenA=\!dimenD
  \multiply\!dimenD \!dimenC
  \advance\!dimenB -\!dimenD
  \!dimenD=\!dimenC
    \ifdim\!dimenD<\!zpt \!dimenD=-\!dimenD 
  \fi
  \ifdim\!dimenD<64pt
    \!divstep[\!tfs]\!divstep[\!tfs]%
  \else 
    \!!divide
  \fi
  #3=\!dimenA\ignorespaces}

\def\!!divide{%
  \ifdim\!dimenD<256pt
    \!divstep[64]\!divstep[32]\!divstep[32]%
  \else 
    \!divstep[8]\!divstep[8]\!divstep[8]\!divstep[8]\!divstep[8]%
    \!dimenA=2\!dimenA
  \fi}

\def\!divstep[#1]{
  \!dimenB=#1\!dimenB
  \!dimenD=\!dimenB
    \divide \!dimenD by \!dimenC
  \!dimenA=#1\!dimenA
    \advance\!dimenA by \!dimenD%
  \multiply\!dimenD by \!dimenC
    \advance\!dimenB by -\!dimenD}
 
\def\Divide <#1> by <#2> forming <#3> {%
  \!divide{#1}{#2}{#3}}

\def\circulararc{%
  \ellipticalarc axes ratio 1:1 }

\def\ellipticalarc axes ratio #1:#2 #3 degrees from #4 #5 center at #6 #7 {%
  \!angle=#3pt\relax
  \ifdim\!angle>\!zpt 
    \def\!sign{}
  \else 
    \def\!sign{-}\!angle=-\!angle
  \fi
  \!xxloc=\!M{#6}\!xunit
  \!yyloc=\!M{#7}\!yunit     
  \!xxS=\!M{#4}\!xunit
  \!yyS=\!M{#5}\!yunit
  \advance\!xxS -\!xxloc
  \advance\!yyS -\!yyloc
  \!divide\!xxS{#1pt}\!xxS 
  \!divide\!yyS{#2pt}\!yyS 
  \let\!MC=\!M
  \!setdimenmode
  \!xS=#1\!xxS  \advance\!xS\!xxloc
  \!yS=#2\!yyS  \advance\!yS\!yyloc
  \!start (\!xS,\!yS)%
  \!loop\ifdim\!angle>14.9999pt
    \!rotate(\!xxS,\!yyS)by(\!cos,\!sign\!sin)to(\!xxM,\!yyM) 
    \!rotate(\!xxM,\!yyM)by(\!cos,\!sign\!sin)to(\!xxE,\!yyE)
    \!xM=#1\!xxM  \advance\!xM\!xxloc  \!yM=#2\!yyM  \advance\!yM\!yyloc
    \!xE=#1\!xxE  \advance\!xE\!xxloc  \!yE=#2\!yyE  \advance\!yE\!yyloc
    \!qjoin (\!xM,\!yM) (\!xE,\!yE)
    \!xxS=\!xxE  \!yyS=\!yyE 
    \advance \!angle -15pt
  \repeat
  \ifdim\!angle>\!zpt
    \!angle=100.53096\!angle
    \divide \!angle 360 
    \!sinandcos\!angle\!!sin\!!cos
    \!rotate(\!xxS,\!yyS)by(\!!cos,\!sign\!!sin)to(\!xxM,\!yyM) 
    \!rotate(\!xxM,\!yyM)by(\!!cos,\!sign\!!sin)to(\!xxE,\!yyE)
    \!xM=#1\!xxM  \advance\!xM\!xxloc  \!yM=#2\!yyM  \advance\!yM\!yyloc
    \!xE=#1\!xxE  \advance\!xE\!xxloc  \!yE=#2\!yyE  \advance\!yE\!yyloc
    \!qjoin (\!xM,\!yM) (\!xE,\!yE)
  \fi
  \let\!M=\!MC
  \ignorespaces}

\def\!rotate(#1,#2)by(#3,#4)to(#5,#6){%
  \!dimenA=#3#1\advance \!dimenA -#4#2
  \!dimenB=#3#2\advance \!dimenB  #4#1
  \divide \!dimenA 32  \divide \!dimenB 32 
  #5=\!dimenA  #6=\!dimenB
  \ignorespaces}
\def\!sin{4.17684}
\def\!cos{31.72624}

\def\!sinandcos#1#2#3{%
 \!dimenD=#1
 \!dimenA=\!dimenD
 \!dimenB=32pt
 \!removept\!dimenD\!value
 \!dimenC=\!dimenD
 \!dimenC=\!value\!dimenC \divide\!dimenC by 64 
 \advance\!dimenB by -\!dimenC
 \!dimenC=\!value\!dimenC \divide\!dimenC by 96 
 \advance\!dimenA by -\!dimenC
 \!dimenC=\!value\!dimenC \divide\!dimenC by 128 
 \advance\!dimenB by \!dimenC%
 \!removept\!dimenA#2
 \!removept\!dimenB#3
 \ignorespaces}

\def\putrule#1from #2 #3 to #4 #5 {%
  \!xloc=\!M{#2}\!xunit  \!xxloc=\!M{#4}\!xunit%
  \!yloc=\!M{#3}\!yunit  \!yyloc=\!M{#5}\!yunit%
  \!dxpos=\!xxloc  \advance\!dxpos by -\!xloc
  \!dypos=\!yyloc  \advance\!dypos by -\!yloc
  \ifdim\!dypos=\!zpt
    \def\!!Line{\!puthline{#1}}\ignorespaces
  \else
    \ifdim\!dxpos=\!zpt
      \def\!!Line{\!putvline{#1}}\ignorespaces
    \else 
       \def\!!Line{}
    \fi
  \fi
  \let\!ML=\!M
  \!setdimenmode
  \!!Line%
  \let\!M=\!ML
  \ignorespaces}

\def\!putsolidhline#1{%
  \ifdim\!dxpos>\!zpt 
    \put{\!hline\!dxpos}#1[l] at {\!xloc} {\!yloc}
  \else 
    \put{\!hline{-\!dxpos}}#1[l] at {\!xxloc} {\!yyloc}
  \fi
  \ignorespaces}
 
\def\!putsolidvline#1{%
  \ifdim\!dypos>\!zpt 
    \put{\!vline\!dypos}#1[b] at {\!xloc} {\!yloc}
  \else 
    \put{\!vline{-\!dypos}}#1[b] at {\!xxloc} {\!yyloc}
  \fi
  \ignorespaces}
 
\def\!hline#1{\hbox to #1{\leaders \hrule height\linethickness\hfill}}
\def\!vline#1{\vbox to #1{\leaders \vrule width\linethickness\vfill}}

\def\!putdashedhline#1{%
  \ifdim\!dxpos>\!zpt 
    \!DLsetup\!Flist\!dxpos
    \put{\hbox to \!totalleaderlength{\!hleaders}\!hpartialpattern\!Rtrunc}
      #1[l] at {\!xloc} {\!yloc} 
  \else 
    \!DLsetup\!Blist{-\!dxpos}
    \put{\!hpartialpattern\!Ltrunc\hbox to \!totalleaderlength{\!hleaders}}
      #1[r] at {\!xloc} {\!yloc} 
  \fi
  \ignorespaces}
 
\def\!putdashedvline#1{%
  \!dypos=-\!dypos
  \ifdim\!dypos>\!zpt 
    \!DLsetup\!Flist\!dypos 
    \put{\vbox{\vbox to \!totalleaderlength{\!vleaders}
      \!vpartialpattern\!Rtrunc}}#1[t] at {\!xloc} {\!yloc} 
  \else 
    \!DLsetup\!Blist{-\!dypos}
    \put{\vbox{\!vpartialpattern\!Ltrunc
      \vbox to \!totalleaderlength{\!vleaders}}}#1[b] at {\!xloc} {\!yloc} 
  \fi
  \ignorespaces}

\def\!DLsetup#1#2{
  \let\!RSlist=#1
  \!countB=#2
  \!countA=\!leaderlength
  \divide\!countB by \!countA
  \!totalleaderlength=\!countB\!leaderlength
  \!Rresiduallength=#2%
  \advance \!Rresiduallength by -\!totalleaderlength
  \!Lresiduallength=\!leaderlength
  \advance \!Lresiduallength by -\!Rresiduallength
  \ignorespaces}
 
\def\!hleaders{%
  \def\!Rule##1{\vrule height\linethickness width##1}%
  \def\!Skip##1{\hskip##1}%
  \leaders\hbox{\!RSlist}\hfill}
 
\def\!hpartialpattern#1{%
  \!dimenA=\!zpt \!dimenB=\!zpt 
  \def\!Rule##1{#1{##1}\vrule height\linethickness width\!dimenD}%
  \def\!Skip##1{#1{##1}\hskip\!dimenD}%
  \!RSlist}
 
\def\!vleaders{%
  \def\!Rule##1{\hrule width\linethickness height##1}%
  \def\!Skip##1{\vskip##1}%
  \leaders\vbox{\!RSlist}\vfill}
 
\def\!vpartialpattern#1{%
  \!dimenA=\!zpt \!dimenB=\!zpt 
  \def\!Rule##1{#1{##1}\hrule width\linethickness height\!dimenD}%
  \def\!Skip##1{#1{##1}\vskip\!dimenD}%
  \!RSlist}
 
\def\!Rtrunc#1{\!trunc{#1}>\!Rresiduallength}
\def\!Ltrunc#1{\!trunc{#1}<\!Lresiduallength}
 
\def\!trunc#1#2#3{%
  \!dimenA=\!dimenB         
  \advance\!dimenB by #1%
  \!dimenD=\!dimenB  \ifdim\!dimenD#2#3\!dimenD=#3\fi
  \!dimenC=\!dimenA  \ifdim\!dimenC#2#3\!dimenC=#3\fi
  \advance \!dimenD by -\!dimenC}

\def\!start (#1,#2){%
  \!plotxorigin=\!xorigin  \advance \!plotxorigin by \!plotsymbolxshift
  \!plotyorigin=\!yorigin  \advance \!plotyorigin by \!plotsymbolyshift
  \!xS=\!M{#1}\!xunit \!yS=\!M{#2}\!yunit
  \!rotateaboutpivot\!xS\!yS
  \!copylist\!UDlist\to\!!UDlist
  \!getnextvalueof\!downlength\from\!!UDlist
  \!distacross=\!zpt
  \!intervalno=0 
  \global\totalarclength=\!zpt
  \ignorespaces}

\def\!ljoin (#1,#2){%
  \advance\!intervalno by 1
  \!xE=\!M{#1}\!xunit \!yE=\!M{#2}\!yunit
  \!rotateaboutpivot\!xE\!yE
  \!xdiff=\!xE \advance \!xdiff by -\!xS
  \!ydiff=\!yE \advance \!ydiff by -\!yS
  \!Pythag\!xdiff\!ydiff\!arclength
  \global\advance \totalarclength by \!arclength%
  \!drawlinearsegment
  \!xS=\!xE \!yS=\!yE
  \ignorespaces}

\def\!linearsolid{%
  \!npoints=\!arclength
  \!countA=\plotsymbolspacing
  \divide\!npoints by \!countA
  \ifnum \!npoints<1 
    \!npoints=1 
  \fi
  \divide\!xdiff by \!npoints
  \divide\!ydiff by \!npoints
  \!xpos=\!xS \!ypos=\!yS
  \loop\ifnum\!npoints>-1
    \!plotifinbounds
    \advance \!xpos by \!xdiff
    \advance \!ypos by \!ydiff
    \advance \!npoints by -1
  \repeat
  \ignorespaces}

\def\!lineardashed{%
  \ifdim\!distacross>\!arclength
    \advance \!distacross by -\!arclength  
  \else
    \loop\ifdim\!distacross<\!arclength
      \!divide\!distacross\!arclength\!dimenA
      \!removept\!dimenA\!t
      \!xpos=\!t\!xdiff \advance \!xpos by \!xS
      \!ypos=\!t\!ydiff \advance \!ypos by \!yS
      \!plotifinbounds
      \advance\!distacross by \plotsymbolspacing
      \!advancedashing
    \repeat  
    \advance \!distacross by -\!arclength
  \fi
  \ignorespaces}

\def\!!advancedashing{%
  \advance\!downlength by -\plotsymbolspacing
  \ifdim \!downlength>\!zpt
  \else
    \advance\!distacross by \!downlength
    \!getnextvalueof\!uplength\from\!!UDlist
    \advance\!distacross by \!uplength
    \!getnextvalueof\!downlength\from\!!UDlist
  \fi}

\def\inboundscheckoff{%
  \def\!plotifinbounds{\!plot(\!xpos,\!ypos)}%
  \def\!initinboundscheck{\relax}\ignorespaces}
 
\inboundscheckoff
 
\def\!!plotifinbounds{%
  \ifdim \!xpos<\!checkleft
  \else
    \ifdim \!xpos>\!checkright
    \else
      \ifdim \!ypos<\!checkbot
      \else
         \ifdim \!ypos>\!checktop
         \else
           \!plot(\!xpos,\!ypos)
         \fi 
      \fi
    \fi
  \fi}

\def\!!initinboundscheck{%
  \!checkleft=\!arealloc     \advance\!checkleft by \!xorigin
  \!checkright=\!arearloc    \advance\!checkright by \!xorigin
  \!checkbot=\!areabloc      \advance\!checkbot by \!yorigin
  \!checktop=\!areatloc      \advance\!checktop by \!yorigin}

\def\!logten#1#2{%
  \expandafter\!!logten#1\!nil
  \!removept\!dimenF#2%
  \ignorespaces}

\def\!!logten#1#2\!nil{%
  \if -#1%
    \!dimenF=\!zpt
    \def\!next{\ignorespaces}%
  \else
    \if +#1%
      \def\!next{\!!logten#2\!nil}%
    \else
      \if .#1%
        \def\!next{\!!logten0.#2\!nil}%
      \else
        \def\!next{\!!!logten#1#2..\!nil}%
      \fi
    \fi
  \fi
  \!next}

\def\!!!logten#1#2.#3.#4\!nil{%
  \!dimenF=1pt 
  \if 0#1%
    \!!logshift#3pt 
  \else 
    \!logshift#2/
    \!dimenE=#1.#2#3pt 
  \fi 
  \ifdim \!dimenE<\!rootten
    \multiply \!dimenE 10 
    \advance  \!dimenF -1pt
  \fi
  \!dimenG=\!dimenE
    \advance\!dimenG 10pt
  \advance\!dimenE -10pt 
  \multiply\!dimenE 10 
  \!divide\!dimenE\!dimenG\!dimenE
  \!removept\!dimenE\!t
  \!dimenG=\!t\!dimenE
  \!removept\!dimenG\!tt
  \!dimenH=\!tt\!tenAe
    \divide\!dimenH 100
  \advance\!dimenH \!tenAc
  \!dimenH=\!tt\!dimenH
    \divide\!dimenH 100   
  \advance\!dimenH \!tenAa
  \!dimenH=\!t\!dimenH
    \divide\!dimenH 100 
  \advance\!dimenF \!dimenH}

\def\!logshift#1{%
  \if #1/%
    \def\!next{\ignorespaces}%
  \else
    \advance\!dimenF 1pt 
    \def\!next{\!logshift}%
  \fi 
  \!next}
 
 \def\!!logshift#1{%
   \advance\!dimenF -1pt
   \if 0#1%
     \def\!next{\!!logshift}%
   \else
     \if p#1%
       \!dimenF=1pt
       \def\!next{\!dimenE=1p}%
     \else
       \def\!next{\!dimenE=#1.}%
     \fi
   \fi
   \!next}

\def\beginpicture{%
  \setbox\!picbox=\hbox\bgroup%
  \!xleft=\maxdimen  
  \!xright=-\maxdimen
  \!ybot=\maxdimen
  \!ytop=-\maxdimen}
 
\def\endpicture{%
  \ifdim\!xleft=\maxdimen
    \!xleft=\!zpt \!xright=\!zpt \!ybot=\!zpt \!ytop=\!zpt 
  \fi
  \global\!Xleft=\!xleft \global\!Xright=\!xright
  \global\!Ybot=\!ybot \global\!Ytop=\!ytop
  \egroup%
  \ht\!picbox=\!Ytop  \dp\!picbox=-\!Ybot
  \ifdim\!Ybot>\!zpt
  \else 
    \ifdim\!Ytop<\!zpt
      \!Ybot=\!Ytop
    \else
      \!Ybot=\!zpt
    \fi
  \fi
  \hbox{\kern-\!Xleft\lower\!Ybot\box\!picbox\kern\!Xright}}
 
\def\endpicturesave <#1,#2>{%
  \endpicture \global #1=\!Xleft \global #2=\!Ybot \ignorespaces}

\def\setcoordinatesystem{%
  \!ifnextchar{u}{\!getlengths }
    {\!getlengths units <\!xunit,\!yunit>}}
\def\!getlengths units <#1,#2>{%
  \!xunit=#1\relax
  \!yunit=#2\relax
  \!ifcoordmode 
    \let\!SCnext=\!SCccheckforRP
  \else
    \let\!SCnext=\!SCdcheckforRP
  \fi
  \!SCnext}
\def\!SCccheckforRP{%
  \!ifnextchar{p}{\!cgetreference }
    {\!cgetreference point at {\!xref} {\!yref} }}
\def\!cgetreference point at #1 #2 {%
  \edef\!xref{#1}\edef\!yref{#2}%
  \!xorigin=\!xref\!xunit  \!yorigin=\!yref\!yunit  
  \!initinboundscheck 
  \ignorespaces}
\def\!SCdcheckforRP{%
  \!ifnextchar{p}{\!dgetreference}%
    {\ignorespaces}}
\def\!dgetreference point at #1 #2 {%
  \!xorigin=#1\relax  \!yorigin=#2\relax
  \ignorespaces}

\long\def\put#1#2 at #3 #4 {%
  \!setputobject{#1}{#2}%
  \!xpos=\!M{#3}\!xunit  \!ypos=\!M{#4}\!yunit  
  \!rotateaboutpivot\!xpos\!ypos%
  \advance\!xpos -\!xorigin  \advance\!xpos -\!xshift
  \advance\!ypos -\!yorigin  \advance\!ypos -\!yshift
  \kern\!xpos\raise\!ypos\box\!putobject\kern-\!xpos%
  \!doaccounting\ignorespaces}
 
\long\def\multiput #1#2 at {%
  \!setputobject{#1}{#2}%
  \!ifnextchar"{\!putfromfile}{\!multiput}}
\def\!putfromfile"#1"{%
  \expandafter\!multiput \input #1 /}
\def\!multiput{%
  \futurelet\!nextchar\!!multiput}
\def\!!multiput{%
  \if *\!nextchar
    \def\!nextput{\!alsoby}%
  \else
    \if /\!nextchar
      \def\!nextput{\!finishmultiput}%
    \else
      \def\!nextput{\!alsoat}%
    \fi
  \fi
  \!nextput}
\def\!finishmultiput/{%
  \setbox\!putobject=\hbox{}%
  \ignorespaces}
 
\def\!alsoat#1 #2 {%
  \!xpos=\!M{#1}\!xunit  \!ypos=\!M{#2}\!yunit  
  \!rotateaboutpivot\!xpos\!ypos%
  \advance\!xpos -\!xorigin  \advance\!xpos -\!xshift
  \advance\!ypos -\!yorigin  \advance\!ypos -\!yshift
  \kern\!xpos\raise\!ypos\copy\!putobject\kern-\!xpos%
  \!doaccounting
  \!multiput}
 
\def\!alsoby*#1 #2 #3 {%
  \!dxpos=\!M{#2}\!xunit \!dypos=\!M{#3}\!yunit 
  \!rotateonly\!dxpos\!dypos
  \!ntemp=#1%
  \!!loop\ifnum\!ntemp>0
    \advance\!xpos by \!dxpos  \advance\!ypos by \!dypos
    \kern\!xpos\raise\!ypos\copy\!putobject\kern-\!xpos%
    \advance\!ntemp by -1
  \repeat
  \!doaccounting 
  \!multiput}
 
\def\accountingon{\def\!doaccounting{\!!doaccounting}\ignorespaces}

\accountingon
\def\!!doaccounting{%
  \!xtemp=\!xpos  
  \!ytemp=\!ypos
  \ifdim\!xtemp<\!xleft 
     \!xleft=\!xtemp 
  \fi
  \advance\!xtemp by  \!wd 
  \ifdim\!xright<\!xtemp 
    \!xright=\!xtemp
  \fi
  \advance\!ytemp by -\!dp
  \ifdim\!ytemp<\!ybot  
    \!ybot=\!ytemp
  \fi
  \advance\!ytemp by  \!dp
  \advance\!ytemp by  \!ht 
  \ifdim\!ytemp>\!ytop  
    \!ytop=\!ytemp  
  \fi}
 
\long\def\!setputobject#1#2{%
  \setbox\!putobject=\hbox{#1}%
  \!ht=\ht\!putobject  \!dp=\dp\!putobject  \!wd=\wd\!putobject
  \wd\!putobject=\!zpt
  \!xshift=.5\!wd   \!yshift=.5\!ht   \advance\!yshift by -.5\!dp
  \edef\!putorientation{#2}%
  \expandafter\!SPOreadA\!putorientation[]\!nil%
  \expandafter\!SPOreadB\!putorientation<\!zpt,\!zpt>\!nil\ignorespaces}
 
\def\!SPOreadA#1[#2]#3\!nil{\!etfor\!orientation:=#2\do\!SPOreviseshift}
 
\def\!SPOreadB#1<#2,#3>#4\!nil{\advance\!xshift by -#2\advance\!yshift by -#3}
 
\def\!SPOreviseshift{%
  \if l\!orientation 
    \!xshift=\!zpt
  \else 
    \if r\!orientation 
      \!xshift=\!wd
    \else 
      \if b\!orientation
        \!yshift=-\!dp
      \else 
        \if B\!orientation 
          \!yshift=\!zpt
        \else 
          \if t\!orientation 
            \!yshift=\!ht
          \fi 
        \fi
      \fi
    \fi
  \fi}

\long\def\!dimenput#1#2(#3,#4){%
  \!setputobject{#1}{#2}%
  \!xpos=#3\advance\!xpos by -\!xshift
  \!ypos=#4\advance\!ypos by -\!yshift
  \kern\!xpos\raise\!ypos\box\!putobject\kern-\!xpos%
  \!doaccounting\ignorespaces}

\def\!setdimenmode{%
  \let\!M=\!M!!\ignorespaces}
\def\!setcoordmode{%
  \let\!M=\!M!\ignorespaces}
\def\!ifcoordmode{%
  \ifx \!M \!M!}
\def\!ifdimenmode{%
  \ifx \!M \!M!!}
\def\!M!#1#2{#1#2} 
\def\!M!!#1#2{#1}
\!setcoordmode
\let\setdimensionmode=\!setdimenmode
\let\setcoordinatemode=\!setcoordmode

\def\!stack[#1]{%
  \let\!lglue=\hfill \let\!rglue=\hfill
  \expandafter\let\csname !#1glue\endcsname=\relax
  \!ifnextchar<{\!!stack}{\!!stack<\stackleading>}}
\def\!!stack<#1>#2{%
  \vbox{\def\!valueslist{}\!ecfor\!value:=#2\do{%
    \expandafter\!rightappend\!value\withCS{\\}\to\!valueslist}%
    \!lop\!valueslist\to\!value
    \let\\=\cr\lineskiplimit=\maxdimen\lineskip=#1%
    \baselineskip=-1000pt\halign{\!lglue##\!rglue\cr \!value\!valueslist\cr}}%
  \ignorespaces}

\def\!lines[#1]#2{%
  \let\!lglue=\hfill \let\!rglue=\hfill
  \expandafter\let\csname !#1glue\endcsname=\relax
  \vbox{\halign{\!lglue##\!rglue\cr #2\crcr}}%
  \ignorespaces}

\def\!Lines[#1]#2{%
  \let\!lglue=\hfill \let\!rglue=\hfill
  \expandafter\let\csname !#1glue\endcsname=\relax
  \vtop{\halign{\!lglue##\!rglue\cr #2\crcr}}%
  \ignorespaces}

\def\setplotsymbol(#1#2){%
  \!setputobject{#1}{#2}
  \setbox\!plotsymbol=\box\!putobject%
  \!plotsymbolxshift=\!xshift 
  \!plotsymbolyshift=\!yshift 
  \ignorespaces}
 
\setplotsymbol({\fiverm .})

\def\!!plot(#1,#2){%
  \!dimenA=-\!plotxorigin \advance \!dimenA by #1
  \!dimenB=-\!plotyorigin \advance \!dimenB by #2
  \kern\!dimenA\raise\!dimenB\copy\!plotsymbol\kern-\!dimenA%
  \ignorespaces}
 
\def\!!!plot(#1,#2){%
  \!dimenA=-\!plotxorigin \advance \!dimenA by #1
  \!dimenB=-\!plotyorigin \advance \!dimenB by #2
  \kern\!dimenA\raise\!dimenB\copy\!plotsymbol\kern-\!dimenA%
  \!countE=\!dimenA
  \!countF=\!dimenB
  \immediate\write\!replotfile{\the\!countE,\the\!countF.}%
  \ignorespaces}

\def\savelinesandcurves on "#1" {%
  \immediate\closeout\!replotfile
  \immediate\openout\!replotfile=#1%
  \let\!plot=\!!!plot}

\def\dontsavelinesandcurves {%
  \let\!plot=\!!plot}
\dontsavelinesandcurves

{\catcode`\%=11\xdef\!Commentsignal{
\def\writesavefile#1 {%
  \immediate\write\!replotfile{\!Commentsignal #1}%
  \ignorespaces}

\def\replot"#1" {%
  \expandafter\!replot\input #1 /}
\def\!replot#1,#2. {%
  \!dimenA=#1sp
  \kern\!dimenA\raise#2sp\copy\!plotsymbol\kern-\!dimenA
  \futurelet\!nextchar\!!replot}
\def\!!replot{%
  \if /\!nextchar 
    \def\!next{\!finish}%
  \else
    \def\!next{\!replot}%
  \fi
  \!next}

\def\!Pythag#1#2#3{%
  \!dimenE=#1\relax                                     
  \ifdim\!dimenE<\!zpt 
    \!dimenE=-\!dimenE 
  \fi
  \!dimenF=#2\relax
  \ifdim\!dimenF<\!zpt 
    \!dimenF=-\!dimenF 
  \fi
  \advance \!dimenF by \!dimenE
  \ifdim\!dimenF=\!zpt 
    \!dimenG=\!zpt
  \else 
    \!divide{8\!dimenE}\!dimenF\!dimenE
    \advance\!dimenE by -4pt
      \!dimenE=2\!dimenE
    \!removept\!dimenE\!!t
    \!dimenE=\!!t\!dimenE
    \advance\!dimenE by 64pt
    \divide \!dimenE by 2
    \!dimenH=7pt
    \!!Pythag\!!Pythag\!!Pythag
    \!removept\!dimenH\!!t
    \!dimenG=\!!t\!dimenF
    \divide\!dimenG by 8
  \fi
  #3=\!dimenG
  \ignorespaces}

\def\!!Pythag{
  \!divide\!dimenE\!dimenH\!dimenI
  \advance\!dimenH by \!dimenI
    \divide\!dimenH by 2}

\def\placehypotenuse for <#1> and <#2> in <#3> {%
  \!Pythag{#1}{#2}{#3}}

\def\!qjoin (#1,#2) (#3,#4){%
  \advance\!intervalno by 1
  \!ifcoordmode
    \edef\!xmidpt{#1}\edef\!ymidpt{#2}%
  \else
    \!dimenA=#1\relax \edef\!xmidpt{\the\!dimenA}%
    \!dimenA=#2\relax \edef\!ymidpt{\the\!dimenA}%
  \fi
  \!xM=\!M{#1}\!xunit  \!yM=\!M{#2}\!yunit   \!rotateaboutpivot\!xM\!yM
  \!xE=\!M{#3}\!xunit  \!yE=\!M{#4}\!yunit   \!rotateaboutpivot\!xE\!yE
  \!dimenA=\!xM  \advance \!dimenA by -\!xS
  \!dimenB=\!xE  \advance \!dimenB by -\!xM
  \!xB=3\!dimenA \advance \!xB by -\!dimenB
  \!xC=2\!dimenB \advance \!xC by -2\!dimenA
  \!dimenA=\!yM  \advance \!dimenA by -\!yS%
  \!dimenB=\!yE  \advance \!dimenB by -\!yM%
  \!yB=3\!dimenA \advance \!yB by -\!dimenB%
  \!yC=2\!dimenB \advance \!yC by -2\!dimenA%
  \!xprime=\!xB  \!yprime=\!yB
  \!dxprime=.5\!xC  \!dyprime=.5\!yC
  \!getf \!midarclength=\!dimenA
  \!getf \advance \!midarclength by 4\!dimenA
  \!getf \advance \!midarclength by \!dimenA
  \divide \!midarclength by 12
  \!arclength=\!dimenA
  \!getf \advance \!arclength by 4\!dimenA
  \!getf \advance \!arclength by \!dimenA
  \divide \!arclength by 12
  \advance \!arclength by \!midarclength
  \global\advance \totalarclength by \!arclength
  \ifdim\!distacross>\!arclength 
    \advance \!distacross by -\!arclength
  \else
    \!initinverseinterp
    \loop\ifdim\!distacross<\!arclength
      \!inverseinterp
      \!xpos=\!t\!xC \advance\!xpos by \!xB
        \!xpos=\!t\!xpos \advance \!xpos by \!xS
      \!ypos=\!t\!yC \advance\!ypos by \!yB
        \!ypos=\!t\!ypos \advance \!ypos by \!yS
      \!plotifinbounds
      \advance\!distacross \plotsymbolspacing
      \!advancedashing
    \repeat  
    \advance \!distacross by -\!arclength
  \fi
  \!xS=\!xE
  \!yS=\!yE
  \ignorespaces}

\def\!getf{\!Pythag\!xprime\!yprime\!dimenA%
  \advance\!xprime by \!dxprime
  \advance\!yprime by \!dyprime}

\def\!initinverseinterp{%
  \ifdim\!arclength>\!zpt
    \!divide{8\!midarclength}\!arclength\!dimenE
    \ifdim\!dimenE<\!wmin \!setinverselinear
    \else 
      \ifdim\!dimenE>\!wmax \!setinverselinear
      \else
        \def\!inverseinterp{\!inversequad}\ignorespaces
         \!removept\!dimenE\!Ew
         \!dimenF=-\!Ew\!dimenE
         \advance\!dimenF by 32pt
         \!dimenG=8pt 
         \advance\!dimenG by -\!dimenE
         \!dimenG=\!Ew\!dimenG
         \!divide\!dimenF\!dimenG\!beta
         \!gamma=1pt
         \advance \!gamma by -\!beta
      \fi
    \fi
  \fi
  \ignorespaces}

\def\!inversequad{%
  \!divide\!distacross\!arclength\!dimenG
  \!removept\!dimenG\!v
  \!dimenG=\!v\!gamma
  \advance\!dimenG by \!beta
  \!dimenG=\!v\!dimenG
  \!removept\!dimenG\!t}

\def\!setinverselinear{%
  \def\!inverseinterp{\!inverselinear}%
  \divide\!dimenE by 8 \!removept\!dimenE\!t
  \!countC=\!intervalno \multiply \!countC 2
  \!countB=\!countC     \advance \!countB -1
  \!countA=\!countB     \advance \!countA -1
  \wlog{\the\!countB th point (\!xmidpt,\!ymidpt) being plotted 
    doesn't lie in the}%
  \wlog{ middle third of the arc between the \the\!countA th 
    and \the\!countC th points:}%
  \wlog{ [arc length \the\!countA\space to \the\!countB]/[arc length 
    \the \!countA\space to \the\!countC]=\!t.}%
  \ignorespaces}
 
\def\!inverselinear{%
  \!divide\!distacross\!arclength\!dimenG
  \!removept\!dimenG\!t}

\def\startrotation{%
  \let\!rotateaboutpivot=\!!rotateaboutpivot
  \let\!rotateonly=\!!rotateonly
  \!ifnextchar{b}{\!getsincos }%
    {\!getsincos by {\!cosrotationangle} {\!sinrotationangle} }}
\def\!getsincos by #1 #2 {%
  \edef\!cosrotationangle{#1}%
  \edef\!sinrotationangle{#2}%
  \!ifcoordmode 
    \let\!ROnext=\!ccheckforpivot
  \else
    \let\!ROnext=\!dcheckforpivot
  \fi
  \!ROnext}
\def\!ccheckforpivot{%
  \!ifnextchar{a}{\!cgetpivot}%
    {\!cgetpivot about {\!xpivotcoord} {\!ypivotcoord} }}
\def\!cgetpivot about #1 #2 {%
  \edef\!xpivotcoord{#1}%
  \edef\!ypivotcoord{#2}%
  \!xpivot=#1\!xunit  \!ypivot=#2\!yunit
  \ignorespaces}
\def\!dcheckforpivot{%
  \!ifnextchar{a}{\!dgetpivot}{\ignorespaces}}
\def\!dgetpivot about #1 #2 {%
  \!xpivot=#1\relax  \!ypivot=#2\relax
  \ignorespaces}

\def\stoprotation{%
  \let\!rotateaboutpivot=\!!!rotateaboutpivot
  \let\!rotateonly=\!!!rotateonly
  \ignorespaces}
 
\def\!!rotateaboutpivot#1#2{%
  \!dimenA=#1\relax  \advance\!dimenA -\!xpivot
  \!dimenB=#2\relax  \advance\!dimenB -\!ypivot
  \!dimenC=\!cosrotationangle\!dimenA
    \advance \!dimenC -\!sinrotationangle\!dimenB
  \!dimenD=\!cosrotationangle\!dimenB
    \advance \!dimenD  \!sinrotationangle\!dimenA
  \advance\!dimenC \!xpivot  \advance\!dimenD \!ypivot
  #1=\!dimenC  #2=\!dimenD
  \ignorespaces}

\def\!!rotateonly#1#2{%
  \!dimenA=#1\relax  \!dimenB=#2\relax 
  \!dimenC=\!cosrotationangle\!dimenA
    \advance \!dimenC -\!rotsign\!sinrotationangle\!dimenB
  \!dimenD=\!cosrotationangle\!dimenB
    \advance \!dimenD  \!rotsign\!sinrotationangle\!dimenA
  #1=\!dimenC  #2=\!dimenD
  \ignorespaces}
\def\!rotsign{}
\def\!!!rotateaboutpivot#1#2{\relax}
\def\!!!rotateonly#1#2{\relax}
\stoprotation

\def\!reverserotateonly#1#2{%
  \def\!rotsign{-}%
  \!rotateonly{#1}{#2}%
  \def\!rotsign{}%
  \ignorespaces}

\def\!getspan span <#1>{%
  \!dshade=#1\relax
  \!ifcoordmode 
    \let\!GRnext=\!GRccheckforAP
  \else
    \let\!GRnext=\!GRdcheckforAP
  \fi
  \!GRnext}
\def\!GRccheckforAP{%
  \!ifnextchar{p}{\!cgetanchor }
    {\!cgetanchor point at {\!xshadesave} {\!yshadesave} }}
\def\!cgetanchor point at #1 #2 {%
  \edef\!xshadesave{#1}\edef\!yshadesave{#2}%
  \!xshade=\!xshadesave\!xunit  \!yshade=\!yshadesave\!yunit
  \ignorespaces}
\def\!GRdcheckforAP{%
  \!ifnextchar{p}{\!dgetanchor}%
    {\ignorespaces}}
\def\!dgetanchor point at #1 #2 {%
  \!xshade=#1\relax  \!yshade=#2\relax
  \ignorespaces}

\def\setshadesymbol{%
  \!ifnextchar<{\!setshadesymbol}{\!setshadesymbol<,,,> }}

\def\!setshadesymbol <#1,#2,#3,#4> (#5#6){%
  \!setputobject{#5}{#6}%
  \setbox\!shadesymbol=\box\!putobject%
  \!shadesymbolxshift=\!xshift \!shadesymbolyshift=\!yshift
  \!dimenA=\!xshift \advance\!dimenA \!smidge
  \!override\!dimenA{#1}\!lshrinkage%
  \!dimenA=\!wd \advance \!dimenA -\!xshift
    \advance\!dimenA \!smidge
    \!override\!dimenA{#2}\!rshrinkage
  \!dimenA=\!dp \advance \!dimenA \!yshift
    \advance\!dimenA \!smidge
    \!override\!dimenA{#3}\!bshrinkage
  \!dimenA=\!ht \advance \!dimenA -\!yshift
    \advance\!dimenA \!smidge
    \!override\!dimenA{#4}\!tshrinkage
  \ignorespaces}
\def\!smidge{-.2pt}%

\def\!override#1#2#3{%
  \edef\!!override{#2}%
  \ifx \!!override\empty
    #3=#1\relax
  \else
    \if z\!!override
      #3=\!zpt
    \else
      \ifx \!!override\!blankz
        #3=\!zpt
      \else
        #3=#2\relax
      \fi
    \fi
  \fi
  \ignorespaces}
\def\!blankz{ z}

\setshadesymbol ({\fiverm .})

\def\!startvshade#1(#2,#3,#4){%
  \let\!!xunit=\!xunit%
  \let\!!yunit=\!yunit%
  \let\!!xshade=\!xshade%
  \let\!!yshade=\!yshade%
  \def\!getshrinkages{\!vgetshrinkages}%
  \let\!setshadelocation=\!vsetshadelocation%
  \!xS=\!M{#2}\!!xunit
  \!ybS=\!M{#3}\!!yunit
  \!ytS=\!M{#4}\!!yunit
  \!shadexorigin=\!xorigin  \advance \!shadexorigin \!shadesymbolxshift
  \!shadeyorigin=\!yorigin  \advance \!shadeyorigin \!shadesymbolyshift
  \ignorespaces}
 
\def\!starthshade#1(#2,#3,#4){%
  \let\!!xunit=\!yunit%
  \let\!!yunit=\!xunit%
  \let\!!xshade=\!yshade%
  \let\!!yshade=\!xshade%
  \def\!getshrinkages{\!hgetshrinkages}%
  \let\!setshadelocation=\!hsetshadelocation%
  \!xS=\!M{#2}\!!xunit
  \!ybS=\!M{#3}\!!yunit
  \!ytS=\!M{#4}\!!yunit
  \!shadexorigin=\!xorigin  \advance \!shadexorigin \!shadesymbolxshift
  \!shadeyorigin=\!yorigin  \advance \!shadeyorigin \!shadesymbolyshift
  \ignorespaces}

\def\!lattice#1#2#3#4#5{%
  \!dimenA=#1
  \!dimenB=#2
  \!countB=\!dimenB
  \!dimenC=#3
  \advance\!dimenC -\!dimenA
  \!countA=\!dimenC
  \divide\!countA \!countB
  \ifdim\!dimenC>\!zpt
    \!dimenD=\!countA\!dimenB
    \ifdim\!dimenD<\!dimenC
      \advance\!countA 1 
    \fi
  \fi
  \!dimenC=\!countA\!dimenB
    \advance\!dimenC \!dimenA
  #4=\!countA
  #5=\!dimenC
  \ignorespaces}

\def\!qshade#1(#2,#3,#4)#5(#6,#7,#8){%
  \!xM=\!M{#2}\!!xunit
  \!ybM=\!M{#3}\!!yunit
  \!ytM=\!M{#4}\!!yunit
  \!xE=\!M{#6}\!!xunit
  \!ybE=\!M{#7}\!!yunit
  \!ytE=\!M{#8}\!!yunit
  \!getcoeffs\!xS\!ybS\!xM\!ybM\!xE\!ybE\!ybB\!ybC
  \!getcoeffs\!xS\!ytS\!xM\!ytM\!xE\!ytE\!ytB\!ytC
  \def\!getylimits{\!qgetylimits}%
  \!shade{#1}\ignorespaces}
 
\def\!lshade#1(#2,#3,#4){%
  \!xE=\!M{#2}\!!xunit
  \!ybE=\!M{#3}\!!yunit
  \!ytE=\!M{#4}\!!yunit
  \!dimenE=\!xE  \advance \!dimenE -\!xS
  \!dimenC=\!ytE \advance \!dimenC -\!ytS
  \!divide\!dimenC\!dimenE\!ytB
  \!dimenC=\!ybE \advance \!dimenC -\!ybS
  \!divide\!dimenC\!dimenE\!ybB
  \def\!getylimits{\!lgetylimits}%
  \!shade{#1}\ignorespaces}
 
\def\!getcoeffs#1#2#3#4#5#6#7#8{%
  \!dimenC=#4\advance \!dimenC -#2
  \!dimenE=#3\advance \!dimenE -#1
  \!divide\!dimenC\!dimenE\!dimenF
  \!dimenC=#6\advance \!dimenC -#4
  \!dimenH=#5\advance \!dimenH -#3
  \!divide\!dimenC\!dimenH\!dimenG
  \advance\!dimenG -\!dimenF
  \advance \!dimenH \!dimenE
  \!divide\!dimenG\!dimenH#8
  \!removept#8\!t
  #7=-\!t\!dimenE
  \advance #7\!dimenF
  \ignorespaces}

\def\!shade#1{%
  \!getshrinkages#1<,,,>\!nil
  \advance \!dimenE \!xS
  \!lattice\!!xshade\!dshade\!dimenE
    \!parity\!xpos
  \!dimenF=-\!dimenF
    \advance\!dimenF \!xE
  \!loop\!not{\ifdim\!xpos>\!dimenF}
    \!shadecolumn%
    \advance\!xpos \!dshade
    \advance\!parity 1
  \repeat
  \!xS=\!xE
  \!ybS=\!ybE
  \!ytS=\!ytE
  \ignorespaces}

\def\!vgetshrinkages#1<#2,#3,#4,#5>#6\!nil{%
  \!override\!lshrinkage{#2}\!dimenE
  \!override\!rshrinkage{#3}\!dimenF
  \!override\!bshrinkage{#4}\!dimenG
  \!override\!tshrinkage{#5}\!dimenH
  \ignorespaces}
\def\!hgetshrinkages#1<#2,#3,#4,#5>#6\!nil{%
  \!override\!lshrinkage{#2}\!dimenG
  \!override\!rshrinkage{#3}\!dimenH
  \!override\!bshrinkage{#4}\!dimenE
  \!override\!tshrinkage{#5}\!dimenF
  \ignorespaces}

\def\!shadecolumn{%
  \!dxpos=\!xpos
  \advance\!dxpos -\!xS
  \!removept\!dxpos\!dx
  \!getylimits
  \advance\!ytpos -\!dimenH
  \advance\!ybpos \!dimenG
  \!yloc=\!!yshade
  \ifodd\!parity 
     \advance\!yloc \!dshade
  \fi
  \!lattice\!yloc{2\!dshade}\!ybpos%
    \!countA\!ypos
  \!dimenA=-\!shadexorigin \advance \!dimenA \!xpos
  \loop\!not{\ifdim\!ypos>\!ytpos}
    \!setshadelocation
    \!rotateaboutpivot\!xloc\!yloc%
    \!dimenA=-\!shadexorigin \advance \!dimenA \!xloc
    \!dimenB=-\!shadeyorigin \advance \!dimenB \!yloc
    \kern\!dimenA \raise\!dimenB\copy\!shadesymbol \kern-\!dimenA
    \advance\!ypos 2\!dshade
  \repeat
  \ignorespaces}
 
\def\!qgetylimits{%
  \!dimenA=\!dx\!ytC              
  \advance\!dimenA \!ytB
  \!ytpos=\!dx\!dimenA
  \advance\!ytpos \!ytS
  \!dimenA=\!dx\!ybC              
  \advance\!dimenA \!ybB
  \!ybpos=\!dx\!dimenA
  \advance\!ybpos \!ybS}
 
\def\!lgetylimits{%
  \!ytpos=\!dx\!ytB
  \advance\!ytpos \!ytS
  \!ybpos=\!dx\!ybB
  \advance\!ybpos \!ybS}
 
\def\!vsetshadelocation{
  \!xloc=\!xpos
  \!yloc=\!ypos}
\def\!hsetshadelocation{
  \!xloc=\!ypos
  \!yloc=\!xpos}

\def\!axisticks {%
  \def\!nextkeyword##1 {%
    \expandafter\ifx\csname !ticks##1\endcsname \relax
      \def\!next{\!fixkeyword{##1}}%
    \else
      \def\!next{\csname !ticks##1\endcsname}%
    \fi
    \!next}%
  \!axissetup
    \def\!axissetup{\relax}%
  \edef\!ticksinoutsign{\!ticksinoutSign}%
  \!ticklength=\longticklength
  \!tickwidth=\linethickness
  \!gridlinestatus
  \!setticktransform
  \!maketick
  \!tickcase=0
  \def\!LTlist{}%
  \!nextkeyword}

\def\ticksout{%
  \def\!ticksinoutSign{+}}

\ticksout

\def\nogridlines{%
  \def\!gridlinestatus{\!gridlinestoofalse}}
\nogridlines

\def\loggedticks{%
  \def\!setticktransform{\let\!ticktransform=\!logten}}
\def\unloggedticks{%
  \def\!setticktransform{\let\!ticktransform=\!donothing}}
\def\!donothing#1#2{\def#2{#1}}
\unloggedticks

\expandafter\def\csname !ticks/\endcsname{%
  \!not {\ifx \!LTlist\empty}
    \!placetickvalues
  \fi
  \def\!tickvalueslist{}%
  \def\!LTlist{}%
  \expandafter\csname !axis/\endcsname}

\def\!maketick{%
  \setbox\!boxA=\hbox{%
    \beginpicture
      \!setdimenmode
      \setcoordinatesystem point at {\!zpt} {\!zpt}   
      \linethickness=\!tickwidth
      \ifdim\!ticklength>\!zpt
        \putrule from {\!zpt} {\!zpt} to
          {\!ticksinoutsign\!tickxsign\!ticklength}
          {\!ticksinoutsign\!tickysign\!ticklength}
      \fi
      \if!gridlinestoo
        \putrule from {\!zpt} {\!zpt} to
          {-\!tickxsign\!xaxislength} {-\!tickysign\!yaxislength}
      \fi
    \endpicturesave <\!Xsave,\!Ysave>}%
    \wd\!boxA=\!zpt}
  
\def\!ticksin{%
  \def\!ticksinoutsign{-}%
  \!maketick
  \!nextkeyword}

\def\!ticksout{%
  \def\!ticksinoutsign{+}%
  \!maketick
  \!nextkeyword}

\def\!tickslength<#1> {%
  \!ticklength=#1\relax
  \!maketick
  \!nextkeyword}

\def\!tickslong{%
  \!tickslength<\longticklength> }

\def\!ticksshort{%
  \!tickslength<\shortticklength> }

\def\!tickswidth<#1> {%
  \!tickwidth=#1\relax
  \!maketick
  \!nextkeyword}

\def\!ticksandacross{%
  \!gridlinestootrue
  \!maketick
  \!nextkeyword}

\def\!ticksbutnotacross{%
  \!gridlinestoofalse
  \!maketick
  \!nextkeyword}

\def\!tickslogged{%
  \let\!ticktransform=\!logten
  \!nextkeyword}

\def\!ticksunlogged{%
  \let\!ticktransform=\!donothing
  \!nextkeyword}

\def\!ticksunlabeled{%
  \!tickcase=0
  \!nextkeyword}

\def\!ticksnumbered{%
  \!tickcase=1
  \!nextkeyword}

\def\!tickswithvalues#1/ {%
  \edef\!tickvalueslist{#1! /}%
  \!tickcase=2
  \!nextkeyword}

\def\!ticksquantity#1 {%
  \ifnum #1>1
    \!updatetickoffset
    \!countA=#1\relax
    \advance \!countA -1
    \!ticklocationincr=\!axisLength
      \divide \!ticklocationincr \!countA
    \!ticklocation=\!axisstart
    \loop \!not{\ifdim \!ticklocation>\!axisend}
      \!placetick\!ticklocation
      \ifcase\!tickcase
          \relax 
        \or
          \relax 
        \or
          \expandafter\!gettickvaluefrom\!tickvalueslist
          \edef\!tickfield{{\the\!ticklocation}{\!value}}%
          \expandafter\!listaddon\expandafter{\!tickfield}\!LTlist%
      \fi
      \advance \!ticklocation \!ticklocationincr
    \repeat
  \fi
  \!nextkeyword}

\def\!ticksat#1 {%
  \!updatetickoffset
  \edef\!Loc{#1}%
  \if /\!Loc
    \def\next{\!nextkeyword}%
  \else
    \!ticksincommon
    \def\next{\!ticksat}%
  \fi
  \next}    
      
\def\!ticksfrom#1 to #2 by #3 {%
  \!updatetickoffset
  \edef\!arg{#3}%
  \expandafter\!separate\!arg\!nil
  \!scalefactor=1
  \expandafter\!countfigures\!arg/
  \edef\!arg{#1}%
  \!scaleup\!arg by\!scalefactor to\!countE
  \edef\!arg{#2}%
  \!scaleup\!arg by\!scalefactor to\!countF
  \edef\!arg{#3}%
  \!scaleup\!arg by\!scalefactor to\!countG
  \loop \!not{\ifnum\!countE>\!countF}
    \ifnum\!scalefactor=1
      \edef\!Loc{\the\!countE}%
    \else
      \!scaledown\!countE by\!scalefactor to\!Loc
    \fi
    \!ticksincommon
    \advance \!countE \!countG
  \repeat
  \!nextkeyword}

\def\!updatetickoffset{%
  \!dimenA=\!ticksinoutsign\!ticklength
  \ifdim \!dimenA>\!offset
    \!offset=\!dimenA
  \fi}

\def\!placetick#1{%
  \if!xswitch
    \!xpos=#1\relax
    \!ypos=\!axisylevel
  \else
    \!xpos=\!axisxlevel
    \!ypos=#1\relax
  \fi
  \advance\!xpos \!Xsave
  \advance\!ypos \!Ysave
  \kern\!xpos\raise\!ypos\copy\!boxA\kern-\!xpos
  \ignorespaces}

\def\!gettickvaluefrom#1 #2 /{%
  \edef\!value{#1}%
  \edef\!tickvalueslist{#2 /}%
  \ifx \!tickvalueslist\!endtickvaluelist
    \!tickcase=0
  \fi}
\def\!endtickvaluelist{! /}

\def\!ticksincommon{%
  \!ticktransform\!Loc\!t
  \!ticklocation=\!t\!!unit
  \advance\!ticklocation -\!!origin
  \!placetick\!ticklocation
  \ifcase\!tickcase
    \relax 
  \or 
    \ifdim\!ticklocation<-\!!origin
      \edef\!Loc{$\!Loc$}%
    \fi
    \edef\!tickfield{{\the\!ticklocation}{\!Loc}}%
    \expandafter\!listaddon\expandafter{\!tickfield}\!LTlist%
  \or 
    \expandafter\!gettickvaluefrom\!tickvalueslist
    \edef\!tickfield{{\the\!ticklocation}{\!value}}%
    \expandafter\!listaddon\expandafter{\!tickfield}\!LTlist%
  \fi}

\def\!separate#1\!nil{%
  \!ifnextchar{-}{\!!separate}{\!!!separate}#1\!nil}
\def\!!separate-#1\!nil{%
  \def\!sign{-}%
  \!!!!separate#1..\!nil}
\def\!!!separate#1\!nil{%
  \def\!sign{+}%
  \!!!!separate#1..\!nil}
\def\!!!!separate#1.#2.#3\!nil{%
  \def\!arg{#1}%
  \ifx\!arg\!empty
    \!countA=0
  \else
    \!countA=\!arg
  \fi
  \def\!arg{#2}%
  \ifx\!arg\!empty
    \!countB=0
  \else
    \!countB=\!arg
  \fi}
 
\def\!countfigures#1{%
  \if #1/%
    \def\!next{\ignorespaces}%
  \else
    \multiply\!scalefactor 10
    \def\!next{\!countfigures}%
  \fi
  \!next}

\def\!scaleup#1by#2to#3{%
  \expandafter\!separate#1\!nil
  \multiply\!countA #2\relax
  \advance\!countA \!countB
  \if -\!sign
    \!countA=-\!countA
  \fi
  #3=\!countA
  \ignorespaces}

\def\!scaledown#1by#2to#3{%
  \!countA=#1\relax
  \ifnum \!countA<0 
    \def\!sign{-}
    \!countA=-\!countA
  \else
    \def\!sign{}%
  \fi
  \!countB=\!countA
  \divide\!countB #2\relax
  \!countC=\!countB
    \multiply\!countC #2\relax
  \advance \!countA -\!countC
  \edef#3{\!sign\the\!countB.}
  \!countC=\!countA 
  \ifnum\!countC=0 
    \!countC=1
  \fi
  \multiply\!countC 10
  \!loop \ifnum #2>\!countC
    \edef#3{#3\!zero}%
    \multiply\!countC 10
  \repeat
  \edef#3{#3\the\!countA}
  \ignorespaces}

\def\!placetickvalues{%
  \advance\!offset \tickstovaluesleading
  \if!xswitch
    \setbox\!boxA=\hbox{%
      \def\\##1##2{%
        \!dimenput {##2} [B] (##1,\!axisylevel)}%
      \beginpicture 
        \!LTlist
      \endpicturesave <\!Xsave,\!Ysave>}%
    \!dimenA=\!axisylevel
      \advance\!dimenA -\!Ysave
      \advance\!dimenA \!tickysign\!offset
      \if -\!tickysign
        \advance\!dimenA -\ht\!boxA
      \else
        \advance\!dimenA  \dp\!boxA
      \fi
    \advance\!offset \ht\!boxA 
      \advance\!offset \dp\!boxA
    \!dimenput {\box\!boxA} [Bl] <\!Xsave,\!Ysave> (\!zpt,\!dimenA)
  \else
    \setbox\!boxA=\hbox{%
      \def\\##1##2{%
        \!dimenput {##2} [r] (\!axisxlevel,##1)}%
      \beginpicture 
        \!LTlist
      \endpicturesave <\!Xsave,\!Ysave>}%
    \!dimenA=\!axisxlevel
      \advance\!dimenA -\!Xsave
      \advance\!dimenA \!tickxsign\!offset
      \if -\!tickxsign
        \advance\!dimenA -\wd\!boxA
      \fi
    \advance\!offset \wd\!boxA
    \!dimenput {\box\!boxA} [Bl] <\!Xsave,\!Ysave> (\!dimenA,\!zpt)
  \fi}

\normalgraphs
\catcode`!=12 

\catcode`@=11 \catcode`!=11
  
\let\!pictexendpicture=\endpicture 
\let\!pictexframe=\frame
\let\!pictexlinethickness=\linethickness
\let\!pictexmultiput=\multiput
\let\!pictexput=\put

\def\beginpicture{%
  \setbox\!picbox=\hbox\bgroup%
  \let\endpicture=\!pictexendpicture
  \let\frame=\!pictexframe
  \let\linethickness=\!pictexlinethickness
  \let\multiput=\!pictexmultiput
  \let\put=\!pictexput
  \let\input=\@@input   
  \!xleft=\maxdimen  
  \!xright=-\maxdimen
  \!ybot=\maxdimen
  \!ytop=-\maxdimen}

\let\frame=\!latexframe

\let\pictexframe=\!pictexframe

\let\linethickness=\!latexlinethickness
\let\pictexlinethickness=\!pictexlinethickness

\let\\=\@normalcr
\catcode`@=12 \catcode`!=12
\catcode `\!=11
\catcode `\@=11

\let\!tacr=\\ 

\newdimen\LineThicknessUnit 
\newdimen\StrutUnit            
\newskip \InterColumnSpaceUnit  
\newdimen\ColumnWidthUnit     
\newdimen\KernUnit

\let\!taLTU=\LineThicknessUnit 
\let\!taCWU=\ColumnWidthUnit   
\let\!taKU =\KernUnit          

\newtoks\NormalTLTU
\newtoks\NormalTSU
\newtoks\NormalTICSU
\newtoks\NormalTCWU
\newtoks\NormalTKU

\NormalTLTU={1in \divide \LineThicknessUnit by 300 }
\NormalTSU ={\normalbaselineskip
  \divide \StrutUnit by 11 }  
\NormalTICSU={.5em plus 1fil minus .25em}  
\NormalTCWU ={.5em}
\NormalTKU  ={.5em}

\def\NormalTableUnits{%
  \LineThicknessUnit   =\the\NormalTLTU
  \StrutUnit           =\the\NormalTSU
  \InterColumnSpaceUnit=\the\NormalTICSU
  \ColumnWidthUnit     =\the\NormalTCWU
  \KernUnit            =\the\NormalTKU}
 
\NormalTableUnits

\newcount\@multicnt
\newcount\LineThicknessFactor    
\newcount\StrutHeightFactor      
\newcount\StrutDepthFactor       
\newcount\InterColumnSpaceFactor 
\newcount\ColumnWidthFactor      
\newcount\KernFactor
\newcount\VspaceFactor

\LineThicknessFactor    =2
\StrutHeightFactor      =8
\StrutDepthFactor       =3
\InterColumnSpaceFactor =3
\ColumnWidthFactor      =10
\KernFactor             =1
\VspaceFactor           =2

\newcount\TracingKeys 
\newcount\TracingFormats  

\def\BeginTableParBox#1{%
  \vtop\bgroup 
    \hsize=#1
    \normalbaselines 
    \let~=\!ttTie
    \let\-=\!ttDH
    \the\EveryTableParBox} 
  
\def\EndTableParBox{%
    \MakeStrut{0pt}{\StrutDepthFactor\StrutUnit}
  \egroup} 

\newtoks\EveryTableParBox
\EveryTableParBox={%
  \parindent=0pt
  \raggedright
  \rightskip=0pt plus 4em 
  \relax}

\newtoks\EveryTable
\newtoks\!taTableSpread

\newskip\LeftTabskip
\newskip\RightTabskip

\newcount\!taCountA
\newcount\!taColumnNumber
\newcount\!taRecursionLevel 

\newdimen\!taDimenA  
\newdimen\!taDimenB  
\newdimen\!taDimenC  
\newdimen\!taMinimumColumnWidth

\newtoks\!taToksA

\newtoks\!taPreamble
\newtoks\!taDataColumnTemplate
\newtoks\!taRuleColumnTemplate
\newtoks\!taOldRuleColumnTemplate
\newtoks\!taLeftGlue
\newtoks\!taRightGlue

\newskip\!taLastRegularTabskip

\newif\if!taDigit
\newif\if!taBeginFormat
\newif\if!taOnceOnlyTabskip

\def\TaBlE{%
  T\kern-.27em\lower.5ex\hbox{A}\kern-.18em B\kern-.1em
    \lower.5ex\hbox{L}\kern-.075em E}

{\catcode`\|=13 \catcode`\"=13
  \gdef\ActivateBarAndQuote{%
    \ifnum \catcode`\|=13
    \else
      \catcode`\|=13
      \def|{%
        \ifmmode
          \vert
        \else
          \char`\|
        \fi}%
    \fi
    \ifnum \catcode`\"=13
    \else
      \catcode`\"=13
      \def"{\char`\"}%
    \fi}}
 
{\catcode `\|=12 \catcode `\"=12 

}

\def\!thMessage#1{\immediate\write16{#1}\ignorespaces}
 
\let\!thx=\expandafter

\def\!thGobble#1{} 

\def\\{\let\!thSpaceToken= }\\ 

\def\!thHeight{height}
\def\!thDepth{depth}
\def\!thWidth{width}

\def\!thToksEdef#1=#2{%
  \edef\!ttemp{#2}%
  #1\!thx{\!ttemp}%
  \ignorespaces}

\def\!thStoreErrorMsg#1#2{%
  \toks0 =\!thx{\csname #2\endcsname}%
  \edef#1{\the\toks0 }}

\def\!thReadErrorMsg#1{%
  \!thx\!thx\!thx\!thGobble\!thx\string #1}

\def\!thError#1#2{%
  \begingroup
    \newlinechar=`\^^J%
    \edef\!ttemp{#2}%
    \errhelp=\!thx{\!ttemp}%
    \!thMessage{%
      ^^J\!thReadErrorMsg\!thErrorMsgA 
      ^^J\!thReadErrorMsg\!thErrorMsgB}%
    \errmessage{#1}%
  \endgroup}

\!thStoreErrorMsg\!thErrorMsgA{%
  TABLE error; see manual for explanation.}
\!thStoreErrorMsg\!thErrorMsgB{%
  Type \space H <return> \space for immediate help.}

\def\!thGetReplacement#1#2{%
   \begingroup
     \!thMessage{#1}
     \endlinechar=-1
     \global\read16 to#2%
   \endgroup}

\def\!thLoop#1\repeat{%
  \def\!thIterate{%
    #1%
    \!thx \!thIterate
    \fi}%
  \!thIterate 
  \let\!thIterate\relax}

\def\Smash{%
  \relax
  \ifmmode
    \expandafter\mathpalette
    \expandafter\!thDoMathVCS
  \else
    \expandafter\!thDoVCS
  \fi}
                      
\def\!thDoVCS#1{%
  \setbox\z@\hbox{#1}%
  \!thFinishVCS}
                      
\def\!thDoMathVCS#1#2{%
  \setbox\z@\hbox{$\m@th#1{#2}$}%
  \!thFinishVCS}
                      
\def\!thFinishVCS{%
  \vbox to\z@{\vss\box\z@\vss}}

\def\Lower{%
  \def\!thSign{-}%
  \!tgGetValue\!thSetDimen}

\def\!thSetDimen{%
  \ifnum \!tgCode=1
    \ifx \!tgValue\empty
      \!taDimenA \StrutHeightFactor\StrutUnit
      \advance \!taDimenA \StrutDepthFactor\StrutUnit
      \divide \!taDimenA 2
    \else
      \!taDimenA \!tgValue\StrutUnit
    \fi
  \else
    \!taDimenA \!tgValue
  \fi
  \!taDimenA=\!thSign\!taDimenA\relax
  \ifmmode
    \expandafter\mathpalette
    \expandafter\!thDoMathRaise
  \else
    \expandafter\!thDoSimpleRaise
  \fi}
                      
\def\!thDoSimpleRaise#1{%
  \setbox\z@\hbox{\raise \!taDimenA\hbox{#1}}%
  \!thFinishRaise} 
                      
\def\!thDoMathRaise#1#2{%
  \setbox\z@\hbox{\raise \!taDimenA\hbox{$\m@th#1{#2}$}}%
  \!thFinishRaise}

\def\!thFinishRaise{%
  \ht\z@\z@ 
  \dp\z@\z@
  \box\z@}

\def\!thKernBack{%
  \kern -
  \ifnum \!tgCode=1 
    \ifx \!tgValue\empty 
      \the\KernFactor
    \else
      \!tgValue    
    \fi
    \KernUnit
  \else 
    \!tgValue      
  \fi
  \ignorespaces}%

\def\Vspace{%
  \noalign
  \bgroup
  \!tgGetValue\!thVspace}

\def\!thVspace{%
  \vskip
    \ifnum \!tgCode=1 
      \ifx \!tgValue\empty 
        \the\VspaceFactor
      \else
        \!tgValue    
      \fi
      \StrutUnit
    \else 
      \!tgValue      
    \fi
  \egroup} 

\def\BeginFormat{%
  \catcode`\|=12 
  \catcode`\"=12 
  \!taPreamble={}%
  \!taColumnNumber=0
  \skip0 =\InterColumnSpaceUnit
  \multiply\skip0 \InterColumnSpaceFactor
  \divide\skip0 2
  \!taRuleColumnTemplate=\!thx{%
    \!thx\tabskip\the\skip0 }%
  \!taLastRegularTabskip=\skip0 
  \!taOnceOnlyTabskipfalse
  \!taBeginFormattrue 
  \def\!tfRowOfWidths{}
  \ReadFormatKeys}

\def\!tfSetWidth{%
  \ifx \!tfRowOfWidths \empty  
    \ifnum \!taColumnNumber>0  
      \begingroup              
         \!taCountA=1          
         \aftergroup \edef \aftergroup \!tfRowOfWidths \aftergroup {%
           \aftergroup &\aftergroup \omit
           \!thLoop
             \ifnum \!taCountA<\!taColumnNumber
             \advance\!taCountA 1
             \aftergroup \!tfAOAO
           \repeat 
           \aftergroup }%
      \endgroup
    \fi
  \fi      
  \ifx [\!ttemp 
    \!thx\!tfSetWidthText
  \else
    \!thx\!tfSetWidthValue
  \fi}

\def\!tfAOAO{%
  &\omit&\omit}

\def\!tfSetWidthText [#1]{
  \def\!tfWidthText{#1}%
  \ReadFormatKeys}

\def\!tfSetWidthValue{%
  \!taMinimumColumnWidth = 
    \ifnum \!tgCode=1 
      \ifx\!tgValue\empty 
        \ColumnWidthFactor
      \else
        \!tgValue 
      \fi
      \ColumnWidthUnit
    \else
      \!tgValue 
    \fi
  \def\!tfWidthText{}
  \ReadFormatKeys}

\def\!tfSetTabskip{%
  \ifnum \!tgCode=1
    \skip0 =\InterColumnSpaceUnit
    \multiply\skip0 
      \ifx \!tgValue\empty
        \InterColumnSpaceFactor         
      \else
       \!tgValue                        
      \fi
  \else
    \skip0 =\!tgValue                   
  \fi
  \divide\skip0 by 2
  \ifnum\!taColumnNumber=0 
    \!thToksEdef\!taRuleColumnTemplate={%
      \the\!taRuleColumnTemplate 
      \tabskip \the\skip0 }
  \else
    \!thToksEdef\!taDataColumnTemplate={%
      \the\!taDataColumnTemplate 
      \tabskip \the\skip0 }
  \fi
  \if!taOnceOnlyTabskip
  \else
    \!taLastRegularTabskip=\skip0 
  \fi                             
  \ReadFormatKeys}

\def\!tfSetVrule{%
  \!thToksEdef\!taRuleColumnTemplate={%
    \noexpand\hfil
    \noexpand\vrule
    \noexpand\!thWidth
    \ifnum \!tgCode=1
      \ifx \!tgValue\empty
        \the\LineThicknessFactor      
      \else
        \!tgValue                     
      \fi
      \!taLTU                         
    \else
      \!tgValue                       
    \fi
    ####%
    \noexpand\hfil
    \the\!taRuleColumnTemplate}       
  \!tfAdjoinPriorColumn}
 
\def\!tfSetAlternateVrule{%
  \afterassignment\!tfSetAlternateA
  \toks0 =}                           

\def\!tfSetAlternateA{%
  \!thToksEdef\!taRuleColumnTemplate={%
    \the\toks0 \the\!taRuleColumnTemplate} 
  \!tfAdjoinPriorColumn}

\def\!tfAdjoinPriorColumn{%
  \ifnum \!taColumnNumber=0
    \!taPreamble=\!taRuleColumnTemplate 
    \ifnum \TracingFormats>0             
      \!tfShowRuleTemplate
    \fi
  \else
    \ifx\!tfRowOfWidths\empty  
    \else
      \!tfUpdateRowOfWidths
    \fi
    \!thToksEdef\!taDataColumnTemplate={%
      \the \!taLeftGlue
      \the \!taDataColumnTemplate
      \the \!taRightGlue}
    \ifnum \TracingFormats>0
      \!tfShowTemplates
    \fi
    \!thToksEdef\!taPreamble={%
      \the\!taPreamble
      &
      \the\!taDataColumnTemplate
      &
      \the\!taRuleColumnTemplate}
  \fi
  \advance \!taColumnNumber 1
  \if!taOnceOnlyTabskip              
    \!thToksEdef\!taDataColumnTemplate={%
       ####\tabskip \the\!taLastRegularTabskip}
  \else
    \!taDataColumnTemplate{##}%
  \fi
  \!taRuleColumnTemplate{}
  \!taLeftGlue{\hfil}
  \!taRightGlue{\hfil}%
  \!taMinimumColumnWidth=0pt
  \def\!tfWidthText{}%
  \!taOnceOnlyTabskipfalse    
  \ReadFormatKeys}

\def\!tfUpdateRowOfWidths{%
  \ifx \!tfWidthText\empty
  \else 
    \!tfComputeMinColWidth
  \fi
  \edef\!tfRowOfWidths{%
    \!tfRowOfWidths
    &%
    \omit                                  
    \ifdim \!taMinimumColumnWidth>0pt
      \hskip \the\!taMinimumColumnWidth
    \fi
    &
    \omit}}                                

\def\!tfComputeMinColWidth{%
  \setbox0 =\vbox{%
    \ialign{
       \span\the\!taDataColumnTemplate\cr
       \!tfWidthText\cr}}%
  \!taMinimumColumnWidth=\wd0 }

\def\!tfShowRuleTemplate{%
  \!thMessage{}
  \!thMessage{TABLE FORMAT}
  \!thMessage{Column: Template}
  \!thMessage{%
    \space *c: ##\tabskip \the\LeftTabskip}
  \!taOldRuleColumnTemplate=\!taRuleColumnTemplate}

\def\!tfShowTemplates{%
  \!thMessage{%
    \space \space r: \the\!taOldRuleColumnTemplate}
  \!taOldRuleColumnTemplate=\!taRuleColumnTemplate
  \!thMessage{%
    \ifnum \!taColumnNumber<10
      \space
    \fi
    \the\!taColumnNumber c: \the\!taDataColumnTemplate}
  \ifdim\!taMinimumColumnWidth>0pt
    \!thMessage{%
      \space \space w: \the\!taMinimumColumnWidth}
  \fi}

\def\!tfFinishFormat{%
  \ifnum \TracingFormats>0
    \!thMessage{%
      \space \space r: \the\!taOldRuleColumnTemplate
        \tabskip \the\RightTabskip}%
    \!thMessage{%
      \space *c: ##\tabskip 0pt}
  \fi
  \ifnum \!taColumnNumber<2
    \!thError{%
      \ifnum \!taColumnNumber=0
        No
      \else
        Only 1
      \fi
      "|"}%
      {\!thReadErrorMsg\!tfTooFewBarsA
       ^^J\!thReadErrorMsg\!tfTooFewBarsB
       ^^J\!thReadErrorMsg\!tkFixIt}%
  \fi
  \!thToksEdef\!taPreamble={%
    ####\tabskip\LeftTabskip 
    &
    \the\!taPreamble \tabskip\RightTabskip
    &
    ####\tabskip 0pt \cr}
  \ifnum \TracingFormats>1
    \!thMessage{Preamble=\the\!taPreamble}
  \fi
  \ifnum \TracingFormats>2
    \!thMessage{Row Of Widths="\!tfRowOfWidths"}
  \fi
  \!taBeginFormatfalse 
  \catcode`\|=13
  \catcode`\"=13
  \!ttDoHalign}

\!thStoreErrorMsg\!tfTooFewBarsA{%
  There must be at least 2 "|"'s (and/or "\string \|"'s)}
\!thStoreErrorMsg\!tfTooFewBarsB{%
  between \string\BeginFormat\space and \string\EndFormat\space (or ".").}

\def\ReFormat[{%
  \omit
  \!taDataColumnTemplate{##}%
  \!taLeftGlue{}%
  \!taRightGlue{}%
  \catcode`\|=12  
  \catcode`\"=12  
  \ReadFormatKeys}

\def\!tfEndReFormat{%
  \ifnum \TracingFormats>0
    \!thMessage{ReF: 
       \the\!taLeftGlue
       \hbox{\the\!taDataColumnTemplate}
       \the\!taRightGlue}
  \fi
  \catcode`\|=13
  \catcode`\"=13
  \!tfReFormat}

\def\!tfReFormat#1{%
  \the \!taLeftGlue
  \vbox{%
    \ialign{%
      \span\the\!taDataColumnTemplate\cr
       #1\cr}}%
  \the \!taRightGlue}

\def\!tgGetValue#1{%
  \def\!tgReturn{#1}
  \futurelet\!ttemp\!tgCheckForParen}

\def\!tgCheckForParen{%
  \ifx\!ttemp (%
    \!thx \!tgDoParen
  \else
    \!thx \!tgCheckForSpace
  \fi}

\def\!tgDoParen(#1){%
  \def\!tgCode{2}%
  \def\!tgValue{#1}
  \!tgReturn}

\def\!tgCheckForSpace{%
  \def\!tgCode{1}%
  \def\!tgValue{}
  \ifx\!ttemp\!thSpaceToken
    \!thx \!tgReturn        
  \else
    \!thx \!tgCheckForDigit         
  \fi}

\def\!tgCheckForDigit{%
  \!taDigitfalse
  \ifx 0\!ttemp
    \!taDigittrue
  \else
    \ifx 1\!ttemp
      \!taDigittrue
    \else
      \ifx 2\!ttemp
        \!taDigittrue
      \else
        \ifx 3\!ttemp
          \!taDigittrue
        \else
          \ifx 4\!ttemp
            \!taDigittrue
          \else
            \ifx 5\!ttemp
              \!taDigittrue
            \else
              \ifx 6\!ttemp
                \!taDigittrue
              \else
                \ifx 7\!ttemp
                  \!taDigittrue
                \else
                  \ifx 8\!ttemp
                    \!taDigittrue
                  \else
                    \ifx 9\!ttemp
                      \!taDigittrue
                    \fi
                  \fi
                \fi
              \fi
            \fi
          \fi
        \fi
      \fi
    \fi
  \fi
  \if!taDigit
    \!thx \!tgGetNumber
  \else
    \!thx \!tgReturn 
  \fi}

\def\!tgGetNumber{%
  \afterassignment\!tgGetNumberA
  \!taCountA=}
\def\!tgGetNumberA{%
  \edef\!tgValue{\the\!taCountA}%
  \!tgReturn}

\def\!tgSetUpParBox{%
  \edef\!ttemp{%
    \noexpand \ReadFormatKeys
    b{\noexpand \BeginTableParBox{%
      \ifnum \!tgCode=1 
        \ifx \!tgValue\empty 
          \the\ColumnWidthFactor
        \else
          \!tgValue    
        \fi
        \!taCWU        
      \else 
        \!tgValue      
      \fi}}}%
  \!ttemp
  a{\EndTableParBox}}

\def\!tgInsertKern{%
  \edef\!ttemp{%
    \kern
    \ifnum \!tgCode=1 
      \ifx \!tgValue\empty 
        \the\KernFactor
      \else
        \!tgValue    
      \fi
      \!taKU         
    \else 
      \!tgValue      
    \fi}%
  \edef\!ttemp{%
    \noexpand\ReadFormatKeys
    \ifh@            
      b{\!ttemp}
    \fi
    \ifv@            
      a{\!ttemp}
    \fi}%
  \!ttemp}

\def\NewFormatKey#1{%
  \!thx\def\!thx\!ttempa\!thx{\string #1}%
  \!thx\def\!thx\!ttempb\!thx{\csname !tk<\!ttempa>\endcsname}%
  \ifnum \TracingKeys>0
    \!tkReportNewKey
  \fi
  \!thx\ifx \!ttempb \relax
    \!thx\!tkDefineKey
  \else 
    \!thx\!tkRejectKey
  \fi}

\def\!tkReportNewKey{%
  \!taToksA\!thx{\!ttempa}%
  \!thMessage{NEW KEY: "\the\!taToksA"}}

\def\!tkDefineKey{%
  \!thx\def\!ttempb}%

\def\!tkRejectKey{%
    \!taToksA\!thx{\!ttempa}%
    \!thError{Key letter "\the\!taToksA" already used}
      {\!thReadErrorMsg\!tkFixIt}
    \def\!tkGarbage}%

\!thStoreErrorMsg\!tkFixIt{%
  You'd better type \space 'E' \space and fix your file.}

\def\ReadFormatKeys#1{%
  \!thx\def\!thx\!ttempa\!thx{\string #1}%
  \!thx\def\!thx\!ttempb\!thx{\csname !tk<\!ttempa>\endcsname}%
  \ifnum \TracingKeys>1
    \!tkReportKey
  \fi
  \!thx\ifx \!ttempb\relax 
    \!thx\!tkReplaceKey
  \else
    \!thx\!ttempb
  \fi}

\def\!tkReportKey{%
  \!taToksA\!thx{\!ttempa}%
  \!thMessage{KEY: "\the\!taToksA"}}

\def\!tkReplaceKey{%
  \!taToksA\!thx{\!ttempa}%
  \!thError {Undefined format key "\the\!taToksA"}
    {\!thReadErrorMsg\!tkUndefined ^^J\!thReadErrorMsg\!tkBadKey}
  \!tkReplaceKeyA}

\def\!tkReplaceKeyA{%
  \!thGetReplacement{\!thReadErrorMsg\!tkReplace}\!tkReplacement
  \!thx\ReadFormatKeys\!tkReplacement}

\!thStoreErrorMsg\!tkUndefined{%
  The format key in " "'s on the next to top line is undefined.}
\!thStoreErrorMsg\!tkBadKey{%
  Type \space E \space to quit now, or
  \space<CR> \space and respond to next prompt.}
\!thStoreErrorMsg\!tkReplace{%
  Type \space<replacement key><CR> \space,
   or simply \space<CR> \space to skip offending key:}

\NewFormatKey b#1{%
  \!thx\!tkJoin\!thx{\the\!taDataColumnTemplate}{#1}%
  \ReadFormatKeys}

\def\!tkJoin#1#2{%
  \!taDataColumnTemplate{#2#1}}%

\NewFormatKey a#1{%
  \!taDataColumnTemplate\!thx{\the\!taDataColumnTemplate #1}%
  \ReadFormatKeys}

\NewFormatKey \{{%
  \!taDataColumnTemplate=\!thx{\!thx{\the\!taDataColumnTemplate}}%
  \ReadFormatKeys}

\NewFormatKey *#1#2{%
  \!taCountA=#1\relax
  \!taToksA={}%
  \!thLoop 
    \ifnum \!taCountA > 0
    \!taToksA\!thx{\the\!taToksA #2}%
    \advance\!taCountA -1
  \repeat 
  \!thx\ReadFormatKeys\the\!taToksA}

\NewFormatKey \LeftGlue#1{%
  \!taLeftGlue{#1}%
  \ReadFormatKeys}

\NewFormatKey \RightGlue#1{%
  \!taRightGlue{#1}%
  \ReadFormatKeys}

\NewFormatKey c{%
  \ReadFormatKeys 
  \LeftGlue\hfil
  \RightGlue\hfil}

\NewFormatKey l{%
  \ReadFormatKeys 
  \LeftGlue{}   
  \RightGlue\hfil}

\NewFormatKey r{%
  \ReadFormatKeys 
  \LeftGlue\hfil
  \RightGlue{}}

\NewFormatKey k{%
  \h@true
  \v@true
  \!tgGetValue{\!tgInsertKern}}

\NewFormatKey i{%
  \h@true
  \v@false
  \!tgGetValue{\!tgInsertKern}}
  
\NewFormatKey j{%
  \h@false
  \v@true
  \!tgGetValue{\!tgInsertKern}}

\NewFormatKey n{%
  \def\!tnStyle{}%
   \futurelet\!tnext\!tnTestForBracket}

\NewFormatKey N{%
  \def\!tnStyle{$}%
   \futurelet\!tnext\!tnTestForBracket}

\NewFormatKey m{%
  \ReadFormatKeys b$ a$}

\NewFormatKey M{%
  \ReadFormatKeys \{ b{$\displaystyle} a$}

\NewFormatKey \m{%
  \ReadFormatKeys l b{{}} m}

\NewFormatKey \M{%
  \ReadFormatKeys l b{{}} M}

\NewFormatKey f#1{%
  \ReadFormatKeys b{#1}}

\NewFormatKey B{%
  \ReadFormatKeys f\bf}

\NewFormatKey I{%
  \ReadFormatKeys f\it}

\NewFormatKey S{%
  \ReadFormatKeys f\sl}

\NewFormatKey R{%
  \ReadFormatKeys f\rm}

\NewFormatKey T{%
  \ReadFormatKeys f\tt}

\NewFormatKey p{%
  \!tgGetValue{\!tgSetUpParBox}}

\NewFormatKey w{%
  \!tkTestForBeginFormat w{\!tgGetValue{\!tfSetWidth}}}

\NewFormatKey s{%
  \!taOnceOnlyTabskipfalse    
  \!tkTestForBeginFormat t{\!tgGetValue{\!tfSetTabskip}}}

\NewFormatKey o{%
  \!taOnceOnlyTabskiptrue
  \!tkTestForBeginFormat o{\!tgGetValue{\!tfSetTabskip}}}

\NewFormatKey |{%
  \!tkTestForBeginFormat |{\!tgGetValue{\!tfSetVrule}}}

\NewFormatKey \|{%
  \!tkTestForBeginFormat \|{\!tfSetAlternateVrule}}

\NewFormatKey .{%
  \!tkTestForBeginFormat.{\!tfFinishFormat}} 

\NewFormatKey \EndFormat{%
  \!tkTestForBeginFormat\EndFormat{\!tfFinishFormat}} 

\NewFormatKey ]{%
  \!tkTestForReFormat ] \!tfEndReFormat}

\def\!tkTestForBeginFormat#1#2{%
  \if!taBeginFormat  
    \def\!ttemp{#2}%
    \!thx \!ttemp    
  \else
    \toks0={#1}%
    \toks2=\!thx{\string\ReFormat}%
    \!thx \!tkImproperUse
  \fi}   

\def\!tkTestForReFormat#1#2{%
  \if!taBeginFormat  
    \toks0={#1}%
    \toks2=\!thx{\string\BeginFormat}%
    \!thx \!tkImproperUse
  \else
    \def\!ttemp{#2}%
    \!thx \!ttemp    
  \fi}   

\def\!tkImproperUse{%
  \!thError{\!thReadErrorMsg\!tkBadUseA "\the\toks0 "}%
    {\!thReadErrorMsg\!tkBadUseB \the\toks2 \space command.
    ^^J\!thReadErrorMsg\!tkBadKey}%
  \!tkReplaceKeyA}
 
\!thStoreErrorMsg\!tkBadUseA{Improper use of key }  
\!thStoreErrorMsg\!tkBadUseB{%
  The key mentioned above can't be used in a }

\def\!tnTestForBracket{%
  \ifx [\!tnext
    \!thx\!tnGetArgument
  \else
    \!thx\!tnGetCode
  \fi}

\def\!tnGetCode#1 {
  \!tnConvertCode #1..!}

\def\!tnConvertCode #1.#2.#3!{%
  \begingroup
    \aftergroup\edef \aftergroup\!ttemp \aftergroup{%
      \aftergroup[%
      \!taCountA #1
      \!thLoop
        \ifnum \!taCountA>0
        \advance\!taCountA -1
        \aftergroup0
      \repeat
      \def\!ttemp{#3}%
      \ifx\!ttemp \empty
      \else
        \aftergroup.
        \!taCountA #2
        \!thLoop 
          \ifnum \!taCountA>0
          \advance\!taCountA -1
          \aftergroup0
        \repeat
      \fi 
      \aftergroup]\aftergroup}%
    \endgroup\relax
    \!thx\!tnGetArgument\!ttemp}
  
\def\!tnGetArgument[#1]{%
  \!tnMakeNumericTemplate\!tnStyle#1..!}

\def\!tnMakeNumericTemplate#1#2.#3.#4!{
  \def\!ttemp{#4}%
  \ifx\!ttemp\empty
    \!taDimenC=0pt
  \else
    \setbox0=\hbox{\m@th #1.#3#1}%
    \!taDimenC=\wd0
  \fi
  \setbox0 =\hbox{\m@th #1#2#1}%
  \!thToksEdef\!taDataColumnTemplate={%
    \noexpand\!tnSetNumericItem
    {\the\wd0 }%
    {\the\!taDimenC}%
    {#1}%
    \the\!taDataColumnTemplate}  
  \ReadFormatKeys}

\def\!tnSetNumericItem #1#2#3#4 {
  \!tnSetNumericItemA {#1}{#2}{#3}#4..!}

\def\!tnSetNumericItemA #1#2#3#4.#5.#6!{%
  \def\!ttemp{#6}%
  \hbox to #1{\hss \m@th #3#4#3}%
  \hbox to #2{%
    \ifx\!ttemp\empty
    \else
       \m@th #3.#5#3%
    \fi
    \hss}}

\def\MakeStrut#1#2{%
  \vrule width0pt height #1 depth #2}

\def\StandardTableStrut{%
  \MakeStrut{\StrutHeightFactor\StrutUnit}
    {\StrutDepthFactor\StrutUnit}}

\def\AugmentedTableStrut#1#2{%
  \dimen@=\StrutHeightFactor\StrutUnit
  \advance\dimen@ #1\StrutUnit
  \dimen@ii=\StrutDepthFactor\StrutUnit
  \advance\dimen@ii #2\StrutUnit
  \MakeStrut{\dimen@}{\dimen@ii}}

\def\Enlarge#1#2{
  \!taDimenA=#1\relax
  \!taDimenB=#2\relax
  \let\!TsSpaceFactor=\empty
  \ifmmode
    \!thx \mathpalette
    \!thx \!TsEnlargeMath
  \else
    \!thx \!TsEnlargeOther
  \fi}

\def\!TsEnlargeOther#1{%
  \ifhmode
    \setbox\z@=\hbox{#1%
      \xdef\!TsSpaceFactor{\spacefactor=\the\spacefactor}}%
  \else
    \setbox\z@=\hbox{#1}%
  \fi
  \!TsFinishEnlarge}
    
\def\!TsEnlargeMath#1#2{%
  \setbox\z@=\hbox{$\m@th#1{#2}$}%
  \!TsFinishEnlarge}

\def\!TsFinishEnlarge{%
  \dimen@=\ht\z@
  \advance \dimen@ \!taDimenA
  \ht\z@=\dimen@
  \dimen@=\dp\z@
  \advance \dimen@ \!taDimenB
  \dp\z@=\dimen@
  \box\z@ \!TsSpaceFactor{}}

\def\OpenUp#1#2{%
  \advance \StrutHeightFactor #1\relax
  \advance \StrutDepthFactor #2\relax}

\def\BeginTable{%
  \futurelet\!tnext\!ttBeginTable}

\def\!ttBeginTable{%
  \ifx [\!tnext
    \def\!tnext{\!ttBeginTableA}%
  \else 
    \def\!tnext{\!ttBeginTableA[c]}%
  \fi
  \!tnext}

\def\!ttBeginTableA[#1]{%
  \if #1u
    \ifmmode                 
      \def\!ttEndTable{
        \relax}
    \else                   
      \bgroup
      \def\!ttEndTable{%
        \egroup}%
    \fi
  \else
    \hbox\bgroup $
    \def\!ttEndTable{%
      \egroup 
      $
      \egroup}
    \if #1t%
      \vtop
    \else
      \if #1b%
        \vbox
      \else
        \vcenter 
      \fi
    \fi
    \bgroup      
  \fi
  \advance\!taRecursionLevel 1 
  \let\!ttRightGlue=\relax  
  \everycr={}
  \ifnum \!taRecursionLevel=1
    \!ttInitializeTable
  \fi}

\bgroup
  \catcode`\|=13
  \catcode`\"=13
  \catcode`\~=13
  \gdef\!ttInitializeTable{%
    \let\!ttTie=~ 
    \let\!ttDH=\- 
    \catcode`\|=\active
    \catcode`\"=\active
    \catcode`\~=\active
    \def |{\unskip\!ttRightGlue&&}
    \def\|{\unskip\!ttRightGlue&\omit\!ttAlternateVrule}%
    \def"{\unskip\!ttRightGlue&\omit&}
    \def~{\kern .5em}
    \def\\{\!ttEndOfRow}%
    \def\-{\!ttShortHrule}%
    \def\={\!ttLongHrule}%
    \def\_{\!ttFullHrule}%
    \def\Left##1{##1\hfill\null}
    \def\Center##1{\hfill ##1\hfill\null}
    \def\Right##1{\hfill##1}%
    \the\EveryTable}
\egroup

\let\!ttRightGlue=\relax  

\def\!ttDoHalign{%
  \baselineskip=0pt \lineskiplimit=0pt \lineskip=0pt %
  \tabskip=0pt
  \halign \the\!taTableSpread \bgroup
   \span\the\!taPreamble
   \ifx \!tfRowOfWidths \empty
   \else 
     \!tfRowOfWidths \cr %
   \fi}

\def\EndTable{%
  \egroup 
  \!ttEndTable}

\def\!ttEndOfRow{%
  \futurelet\!tnext\!ttTestForBlank}

\def\!ttTestForBlank{%
  \ifx \!tnext\!thSpaceToken  
    \!thx\!ttDoStandard
  \else
    \!thx\!ttTestForZero
  \fi}
  
\def\!ttTestForZero{%
  \ifx 0\!tnext
    \!thx \!ttDoZero
  \else
    \!thx \!ttTestForPlus
  \fi}

\def\!ttTestForPlus{%
  \ifx +\!tnext
    \!thx \!ttDoPlus
  \else
    \!thx \!ttDoStandard
  \fi}

\def\!ttDoZero#1{
  \cr} 

\def\!ttDoPlus#1#2#3{
  \AugmentedTableStrut{#2}{#3}%
  \cr} 

\def\!ttDoStandard{%
  \StandardTableStrut
  \cr}

\def\!ttAlternateVrule{%
  \!tgGetValue{\!ttAVTestForCode}}  

\def\!ttAVTestForCode{%
  \ifnum \!tgCode=2              
    \!thx\!ttInsertVrule         
  \else
    \!thx\!ttAVTestForEmpty
  \fi}

\def\!ttAVTestForEmpty{%
  \ifx \!tgValue\empty           
    \!thx\!ttAVTestForBlank
  \else
    \!thx\!ttInsertVrule         
  \fi}

\def\!ttAVTestForBlank{%
  \ifx \!ttemp\!thSpaceToken     
    \!thx\!ttInsertVrule
  \else
    \!thx\!ttAVTestForStar 
  \fi}

\def\!ttAVTestForStar{%
  \ifx *\!ttemp                  
    \!thx\!ttInsertDefaultPR     
  \else
    \!thx\!ttGetPseudoVrule       
  \fi}

\def\!ttInsertVrule{%
  \hfil 
  \vrule \!thWidth
    \ifnum \!tgCode=1
      \ifx \!tgValue\empty 
        \LineThicknessFactor
      \else
        \!tgValue
      \fi
      \LineThicknessUnit
    \else
      \!tgValue
    \fi
  \hfil
  &}

\def\!ttInsertDefaultPR*{%
  \PseudoVrule    
  &}

\def\!ttGetPseudoVrule#1{%
  \toks0={#1}%
  #1&}

\def\PseudoVrule{}

\def\ifundefined#1{\expandafter\ifx\csname#1\endcsname\relax}
\def\!ttuse#1{%
  \ifnum #1>\@ne 
    \omit 
    \ifundefined{mscount}
       \@multicnt=#1
       \advance\@multicnt by \m@ne
       \advance\@multicnt by \@multicnt
       \!thLoop 
         \ifnum\@multicnt>\@ne 
         \sp@n %
       \repeat 
    \else
       \mscount=#1      
       \advance\mscount by \m@ne
       \advance\mscount by \mscount
       \!thLoop 
         \ifnum\mscount>\@ne 
         \sp@n %
       \repeat 
    \fi
    \span 
  \fi}

\def\!ttUse#1[{%
  \!ttuse{#1}%
  \ReFormat[}

\def\!ttFullHrule{%
  \noalign
  \bgroup
  \!tgGetValue{\!ttFullHruleA}}

\def\!ttFullHruleA{%
  \!ttGetHalfRuleThickness 
  \hrule \!thHeight \dimen0 \!thDepth \dimen0
  \penalty0 
  \egroup} 

\def\!ttShortHrule{%
  \omit
  \!tgGetValue{\!ttShortHruleA}}

\def\!ttShortHruleA{%
  \!ttGetHalfRuleThickness 
  \leaders \hrule \!thHeight \dimen0 \!thDepth \dimen0 \hfill
  \null    
  \ignorespaces} 

\def\!ttLongHrule{%
  \omit\span\omit\span \!ttShortHrule}

\def\!ttGetHalfRuleThickness{%
  \dimen0 =
    \ifnum \!tgCode=1
      \ifx \!tgValue\empty
        \LineThicknessFactor
      \else
        \!tgValue    
      \fi
      \LineThicknessUnit
    \else
      \!tgValue      
    \fi
  \divide\dimen0 2 }

\def\WidenTableBy#1{%
  \ifdim #1=0pt
    \!taTableSpread={}%
  \else
    \!taTableSpread={spread #1}%
  \fi}

\def\JustLeft{%
  \omit \let\!ttRightGlue=\hfill}
\def\JustCenter{%
  \omit \hfill\null \let\!ttRightGlue=\hfill}

\let\\=\!tacr
\catcode`\!=12
\catcode`\@=12

\def\strich{\vskip0.5cm\hrule\vskip3ptplus12pt\null}
\newenvironment{Defliste}[1]%
{\begin{list}{}{ \settowidth{\labelwidth}{\it{#1}\quad}
                      \setlength{\leftmargin}{\labelwidth}}}
{\end{list}}

 \includeonly{Zwi92_Intro,Zwi92_Represent_sans_N,Zwi92_Algo_sans_N,Zwi92_Results_sans_N,Zwi92_fonctionnel,Zwi92_Conclusion_sans_N,Annexe1,Annexe2,Biblio}

\usepackage[latin1]{inputenc}

\usepackage{graphicx}
\graphicspath{%
    {converted_graphics/}
    {/}
}
\begin{document}

\title{R\'esolution num\'erique\\[2mm] du probl\`eme de Dirichlet $\Delta u = a\,u^3$\\[2mm] \`a l'aide du mouvement brownien}        
\author{Jean-Paul MORILLON\\[4mm]        
Rapport de recherche \\[2mm]
Laboratoire PIMENT \\[1mm]
Universit\'e de La R\'eunion\\[1mm]
97487 Saint-Denis C\'edex\\[2mm]
La R\'eunion}
\date{\today}          
\maketitle

\abstract{Des représentations stochastiques de solutions de problèmes de Dirichlet déterministes linéaires et non linéaires sont déduites de l'application de la formule de Itô. Ces représentations sont utilisées pour établir des algorithmes de calcul par simulation de marches aléatoires. Les méthodes numériques associées sont appliquées à des exemples de problèmes linéaires et non linéaires. Les résultats des essais avec une fonction source, des estimations de temps d'arrêt et des courbes de régression quadratique sont présentés. Des solutions $u$ du problème de Dirichlet $\Delta u = a\,u^3$ sont calculées par cette méthode purement stochastique pour des valeurs de $a$ positives et négatives.

\noindent{\bf Mots clés\,:} Problèmes de Dirichlet déterministes linéaires et non linéaires --- \'Equations différentielles stochastiques --- Fonctionnelles intégrales --- Résolu\-tion numérique stochastique --- Mar\-che aléatoire
}

\strich

\selectlanguage{english}
\abstract{In this paper, we are interested in numerical solution of some linear boundary value problems with Dirichlet boundary part, by the means of simulation of random walks. We use a probabilistic interpretation of solution $u$, assuming that the coefficient and the boundary data are sufficiently smooth, and applying It\^o's formula. From these stochastic representations of solution, we extend some algorithms obtained for standard boundary conditions to the quasi-linear source of the type $f(u)= a\,u^3$. For positive and negative parameter $a$, we then obtain numerical results by applying the stochastic methods based upon these generalized algorithms.}

\noindent{\bf Keywords:} Linear and nonlinear Dirichlet BVP; SDE; Probabilistic representation; Numerical stochastic method; Random walks

\strich

\clearpage
\selectlanguage{francais}

\tableofcontents
\setcounter{tocdepth}{2}

\chapter{Introduction}
\label{Introduction}

Le calcul des solutions d'une équation différentielle du type $u"=f(u(t))$ est d'un intérêt constant (voir le récent article \cite{ambrosio}). Notre objectif est principalement de calculer des solutions de problèmes de Dirichlet $u"=a u^3$ mis en exergue dans \cite{Zwi98} où $a=1$.

La résolution numérique de problèmes aux limites peut être faite, de manière classique, en mettant en \oe uvre des méthodes de différences finies ou d'élé\-ments finis, associées à un maillage de l'espace.

Des méthodes de Monte-Carlo peuvent également être appliquées; elles conduisent, après discrétisation des équations, à des traitements particuliers liés à la géométrie locale des domaines (voir, par exemple, \cite{haji}, \cite{dautray}, \cite{marshall} et \cite{kushner:92}).

Cette étude explore une autre voie en utilisant la théorie des processus de Markov qui fournit des représentations intégrales pour les solutions des problèmes aux limites linéaires et non linéaires déterministes stationnaires. On sait, en particulier, que le problème de Dirichlet admet une représentation intégrale \cite{freidlin} utilisée pour établir des algorithmes de calcul de la solution par simulation de marches aléatoires \cite{souza}.

Les méthodes numériques qui en découlent ne nécessitent pas d'entrer en mé\-moi\-re un maillage de discrétisation dans le cas linéaire. Dans le cas non linéaire, on utilise le maillage élémentaire support des marches aléatoires. Ceci conduit à une programmation courte et facile à vérifier pas à pas.

De manière plus précise, on résout un problème de Dirichlet déterministe et stationnaire dans des cas linéaires et non linéaires. On montre notamment, à partir d'exemples, que l'algorithme stochastique converge plus rapidement par relaxation dans le cas non linéaire.

Le chapitre 2 est consacré aux représentations stochastiques des solutions de problèmes de Dirichlet linéaires et non linéaires. Les approches des représentations par réalisations des processus sont présentées dans le chapitre 3 et les algorithmes de calcul correspondants sont établis. Des simulations numériques des marches aléatoires sont effectuées dans le chapitre 4 pour calculer des valeurs des solutions. Une approche fonctionnelle est présentée dans le chapitre 5.

\chapter{Représentations des solutions}
\label{Representation}
Ce chapitre est consacré à la représentation des solutions de problèmes déterministes de Dirichlet linéaires ou non linéaires par des intégrales stochastiques. 

Considérons un ouvert borné $G$ de ${\mathbb{R}}^d$ (la dimension $d$ étant fixée). On note $x \in \mathbb{R}^d$ la variable d'espace et $n$ la normale unitaire intérieure définie sur la frontière $\partial G$.

\`A partir d'équations différentielles stochastiques (EDS) (pour ces der\-nières, voir par exemple, \cite{pardoux,kloeden:92} et les références associées), nous pouvons représenter les solutions de problèmes aux limites en faisant intervenir l'espé\-rance d'inté\-grales\,: fonctionnelles de trajectoires et de fonctions aléatoires solutions des EDS. L'application de la formule de Itô en liaison avec les EDS conduit à la représentation des solutions.

Sur le plan fonctionnel, les représentations que nous obtenons n'ont été établies en général que sous des hypothèses de régularité sur le domaine $G$ et sur les données aux limites (se reporter aux références mentionnées dans la suite de ce chapitre). Néanmoins, nous utilisons ici les représentations stochastiques lorsque les données ne sont pas nécessairement régulières\,: la frontière $\partial G$ du domaine peut, par exemple, présenter des coins ou des arêtes --- cette procédure s'est déjà montrée efficace par exemple dans \cite{haji,dautray,kushner:92} pour les équations discrétisées et dans \cite{souza} pour le problème de Dirichlet.

Dans le paragraphe 1, nous représentons la solution du problème de Dirichlet. Le paragraphe 2 est consacré à des problèmes de Dirichlet quasi-linéaires où le second membre dépend de la solution.
\section{Problème de Dirichlet}
Considérons le problème de Dirichlet\,:
\begin{eqnarray}
\left\{
       \begin{array}{rcll}
          -\frac12 \, \Delta u & = & f   & \mbox{dans } G      \\
                             u & = & g & \mbox{sur } \partial G      \end{array}
\right.  & \mbox{} & \mbox{}
\label{Dirichlet}
\end{eqnarray}
où $u$ est la fonction inconnue, définie sur $G$, à valeurs réelles 
et les données sont les fonctions scalaires $f$ et $g$, définies respectivement sur $G$ et $\partial G$.

Introduisons le processus de Markov\,:
$
X_t^x = x + W_t\,,~ t \geq 0\,,
$
où $W_t$ est un processus de Wiener standard à valeurs dans $\mathbb{R}^d$ tel que l'on ait\,: $X_0^x = x$.

Soit $\overline{G}$ le complémentaire de $G$ dans $\mathbb{R}^d$. Le temps d'atteinte de $\overline{G}$ est défini par\,:
$$
\tau^x = \inf \left\{\, t > 0 \mid X_t^x \in \overline{G} \,\right\}.
$$

Lorsque la condition\,: $E[\tau^x]<+\infty$ est vérifiée en tout point $x$ de $G$, la formule de Itô appliquée à $\,u \left( \, X_t^x \, \right)$ conduit à la représentation de la solution de (\ref{Dirichlet}), que nous écrivons, pour des raisons techniques, sous la forme\,:
\begin{equation}
u(x) = E \left[\, Y \,\right], \quad x \in G,
\label{first representation for u}
\end{equation}
avec la variable aléatoire\,:
\begin{equation}
Y = \int_0^{\tau^x} f\left( \, X_t^x \, \right) \, dt + g \left( \, X_{\tau^x}^x \, \right) .
\label{first representation for Y}
\end{equation}
Des résultats fonctionnels sur ce problème et la représentation associée  peuvent, par exemple, être trouvés dans \cite{bensoussan}, \cite{freidlin} et \cite{karatzas}.

Considérons le temps d'atteinte $\tau$ et montrons que $\tau$ dépend du coefficient de diffusion $a$.

\noindent Soit $B(O,r)$ une boule ouverte de $\mathbb{R}^d$, de centre $O$ et de rayon $r$, de frontière ${\cal C}(O,r)$.\\
Considérons $u\,: B(O,r) \longrightarrow \mathbb{R}$ la solution unique du problème\,:
\begin{eqnarray}
\left\{
       \begin{array}{rcll}
 -\frac12 \, a \,\Delta u & = & 1 & \mbox{dans } B(O,r)       \\[2mm]
                        u & = & 0 & \mbox{sur  } {\cal C}(O,r)
       \end{array}
\right.  & \mbox{} & \mbox{}
\label{a tend vers zero}
\end{eqnarray}
où $a > 0$ est supposé constant.

Alors $u$ admet la représentation\index{représentation} (\ref{first representation for u}--\ref{first representation for Y}) avec le processus\index{processus stochas}\,:
\begin{equation}
X_t^x = x + \sqrt{a}\,W_t\,.
\label{processus_diffusion}
\end{equation}

D'après le raisonnement appliqué par \cite[page 253, avec $a \equiv 1$]{karatzas}, on en déduit que le temps d'atteinte\index{tps:atteinte} $\tau^x$ du complémentaire de $B(O,r)$ par le processus $X_t^x$, à partir d'un point intérieur $x$, vérifie\,:
$$
u(x) = E\left[\, \tau^x \, \right] = \frac{r^2 - |x|^2}{d\,a} \,, \quad x \in B(O,r).
$$
Par conséquent, lorsque $a$ tend vers zéro par valeurs positives, $E\left[ \,\tau^x \, \right]$ tend vers $+\infty$.

De plus, lorsque $a$ décroît, les accroissements\index{accroissements du processus} du processus (\ref{processus_diffusion}) décrois\-sent avec $\sqrt{a}$ car $dX_t^x=\sqrt{a}~dW_t$\,; ceci limite notre champ d'investigation lors des essais numériques. C'est pourquoi, avant toute résolution numérique, on s'intéresse au comportement de la solution $u$ lorsque $a$ tend vers zéro sur $G$ pour éviter les éventuelles difficultés numériques.

\section{Problème de Dirichlet quasi-linéaire}

\indent Considérons le problème\index{pb:Dir:non lin}\,:
\begin{eqnarray}
\left\{
       \begin{array}{rcll}
-\frac12 \, \Delta u  & = & f(x,u)   & \mbox{dans } G       \\[2mm]
                   u  & = & g & \mbox{sur } \partial G       \end{array}
\right.  & \mbox{} & \mbox{}
\label{Dirichlet-non lineaire}
\end{eqnarray}
où $f$ est une fonction donnée de $(x,u)$, $f$ non linéaire en $u$, supposée régulière.

Appliquons la méthode des approximations successives\index{méthode:approx:successives}\,:\\
{\it \`A l'ordre $m$}, substituons $f$ par $f\left(\,u^{(m-1)}\,\right)$ dans (\ref{Dirichlet}),  alors (\ref{first representation for u}) et (\ref{first representation for Y}) donnent la représentation\index{représentation} de la solution de (\ref{Dirichlet-non lineaire}) sous la forme itérative\index{forme itérative}\,:
\begin{equation}
u^{(m)}(x)  =  E \left[\, Y^{(m)} \,\right], \quad x \in G,
\label{representation potentiel non lineaire de u}
\end{equation}
avec la variable aléatoire\index{v.a.}\,:
\begin{eqnarray}
Y^{(m)} & = & \int_0^{\tau} f\biggl( \, X_t^x ,\, u^{(m-1)}\left(\, X_t^x \,\right) \, \biggr) \, dt + g \left( \, X_{\tau}^x \, \right).
\label{representation potentiel non lineaire de Y}
\end{eqnarray}
Sous certaines conditions, en particulier sur la régularité du second membre $f$, la suite $Y^{(m)}$ tend vers la représentation stochastique de la solution $u$ du problème \ref{Dirichlet-non lineaire} (voir le chapitre \ref{proba:fonctionnel}).

L'objectif principal de ce travail est le calcul de solutions, supposées suffi\-samment régulières et bornées, du système non linéaire :
\begin{eqnarray}
\left\{
       \begin{array}{rcll}
\Delta u  & = & a\, u^3   & \mbox{dans } G=]0;1[       \\[2mm]
       u(0) & = & 0 &  \\[2mm]
       u(1) & = & 1 &  \end{array}
\right.  & \mbox{} & \mbox{}
\label{Sys_u_trois}
\end{eqnarray}
avec le paramètre $a\in \mathbb{R}$ et où le second membre $f$ est une fonction non linéaire de $u$.

\section{Commentaires}
Les représentations des solutions ont été établies à partir d'équations différentielles stochastiques (EDS\index{EDS}) (pour ces dernières, voir par exemple, \cite{pardoux,kloeden:92} et les références associées).

Sur le plan fonctionnel, les représentations que nous avons obtenues, n'ont été établies en général que sous des hypothèses de régularité\index{hypo:régularité} sur le domaine $G$ et sur les données aux limites. Néanmoins, nous utilisons les représentations stochastiques\index{repr:stochas} lorsque les données ne sont pas nécessairement régulières\,: la frontière $\partial G$ du domaine peut, par exemple, présenter des coins\index{coins:aretes} ou des arêtes --- cette procédure s'est déjà montrée efficace par exemple dans \cite{haji,dautray,kushner:92} pour les équations discrétisées\index{equations discretisees} et dans \cite{souza} pour le problème de Dirichlet\index{pb de Dir sans pot} sans terme de potentiel (pour cette notion, voir \cite{morillon:95}).

Les représentations associées aux problèmes de Dirichlet (\ref{Dirichlet}) peuvent être trouvées, par exemple, dans \cite{bensoussan}, \cite[page 127]{freidlin} et \cite[pages 244 et 364]{karatzas}, où la solution est supposée de classe ${\cal C}^2\left(\,G\,\right)$. Des résultats fonctionnels concernant la régularité\index{régularite de la solution} de la solution des problèmes de Dirichlet (\ref{Dirichlet}), et la régularité de ses dérivées, sont rappelés dans \cite{dautray} où les conditions sur les données sont telles que la solution est de classe ${\cal C}^2\left(\,G\,\right)$.

Dans le cas non linéaire, des représentations itératives du type (\ref{representation potentiel non lineaire de u}--\ref{representation potentiel non lineaire de Y}) peuvent être trouvées dans \cite{freidlin}, mais seulement dans le cas d'une équation parabolique semi-linéaire\index{eq:parabol:semilin}.

\chapter{Algorithmes de r\'esolution}
\label{Algorithmes}
Ce chapitre a pour objet de donner des réalisations des processus associés aux représentations du chapitre précédent. Ces réalisations permettent de calculer effectivement les solutions des problèmes aux limites. Le traitement des conditions aux limites se ramène à un comptage des nombres d'absorptions sur les frontières du domaine. Le traitement de la non-linéarité conduit à des itérations des algorithmes de calcul des solutions de problèmes linéaires.
\section{Probl\`eme de Dirichlet}
\label{Dirichlet_lineaire}
La solution du problème (\ref{Dirichlet}) est représentée par la moyenne d'une varia\-ble aléatoire $Y$ fonction d'un processus de Wiener standard \cite{souza}.\par
La représentation de ce problème par les équations (\ref{first representation for u}--\ref{first representation for Y}), page \pageref{first representation for Y}, montre qu'il suffit de calculer $NT$ valeurs approchées de $Y$, notées $Y_1,\ldots,Y_{NT}$, pour obtenir $u(x)$\,:
$$u(x)=\frac{1}{NT} \sum_{n=1}^{NT} Y_n \,.$$
Chaque valeur $Y_n$ est obtenue à partir d'une simulation d'un processus de Wiener de la manière suivante\,: le pas $h>0$ étant fixé, on simule le processus $X_t^x$ par la suite de vecteurs aléatoires $X_0,\ldots,X_k,\ldots$ définis par les formules de récurrence\,:
\begin{eqnarray*}
X_0&=&x\: \in \: G\,,   \qquad \mbox{initialisation,}\\
X_{k+1}&=&X_k+h\,D_k\,, \qquad k\: \in \: \mathbb{N},
\end{eqnarray*}
où $D_k$ est un vecteur aléatoire tel que, si $e_i$, $i=1,\ldots,d$, est la base cano\-nique de $\mathbb{R}^d$, alors 
$$D_k=\sum_{i=1}^d D_{k,i} \, e_i$$
vérifie\,:
$$\mbox{Prob}(D_k=e_i)=\mbox{Prob}(D_k=-e_i)=\frac{1}{2d},\qquad i=1,\ldots,d\;.$$
\indent En dimension trois, par exemple, on engendre une suite, nécessairement finie, de nombres\,: $$U_0,\ldots,U_k,\ldots,U_N\,,$$ simulant une variable aléatoire $U$ uniformément distribuée sur l'intervalle $[0,1]$ à l'aide d'un générateur de nombres pseudo-aléatoires\,; puis on pose, pour tout $k$\,:
$$\frac{i}{6}\, \leq \, U_k \, \leq \, \frac{i+1}{6}~~\Longrightarrow~~D_k=(-1)^i e_j\,,\quad i=0,\ldots,5,$$
avec
$$j=1+\left[ \, \frac{i}{2} \, \right] $$
où les crochets désignent la partie entière de l'argument.\par
\begin{figure}[h]
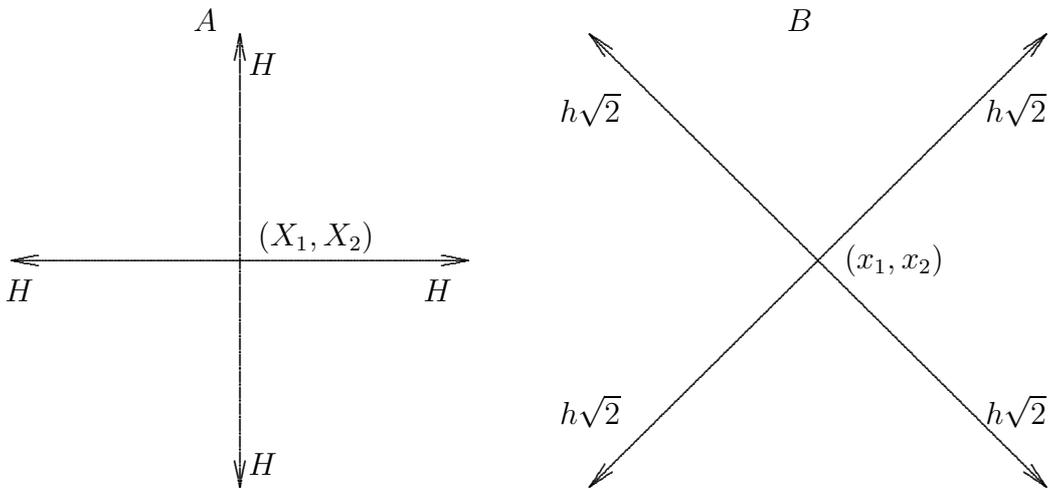

\begin{minipage}[c]{.46\linewidth}
\begin{center} $A\qquad $
\beginpicture
\setcoordinatesystem units <1cm,1cm>
\setplotarea x from 0 to 6 , y from 0 to 6
\arrow <10pt> [.2, .4]  from 3 3 to 6 3 
\arrow <10pt> [.2, .4]  from 3 3 to 0 3 
\arrow <10pt> [.2, .4]  from 3 3 to 3 6 
\arrow <10pt> [.2, .4]  from 3 3 to 3 0 
\put {$H$} at 3.3 5.6
\put {$H$} at 3.3 0.3
\put {$H$} at 5.6 2.6
\put {$H$} at 0.1 2.6
\put {$(X_1,X_2)$} at 4 3.3
\endpicture
\end{center}
\end{minipage}\hfill
\begin{minipage}[c]{.46\linewidth}
\begin{center} $B $
\beginpicture
\setcoordinatesystem units <1cm,1cm>
\setplotarea x from 0 to 6 , y from 0 to 6
\arrow <10pt> [.2, .4]  from 3 3 to 6 6
\arrow <10pt> [.2, .4]  from 3 3 to 6 0 
\arrow <10pt> [.2, .4]  from 3 3 to 0 6 
\arrow <10pt> [.2, .4]  from 3 3 to 0 0 
\put {$h\sqrt{2}$} at 5.6 5
\put {$h\sqrt{2}$} at 5.6 1
\put {$h\sqrt{2}$} at 0   5
\put {$h\sqrt{2}$} at 0   1
\put {$(x_1,x_2)$} at 4 3
\endpicture
\end{center}
\end{minipage}
\caption{Pas élémentaires d'une marche aléatoire dans le plan, avec une probabilité de 1/4 chacun}
\label{pas elementaire}
\end{figure}

En dimension deux, à chaque intervalle de temps, le déplacement se fait :  
\begin{enumerate}
  \item soit à pas constant $h$ dans l'une des deux directions d'un repère orthonormé, direction choisie au hasard avec une probabilité 1/2, et avec un des deux sens choisi au hasard avec une probabilité 1/2, comme le montrent la figure \ref{pas elementaire}-$A$ et l'algorithme~1 du pas élémentaire de la figure~\ref{algorithme du pas elementaire}, où la fonction \verb!random! renvoie un réel appartenant à $[0,1[$;
  \item soit à pas constant $h\sqrt{2}$ dans l'une des deux bissectrices d'un repère orthonormé, direction choisie au hasard avec une probabilité 1/2, et avec un des deux sens choisi au hasard avec une probabilité 1/2, comme le montrent la figure \ref{pas elementaire}-$B$ et l'algorithme~2 du pas élémentaire de la figure~\ref{algorithme du pas elementaire}, où la fonction \verb!randi! renvoie un entier pair ou impair pour chaque coordonnée.
\end{enumerate}

\mbox{ }
\\
\begin{figure}[h]
\begin{minipage}[c]{.46\linewidth}
{Algorithme 1 du pas élémentaire}
\begin{verbatim}
U1:=random;     U2:=random;
if U1<0.5 then
 if U2<0.5 then X2:=X2 + H
           else X2:=X2 - H
 else
 if U2<0.5 then X1:=X1 + H
           else X1:=X1 - H;
\end{verbatim}\end{minipage}\hfill
\begin{minipage}[c]{.46\linewidth}
{Algorithme 2 utilisé sous MATLAB}
\begin{verbatim}



m =[x1;x2] ;
m =m+h.*(-1).^randi(2,[2, 1]);


\end{verbatim}\end{minipage}
\caption{Algorithmes du pas élémentaire}
\label{algorithme du pas elementaire}
\end{figure}

Chaque simulation donne une réalisation $X_k$, $0 \leq k \leq N$, d'une marche aléatoire dans $G$ telle que\,:
$$X_0=x\: \in \: G\,, \ldots, X_{N-1} \: \in \: G~~~~~\mbox{et}~~~~~X_N \: \not \in \: G\,,$$
dont la trajectoire est la courbe polygonale reliant les points $X_k$ et $X_{k+1}$, $k=0,\ldots,N-1\,,$ successivement.\par
Le processus aléatoire $X_t^x$ est alors discrétisé selon la méthode d'Euler\,; en particulier, si $h$ est le pas d'une marche aléatoire simulée, alors on a les accroissements en pas et en temps suivants\,:
$$\delta X_t^x=h~~~~~\mbox{et}~~~~~\delta t = h^2/d \,.$$

\`A chaque marche aléatoire simulée $X_k$, on définit et on calcule l'encaissement $Y_n$\,:
$$Y_n=\delta t \cdot \sum_{k=1}^{N} f(X_k) + g_1(X_N)\,.$$
Puis on répète cette démarche $NT$ fois pour obtenir un échantillon de taille $NT$ pour la variable aléatoire $Y$ et on en déduit $u(x)$.\par
L'entier $NT$ correspond au nombre de marches aléatoires simulées, donc au nombre d'absorptions par $\overline{G}$. Cet entier $NT$ est l'indice maximum de la boucle principale du programme. La deuxième et dernière boucle du programme suit pas à pas chaque marche aléatoire en additionnant au fur et à mesure la valeur de la fonction source $f$, puis, pour finir, l'effet de l'absorption avant de stocker l'encaissement $Y_n$ correspondant.\par
Considérons un domaine $G$  de $\mathbb{R}^2$ et un point \verb!(X1D,X2D)! donné dans $G$. Alors un algorithme possible $A1$ du problème de Dirichlet se présente sous la forme\,:
\begin{verbatim}
    YN := 0.0;
    for i:=1 to NT do
    begin
	     X1:= X1D;     X2:= X2D;     YA:= 0.0;
	     while INTER(X1,X2) = true  do
	     begin
	          YA:= YA + F(X1,X2);
	          X1:= X1 + pas élémentaire aléatoire;
	          X2:= X2 + pas élémentaire aléatoire;
	     end;
	     YN:= YN + H * H * YA/2 + G1(X1,X2);
    end;
    U:= YN / NT;
\end{verbatim}
\begin{figure}[h]
\caption{Boucle du problème de Dirichlet}
\label{boucle Dirichlet}
\end{figure}
La fonction \verb!INTER! prend la valeur {\em vraie\/} lorsque le point \verb!(X1,X2)! est dans $G$ et {\em fausse\/} sinon. La fonction \verb!F! définit l'effet de la source en chaque point de la marche dans $G$ et \verb!G1! définit l'effet d'absorption sur la frontière du domaine.
%
\section{Probl\`eme de Dirichlet quasi-linéaire}
La solution du problème (\ref{Dirichlet-non lineaire}) est représentée par la moyenne d'une varia\-ble aléatoire $Y$ fonction d'un processus de Wiener standard.\par
\noindent La représentation de ce problème par les équations (\ref{representation potentiel non lineaire de u}--\ref{representation potentiel non lineaire de Y}), page \pageref{representation potentiel non lineaire de u}, montre qu'il suffit de calculer $NT$ valeurs approchées de $Y$, notées $Y_1,\ldots,Y_{NT}$, pour obtenir $u(x)$\,:
$$u(x)=\frac{1}{NT} \sum_{n=1}^{NT} Y_n \,.$$
Chaque valeur $Y_n$ est obtenue à partir d'une simulation d'une marche aléatoire selon l'algorithme \ref{algorithme du pas elementaire} comme pour la résolution du problème (\ref{Dirichlet}). 

Comme indiqué dans le paragraphe \ref{Dirichlet_lineaire}, chaque simulation donne une réalisation $X_k$, $0 \leq k \leq N$, d'une marche aléatoire dans $G$, le processus aléatoire $X_t^x$ étant discrétisé selon la méthode d'Euler.

\`A chaque marche aléatoire simulée $X_k$, on définit et on calcule l'encaissement $Y_n$\,:
$$Y_n=\delta t \cdot \sum_{k=1}^{N} f\left(\,u\left(X_k\right)\,\right) + g(X_N)\,.$$

Puis on répète cette démarche $NT$ fois pour obtenir un échantillon de taille $NT$ pour la variable aléatoire $Y$ et on en déduit $u(x)$. Dans le cas non linéaire, le calcul est effectué aux points de discrétisation du domaine (points supports des marches aléatoires), puis réitéré à partir des valeurs obtenues en ces points.

Le cas non linéaire demande quelques explications. Non seulement la solution est calculée sur l'ensemble des points de discrétisation du domaine, mais encore la non-linéarité impose des itérations\index{itérations} des algorithmes précédents.

Considérons un domaine $G=\left]0,L\right[$. Alors un algorithme possible $A2$ du problème de Dirichlet non linéaire\index{pb:Dir:non:lin} (\ref{Dirichlet-non lineaire}) se présente sous la forme\,:
\begin{verbatim}
% Problème non linéaire avec condition de Dirichlet
% Subdivision en maxpt intervalles du domaine G = ] 0, L [
maxpt1 = maxpt + 1 ;
h = L / maxpt ; % Pas de la marche aléatoire
% miter % Nombre d'itérations
% nt    % Nombre de marches aléatoires par point
% U     % Initialisation de U
% Itération du problème non linéaire
for iter = 1:miter
    for J = 2:maxpt
        ys  = 0.0 ; dir = 0.0 ;
        for I = 1 : nt % Nb de marches aléatoires par point
            m = J;  ya = 0.0;
            while ( Inter(m)==1 )
                ya = ya + F( U( m ) ) ;
                m  = m + (-1)^randi([0, 1]) ;
            end
            dir = dir + G( m ) ;
            ys  = ys  + ya ;
        end
        U(J) = ( h*h*ys + dir ) /nt ;
    end
end
% Affichage des valeurs obtenues
\end{verbatim}

La variable \,\verb!maxpt! est le nombre d'intervalles de la subdivision\index{subdivision} de $G$\,, \verb!iter! le nombre d'itérations, \verb!miter! le maximum d'itérations\index{itérations}, \verb!U! le tableau d'argument\/ \,\verb!J! des valeurs successives de \verb!U! obtenues à chaque itération.

La marche\index{m.a.} aléatoire se fait d'un point $x_m$ à un point voisin $x_{m-1}$ ou $x_{m+1}$ de la subdivision du domaine\,; le pas aléatoire\index{pas:a.:élément} élémentaire est donc discret\,: on passe d'un indice \verb!M! à un indice voisin \verb!M-1! ou \verb!M+1!, compris entre $1$ et \verb!maxpt!$\mbox{}+1$.
L'entier \verb!nt! correspond au nombre de marches aléatoires simulées à chaque itération, donc au nombre d'absorptions par $\overline{G}$ à chaque itération.

En chaque point $x_j$ de la subdivision, une valeur de la solution $u(x_j)$ est estimée, intégrée à la boucle de subdivision, donc prise en compte pour l'estimation des valeurs de la solution aux points suivants $x_k$, $k>j$. Plus précisément, à l'étape $m$, on calcule une valeur de $u_{j+1}^m$ à l'aide des valeurs de l'étape $m$ en cours : $u_1^m, \ldots, u_j^m$, et des valeurs de l'étape précédente $m-1$ : $u_{j+1}^{m-1},\ldots,u_{\mbox{\tiny maxpt}}^{m-1}$.

Dans le paragraphe \ref{Exemple de probleme de Dirichlet non lineaire} du chapitre suivant, nous nous intéresserons plus particulièrement au système (\ref{Sys_u_trois}) et, compte tenu des essais numériques et de la convergence numérique, nous serons amenés à modifier cet algorithme.

\section{Commentaires}

La vitesse de convergence des algorithmes dépend du pas $h$ de la marche aléatoire\,; le processus brownien $X$ discrétisé et le temps d'atteinte $\tau$ de la frontière $\partial G$ en dépendent également. Lors des essais numériques, notre champ d'investigation est limité à des valeurs de $h$ qui permettent au processus d'atteindre la frontière en des temps raisonnables, donc aux vitesses de calcul des ordinateurs utilisés ; pour éviter d'éventuelles difficultés numé\-riques, on s'intéresse au comportement des temps d'atteinte avant toute résolution numérique.

\chapter{Essais numériques}
\label{Essais}

Des simulations numériques de la méthode stochastique introduite précé\-dem\-ment sont présentées dans ce chapitre. Ces essais sont destinés à valider les représentations du chapitre 2 ainsi que les méthodes ap\-prochées du cha\-pitre 3. Enfin, le système de Dirichlet non linéaire (\ref{Sys_u_trois}) est résolu à l'aide du mouvement brownien approchée par des marches aléatoires.

%
\section{Problème de Dirichlet}
\label{Exemple de probleme de Dirichlet}

\begin{figure}
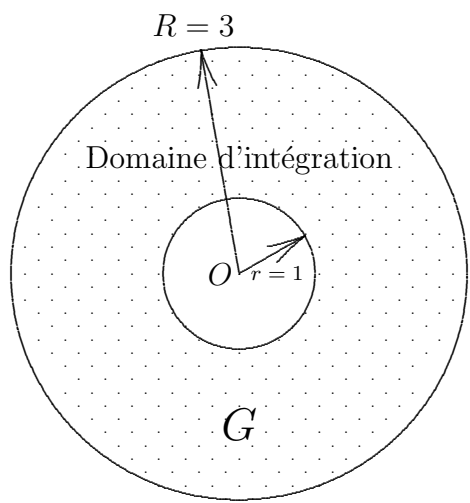

\beginpicture
\setcoordinatesystem units <1cm,1cm>  point at 0 0
\setplotarea x from -7 to 10 , y from -3 to 3
\circulararc 360 degrees from 1 0 center at 0 0
\circulararc 360 degrees from 3 0 center at 0 0
\put {$O$} at -0.25 0
\put{Domaine d'intégration} at 0 1.5
\arrow <.5cm> [.2,.5] from 0 0 to 0.866   .5
 \put {\scriptsize $r=1$} at 0.5 0.0
\arrow <.5cm> [.2,.5] from 0 0 to -0.5   2.958
 \put {$R=3$} at -0.6 3.3
\put {\Large $G$} at 0 -2
  \setquadratic
  \vshade  -3 0 0 <,,z,> -1.5 0 2.598 -1 0 2.828 -0.5 0.866 2.958  0.0 1.0 3.0 <z,,,> 0.5 0.866 2.958 1.0 0.0 2.828 <,,z,> 1.5 0 2.598  3 0 0 <,,z,z> 1.5 0 -2.598 1 0 -2.828 0.5 -0.866 -2.958 0 -1 -3  /
\startrotation by -1 0 about 0 0
  \setquadratic
  \vshade  -3 0 0 <,,,z> -1.5 0 2.598  -1 0 2.828 -0.5 0.866 2.958 0.0 1.0 3.0 <z,,,> 0.5 0.866 2.958 1.0 0.0 2.828 1.5  0 2.598 3 0 0 1.5 0 -2.598 1 0 -2.828 0.5 -0.866 -2.958 0 -1 -3  /
\stoprotation
\endpicture
\caption{Couronne circulaire $G$ de frontière ${\cal C}(O;r=1) \cup {\cal C}(O;R=3)$}
\label{Ring}
\end{figure}

Considérons le problème (\ref{Dirichlet}) dans la couronne $G$ définie par\,:
$$
G=\left\{ (x,y) \in \mathbb{R}^2 \mid 1< \sqrt{x^2+y^2} < 3\right\}
$$
(cf. figure~\ref{Ring}), 
avec les données suivantes\,:
\begin{eqnarray*}
f   & \equiv & 0 \qquad \mbox{dans } G\\\
g & \equiv & 4 \qquad \mbox{sur  } {\cal C}(0;1) = \left\{ (x,y) \in \mathbb{R}^2 \mid x^2+y^2 = 1 \right\}\\
g & \equiv & 6 \qquad \mbox{sur  } {\cal C}(0;3) = \left\{ (x,y) \in \mathbb{R}^2 \mid x^2+y^2 = 9 \right\}.
\end{eqnarray*}
Ce problème de Dirichlet admet pour solution exacte\,:
\begin{eqnarray}
u(x,y) & = & 4 + 2 ~ {\displaystyle \frac{\ln \sqrt{x^2+y^2}}{\ln 3}}\;.
\label{solution du probleme de Dirichlet}
\end{eqnarray}
\newpage
\`A l'aide de l'algorithme $A1$, on a obtenu, dans un premier temps, les résultats du tableau~\ref{Tableau I}, page \pageref{Tableau I}, avec les données suivantes\,: \\[3mm]
\begin{tabular}{lrcl}
Point considéré\,:  & $  (x,y)    $ & $=$ & $ ( 2 ; 0 ) $\\
Nombre de tirages\,:& $   NT      $ & $=$ & $ 10^5   $\\
Valeur exacte\,:    & $u( 2 ; 0 ) $ & $\approx$ & $ 5,\!26186.$
\end{tabular}

\begin{figure}[h]
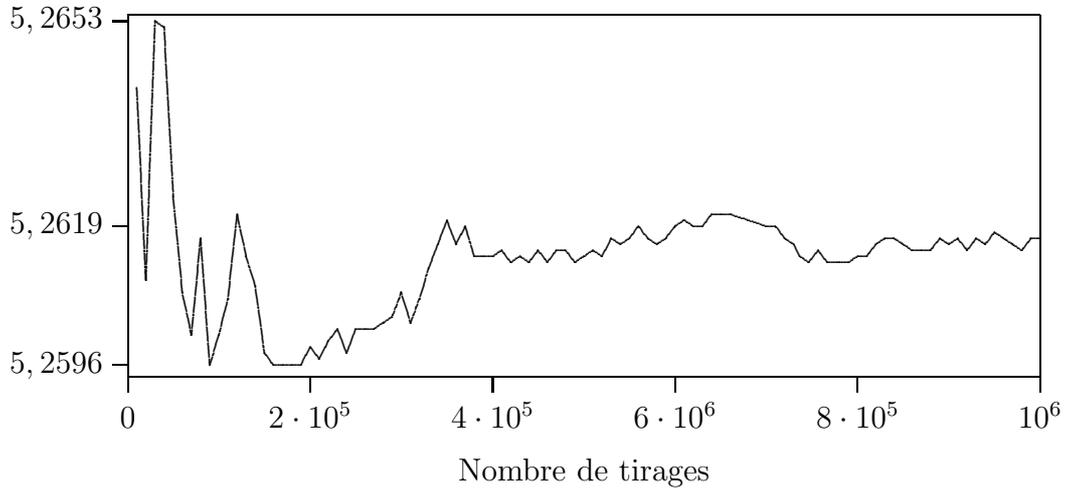

$$\beginpicture
\setcoordinatesystem units <8mm,8mm> point at 0 0
\setplotarea x from 0 to 15, y from 0 to 6
\grid 1 1
\plotheading {\'Evolution de la solution $u$ en fonction du nombre de tirages}
\axis bottom label {Nombre de tirages} ticks
 withvalues $0$ $2\cdot 10^5$ $4\cdot 10^5$ $6\cdot 10^6$ $8\cdot 10^5$ $10^6$ / 
quantity 6 /
%
\axis left             ticks
      withvalues {$5,2596$} {$5,2619$} {$5,2653$} /
at 0.2  2.5  5.9 /
 /
\setlinear
\plot  0.15 4.8   0.3 1.6   0.45 5.9   0.6 5.8   0.75 3   0.9 1.4   1.05 0.7   1.2 2.3   1.35 0.2   1.5 0.7  1.65 1.3  1.8 2.7  1.95 2.0  2.1 1.5  2.25 0.4  2.4 0.2  2.55 0.2  2.7 0.2  2.85 0.2   3 0.5  3.15 0.3  3.3 0.6  3.45 0.8  3.6 0.4  3.75 0.8  3.9 0.8  4.05 0.8  4.2 0.9  4.35 1.0  4.5 1.4  4.65 0.9  4.8 1.3  4.95 1.8  5.1 2.2  5.25 2.6  5.4 2.2  5.55 2.5  5.7 2.0   5.85 2.0  6 2.0  6.15 2.1  6.3 1.9  6.45 2.0  6.6 1.9  6.75 2.1  6.9 1.9  7.05 2.1  7.2 2.1  7.35 1.9   7.5 2  7.65 2.1  7.8 2.0  7.95 2.3  8.1 2.2  8.25 2.3  8.4 2.5  8.55 2.3  8.7 2.2  8.85 2.3  9 2.5  9.15 2.6  9.3 2.5  9.45 2.5  9.6 2.7  9.75 2.7  9.9 2.7  10.5 2.5  10.65 2.5 10.8 2.3  10.95 2.2  11.05 2.0  11.2 1.9  11.35 2.1  11.5 1.9  11.7 1.9  11.85 1.9   12 2  12.15 2 12.3 2.2  12.45 2.3  12.6 2.3  12.75 2.2  12.9 2.1  13.05 2.1  13.2 2.1  13.35 2.3   13.5 2.2  13.65 2.3  13.8 2.1  13.95 2.3  14.1 2.2  14.25 2.4  14.4 2.3  14.55 2.2  14.7 2.1  14.85 2.3   15 2.3  /%
\endpicture$$
\caption{Courbe d'évolution de la solution jusqu'à $10^6$ tirages --- Valeur théorique\,: $u=5,2619$ --- Pas de la marche simulée\,: $h=10^{-3}$}
\label{Courbe}
\end{figure}

\begin{figure}[tbp] 
  \centering
  \includegraphics[bb=88 265 505 575,width=14cm,height=10.4cm,keepaspectratio]{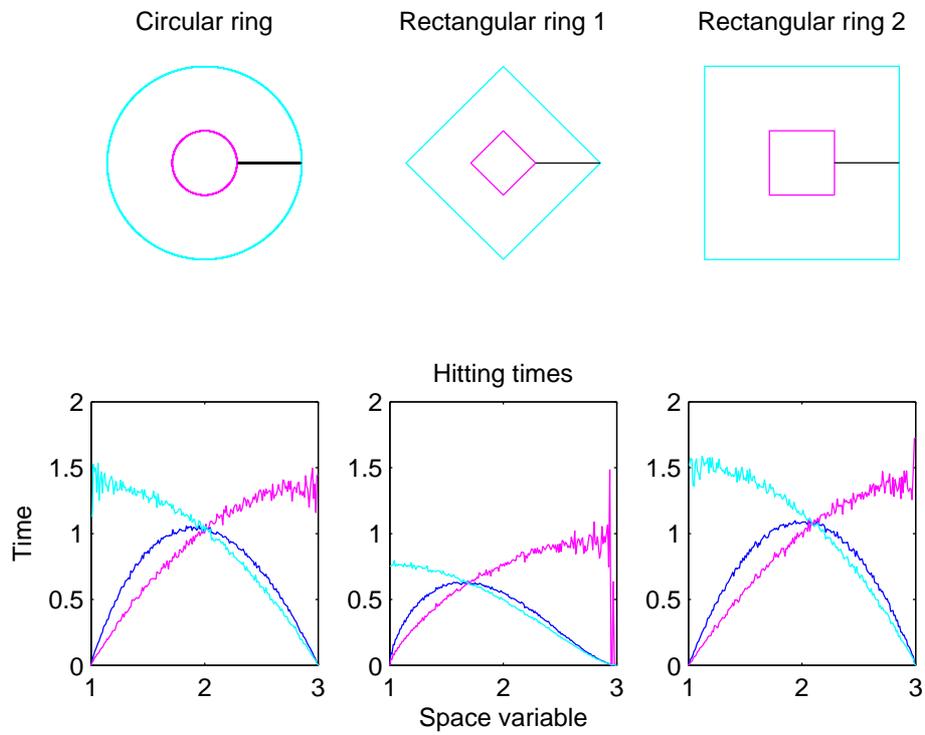}
  \caption{Temps d'atteinte de la fronti\`ere de couronne et des parties internes et externes}
  \label{Fig-C-R1-R2-HittingTimes}
\end{figure}

\begin{table}[h]
\caption{Problème de Dirichlet --- Calcul de la solution en $(2;0)$}
\label{Tableau I}
$$\BeginTable
    \OpenUp11
    \def\C{\JustCenter}
\def\H#1{\JustCenter \Lower{\it #1}}
\BeginFormat
|4 r | r  | r  |4
\EndFormat
\_4
"\JustLeft{Pas d'une marche}| \JustLeft{Valeur calculée} |\JustLeft{Erreur relative} "\\
"          $h$ |            $u_c$ |      $\vert (u_c-u)/u\vert $   "\\+22
\_4
" 0,1~~        | 5,27040          | 1,6 $\cdot 10^{-3}$   "\\
" 0,05~        | 5,26876          | 1,3 $\cdot 10^{-3}$   "\\
" 0,025        | 5,26730          | 1,0 $\cdot 10^{-3}$   "\\
" 0,005        | 5,25766          | 0,8 $\cdot 10^{-3}$   "\\
\_4
\EndTable$$
\end{table}
\noindent et, dans un second temps, l'évolution de la solution calculée jusqu'à $10^6$ tirages avec $h=10^{-3}$ de la figure~\ref{Courbe} où la solution calculée est portée en ordonnée tous les $10^4$ tirages.
\newpage 
Considérons des couronnes circulaires et rectangulaires de la figure \ref{Fig-C-R1-R2-HittingTimes} de la page \pageref{Fig-C-R1-R2-HittingTimes}.

\begin{figure}[tbp] 
  \centering
  \includegraphics[bb=-223 145 819 696,width=14cm,height=7.41cm,keepaspectratio]{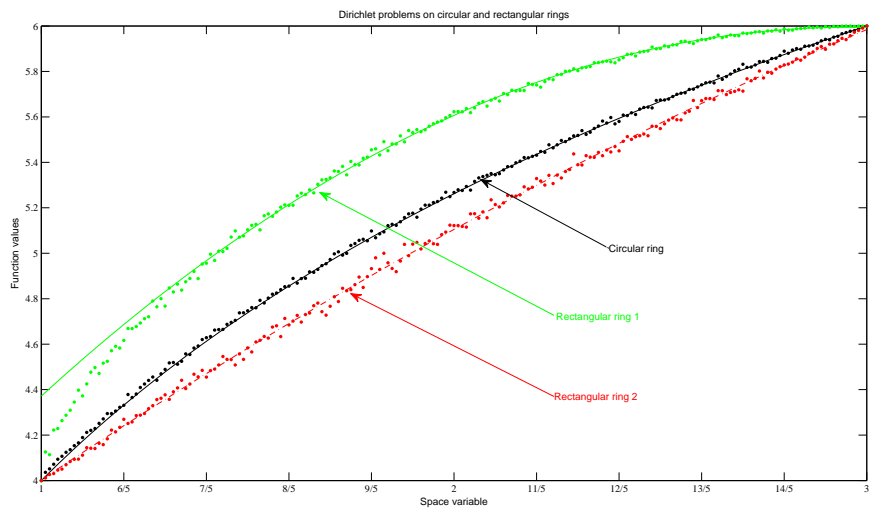}
  \caption{Solutions approchées des trois couronnes}
  \label{fig:Fig3rings}
\end{figure}

\begin{figure}[tbp] 
  \centering
  \includegraphics[bb=-223 145 819 696,width=14cm,height=7.41cm,keepaspectratio]{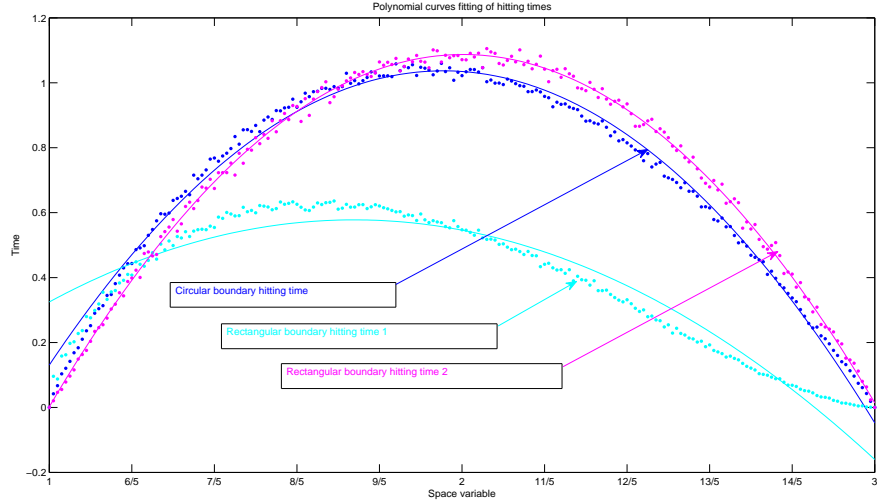}
  \caption{Temps d'atteinte des frontières des trois couronnes}
  \label{fig:Fig3hittingtimes}
\end{figure}

Des résultats des essais numériques sont représentés dans les figures \ref{Fig-C-R1-R2-HittingTimes}, \ref{fig:Fig3rings} et \ref{fig:Fig3hittingtimes} ; ils ont été obtenus avec un pas $h=0.1$ et un nombre de marches aléatoires égal à $5000$ par point.

Dans un premier temps, des évaluations des temps d'atteinte ont été effectuées et reportées sous chaque couronne de la figure \ref{Fig-C-R1-R2-HittingTimes}.

Dans un second temps, les essais numériques ont conduit aux représen\-tations de la solution de la figure \ref{fig:Fig3rings} de la page \pageref{fig:Fig3rings}.

\`A l'aide de l'outil "Curve Fitting Tool" de MATLAB, l'approximation des solutions de la solution conduit aux coefficients suivants (avec un intervalle de confiance à 95\%) :
\begin{enumerate}
  \item Couronne circulaire : $f(x) = a+b\ln(x)/\ln(3)$,\\
   $a = 4.014\quad (4.011, 4.016)$, $b = 1.986\quad (1.981, 1.99)$
  \item Couronne rectangulaire 1 : $f(x) = p_1x^2 + p_2x + p_3$,\\
   $p_1 = -0.4289\quad (-0.4351, -0.4227)$, $p_2 = 2.531\quad (2.505, 2.556)$,\\ $ p_3 = 2.26\quad (2.237, 2.284)$
  \item Couronne rectangulaire 2 : $f(x) = p_1x^2 + p_2x + p_3$,\\
   $p_1 = -0.1137\quad (-0.1201, -0.1072)$, $p_2 = 1.445\quad (1.418, 1.471)$,\\ $ p3 = 2.673\quad (2.648, 2.698)$
\end{enumerate}

Les estimations des temps d'atteinte des frontières des couronnes et leurs courbes de régression quadratique sont données sur la figure \ref{fig:Fig3hittingtimes} de la page \pageref{fig:Fig3hittingtimes}. Les courbes de régression admettent pour équations :
\begin{enumerate}
  \item Couronne circulaire : $f(x) = -0.9932 x^2 + 3.883 x  -2.76$
  \item Couronne rectangulaire 1 : $f(x) = -0.4609 x^2 + 1.6   x -0.811 $
  \item Couronne rectangulaire 2 : $f(x) = -1.081  x^2 + 4.329 x -3.244$
\end{enumerate}

Les temps d'atteinte des frontières et des parties intérieures et extérieures aux couronnes sont précisés dans l'annexe \ref{Annexe2}.


\section{Problème de Dirichlet quasi-linéaire}
\label{Exemple de probleme de Dirichlet non lineaire}

Dans un premier temps, pour valider l'algorithme $A2$, résolvons le pro\-blème\index{pb:Dir:lin} (\ref{Dirichlet-non lineaire}) dans l'intervalle $G = \left] 0 , L \right[$, avec les données suivantes\,:
\begin{eqnarray*}
f(u) & = & \frac{1}{2}\,\exp\left(\frac{u}{1+x}-1\right)  \quad \mbox{dans $G$} \\[3mm] 
g_1(0)& = & 1\,, \qquad g_1(L) = (1+L)\,\biggl(1-\ln(1+L)\biggr)
\end{eqnarray*}
où le second membre $f$ est non linéaire en $u$.

\noindent Ce problème de Dirichlet admet pour solution exacte\index{exacte}\,:
$$
 u(x) = (1+x)\,\biggl(1-\ln(1+x)\biggr)\,.
$$

\`A l'aide de l'algorithme $A2$, les résultats du tableau~\ref{Tableau Inon2} ont été obtenus avec les données suivantes\,:
\begin{center}
\begin{tabular}{lrcl}
Longueur de l'intervalle\,: & $L$   & $=$ & 1  \\
Nombre de points\,: & \verb!MAXPT!$\mbox{}+1$  & $=$ & 51 \\
Nombre de tirages\,:& $   NT      $ & $=$ & $ 10^3  $\\
Pas d'une marche\,: & $ h $ & $=$ & $2\times 10^{-2}$\\
Valeur initiale\,:  & $ u $             &$\equiv$ & 0\\
Nombre d'itérations\,:& \verb!MITER! & $=$ & 10 
\end{tabular}
\begin{equation}\mbox{}\label{data test}\end{equation}
\end{center}
L'intervalle $G$ a été subdivisé en \verb!MAXPT!$\mbox{}=50$ sous-intervalles de même longueur. La solu\-tion a été calculée, à chaque itération\index{itération}, en chacun des points de discré\-tisation. Après dix itérations\index{itérations}, le maximum de l'erreur relative\index{erreur relative} aux points de discrétisation a été calculé et vaut\,: $1,\!62\times 10^{-2}$ au point $0,\!72$.
\begin{table}[hbt]{\small 
\caption{Problème de Dirichlet non linéaire}
\label{Tableau Inon2}
$$\BeginTable    \OpenUp11    \def\C{\JustCenter}
\def\H#1{\JustCenter \Lower{\it #1}}
\BeginFormat |4  r  | r  | r  | r  |4 \EndFormat
\_4
|\JustLeft{Point}  | \JustLeft{Valeur calculée} | \JustLeft{Valeur approchée}| \JustLeft{Erreur relative}|\\+10
| $x$   | $u_c$ | $u$ | $\vert (u_c-u)/u\vert $  |\\+12
\_4
| $0,\!2~$ | $1,\!0013~$ | $0,\!9998~$ | $0,\!15 \times 10^{-2}$ |\\
| $0,\!5~$ | $0,\!8920~$ | $0,\!8918~$ | $0,\!02 \times 10^{-2}$ |\\
| $0,\!72$ | $0,\!79999$ | $0,\!78720$ | $1,\!62 \times 10^{-2}$ |\\
| $0,\!8~$ | $0,\!7415~$ | $0,\!7420~$ | $0,\!07 \times 10^{-2}$ |\\
\_4
\EndTable$$ }
\end{table}

Dans un second temps, considérons le problème\index{pb:Dir:non:lin} (\ref{Dirichlet-non lineaire}), correspondant à (\ref{Sys_u_trois}) issu de \cite{Zwi98}, dans l'intervalle $G = \left] 0 , 1 \right[$, avec les données suivantes\,:
\begin{eqnarray*}
f(u)   & \equiv & -\frac 12 \, a\,u^3 \qquad \mbox{dans } G\\
g(0)   & \equiv & 0 \\
g(1)   & \equiv & 1 
\end{eqnarray*}

\begin{figure}[hbt]
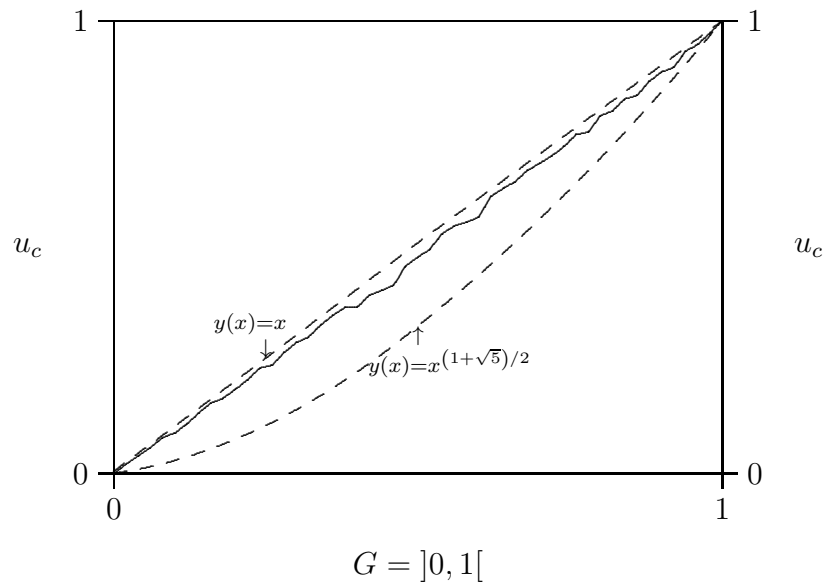

$$\beginpicture
\setcoordinatesystem units <8mm,60mm> point at 0 0
\setplotarea x from 0 to 10, y from 0 to 1
\grid 1 1
\plotheading {Courbe représentative de la solution calculée $u_c$ à l'itération 39}
\put {$y(x)=x~~~\atop \downarrow$} [b] at 2.5 0.255
\put {$\uparrow \atop ~~~~~~ y(x)=x^{\left(1+\sqrt{5}\right)/2}$} [t] at 5 0.326 
\axis bottom label {$G=\left]0,1\right[$} ticks
      withvalues {$0$}  {$1$} /
      at 0 10 /
 / 
\axis left label {$u_c$}  ticks 
      withvalues {$0$} {$1$}  /
      at 0 1  /
 /
\axis right label {$u_c$} ticks 
      withvalues {$0$} {$1$} /
      at 0 1 /
 /
\setlinear
\plot 0 0 0.2 0.02   0.4 0.039  0.6  0.057 0.8 0.079 
1 0.089 1.2 0.109 1.4 0.134 1.6 0.156 1.8 0.166 
2 0.186 2.2 0.207 2.4 0.233 2.6 0.239 2.8 0.267 
3 0.289 3.2 0.301 3.4 0.328 3.6 0.349 3.8 0.367 
4 0.367 4.2 0.394 4.4 0.404 4.6 0.416 4.8 0.459 
5 0.478 5.2 0.496 5.4 0.53  5.6 0.547 5.8 0.555 
6 0.567 6.2 0.612 6.4 0.63  6.6 0.645 6.8 0.67  
7 0.685 7.2 0.702 7.4 0.722 7.6 0.749 7.8 0.754 
8 0.791 8.2 0.801 8.4 0.828 8.6 0.835 8.8 0.866 
9 0.887 9.2 0.897 9.4 0.933 9.6 0.948 9.8 0.973 10 1.0   /
\setdashes
\setquadratic
\plot 0 0.005   2.5 .255   5 .51   7.5 .755   10 1 /
\plot 0 0 2.5 0.106 5 0.326 7.5 0.628 10 1  /
\endpicture$$
\caption{Représentation de la solution du problème de Dirichlet non linéaire\,: $u^{\prime\prime}=u^3$, $u(0)=0$, $u(1)=1$ en trait continu}
\label{Courbe:u:trois}
\end{figure}

\noindent Lorsque $a=  1$, ce problème de Dirichlet admet une solution qui vérifie \cite{Zwi98}\,:
\begin{equation}
\forall\,x\in G, ~ x^{\left(1+\sqrt{5}\right)/2} \leq u(x) \leq x \,.
\label{test:arret}
\end{equation}

\`A l'aide de l'algorithme $A2$ (mais sans relaxation), les résultats de la figure~\ref{Courbe:u:trois} ont été obtenus avec les données suivantes\,:
\begin{center}\begin{tabular}{lrcl}
Nombre de points\,: & \verb!MAXPT!$\mbox{}+1$ & $=$ & 51\\
Nombre de tirages\,:& $   NT      $ & $=$ & $ 2\times 10^3   \,$\\
Pas d'une marche\,: & $ h $ & $=$ & $2\times 10^{-2}$\,\\
Valeur initiale\,:  & $ u $ &$\equiv$ & 0\\
\end{tabular}\end{center}

Le nombre d'itérations\index{itérations} \verb!MITER! a été remplacé par le test d'arrêt\index{test arret}\,:\\
$x_i^{\left(1+\sqrt{5}\right)/2} \leq u(x_i) \leq x_i $ est réalisé (\mbox{cf.} (\ref{test:arret})) {\bf et} $u(x_{i-1}) < u(x_i)$, pour tout $i=2$,\ldots, \verb!MAX!, c'est-à-dire que la solution calculée aux points $x_i$ de discrétisation est bien encadrée et est strictement croissante. Sur la figure \ref{Courbe:u:trois}, le test est satisfait à l'itération 39.

\begin{figure}[tbp] 
  \centering
  \begin{minipage}[c]{.46\linewidth}
   \includegraphics[bb=88 264 505 577,width=5.25cm,height=5.25cm,keepaspectratio]{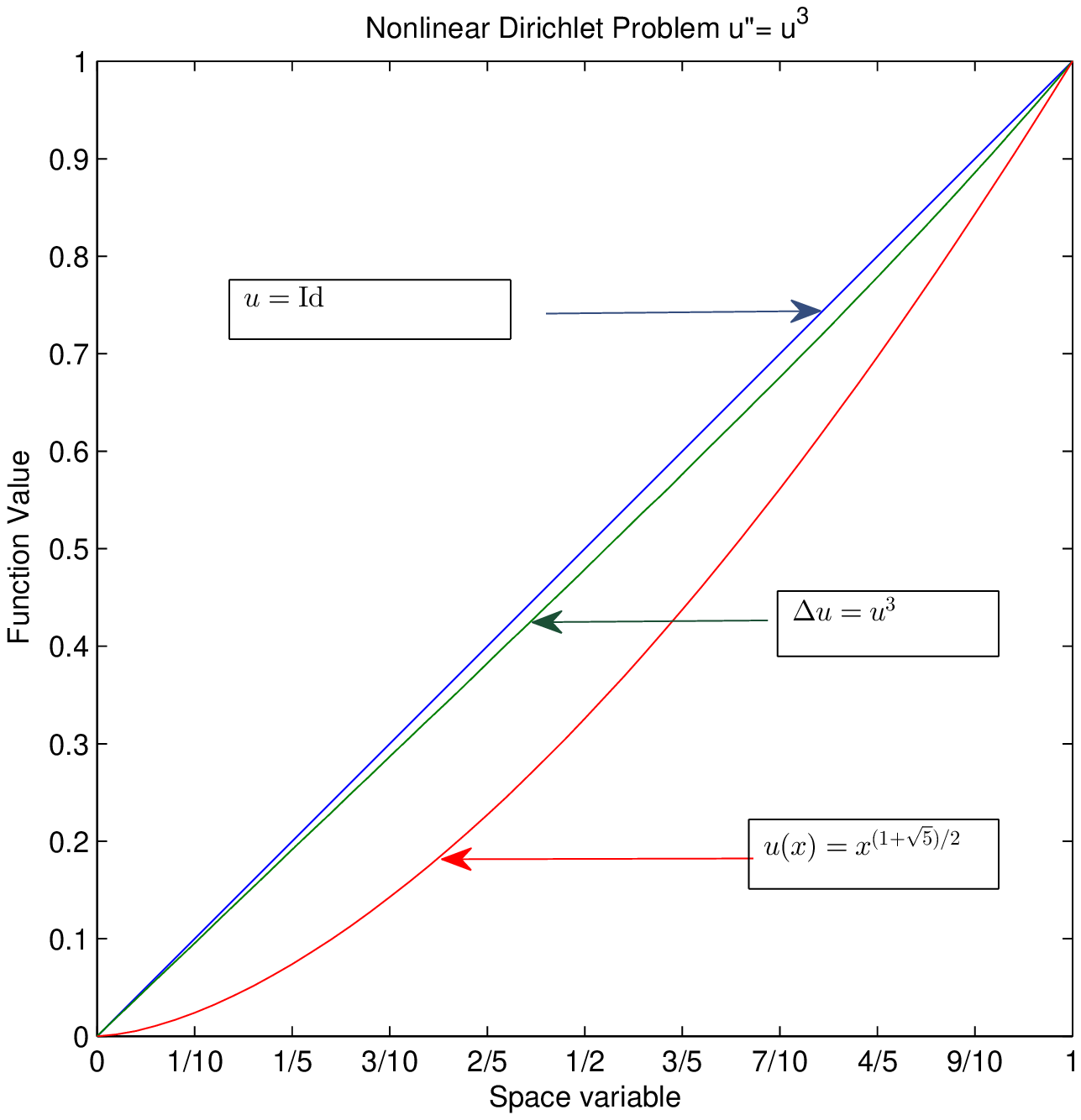}
   \end{minipage} \hfill
   \begin{minipage}[c]{.46\linewidth}
   \includegraphics[bb=88 264 505 577,width=5.25cm,height=5.25cm,keepaspectratio]{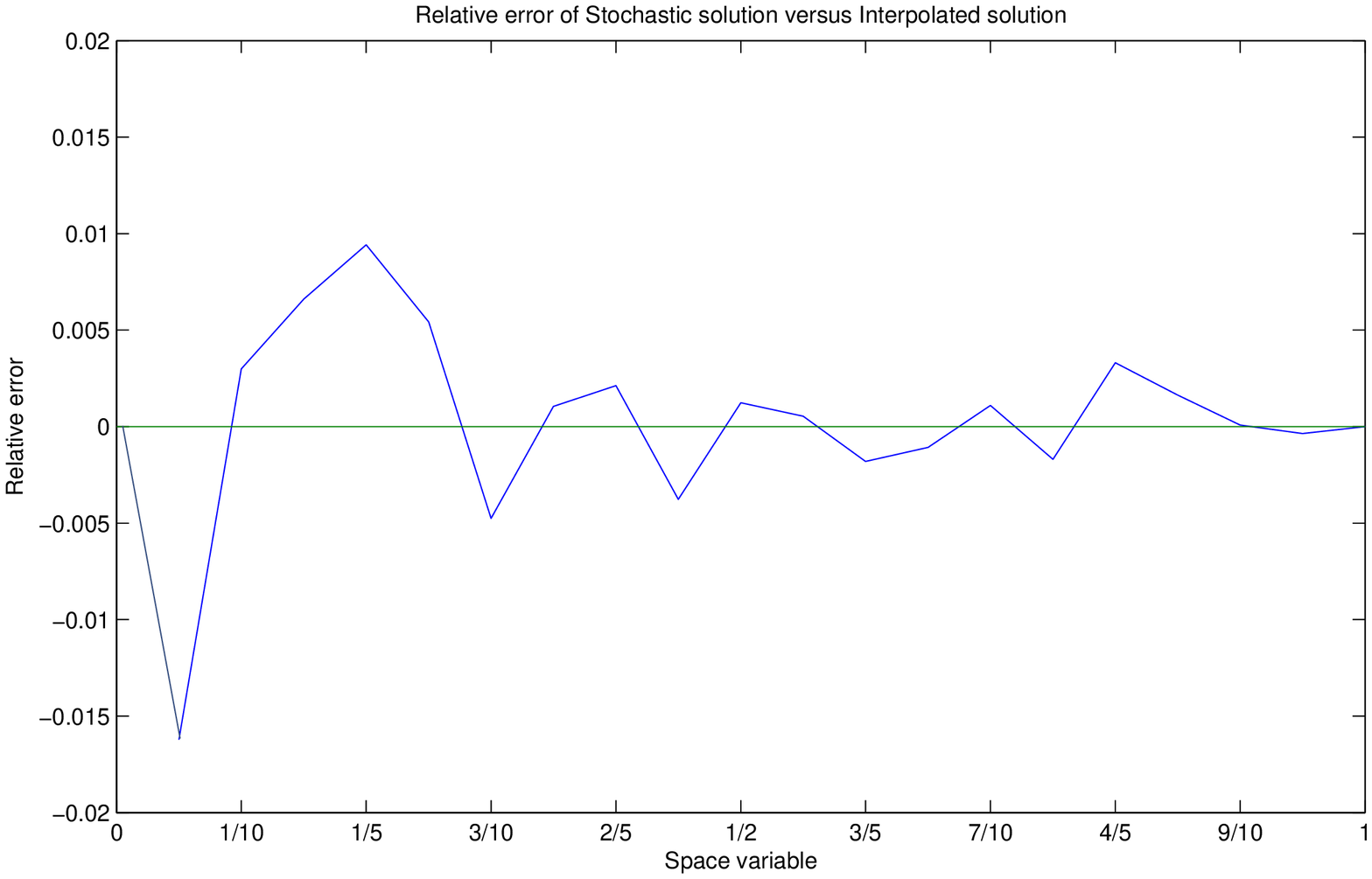}
\end{minipage}
  \caption{Solution pour $a=1$ - Erreur relative empirique}
  \label{fig_a_vaut_1}
\end{figure}

Avec le logiciel MATLAB, le programme prend la forme suivante\,:
\begin{verbatim}
% Problème non linéaire avec condition de Dirichlet
clear all
rng('shuffle','v5uniform');
tic
% Domaine G = ] 0, 1 [
l = 1.0 ;
% Subdivision en maxpt intervalles
maxpt  = 20 ;
maxpt1 = maxpt + 1 ;
h = l / maxpt ; % Pas de la marche aléatoire
% Nombre d'itérations
miter = 500 ;
% Nombre de marches aléatoires
nt = 20000 ;
% Initialisation de U
U  = 0 : h : 1 ;
% Fonctions encadrantes et lissage
x = 0:maxpt;
y = (h.*x).^(0.5*(1+sqrt(5)));
z = h.* x ;
hspline = 0.25*h;
% Itération du problème non linéaire
for iter = 1:miter
    for J = 2:maxpt
        ys  = 0.0 ; dir = 0.0 ;
        for I = 1 : nt % Nb de marches aléatoires par point
            m = J;  ya = 0.0;
            while ( (m > 1) && (m < maxpt1 ) )
                ya = ya - 0.5* U( m )* U( m )* U( m ) ;
                m  = m + (-1)^randi([0, 1]) ;
            end
            if m > 1 dir = dir + 1.0 ; end
            ys  = ys  + ya ;
        end
        U(J)  = ( h*h*ys + dir ) /nt ;
    end
    figure(iter) hold all cla
    xx = 0:hspline:maxpt;
    yy = spline(x,U,xx);
    plot(xx,yy)
    plot(x,y);plot(x,z)
end
toc
\end{verbatim}
Pour la valeur $a=1$, nous avons obtenu une courbe représentative de la solution calculée numériquement par la méthode stochastique et l'erreur relative par rapport à des calculs par différences finies suivant la fonction MATLAB 'bvpc' appliquée à de nombreuses valeurs initiales dans la commande MATLAB 'bvpinit'. Des courbes sont données dans les figures \ref{fig_a_vaut_1} où l'erreur relative empirique atteint un maximum inférieur à 2\% au voisinage de zéro et où la méthode stochastique est régularisante par rapport à la méthode aux différences finies (à condition de prendre des nombres de marches aléatoires et d'itérations suffisamment grands, ce qui a un coût en temps de calcul).


\begin{figure}[tbp] 
  \centering
  \begin{minipage}[c]{.46\linewidth}
   \includegraphics[bb=88 264 505 577,width=5.25cm,height=5.25cm,keepaspectratio]{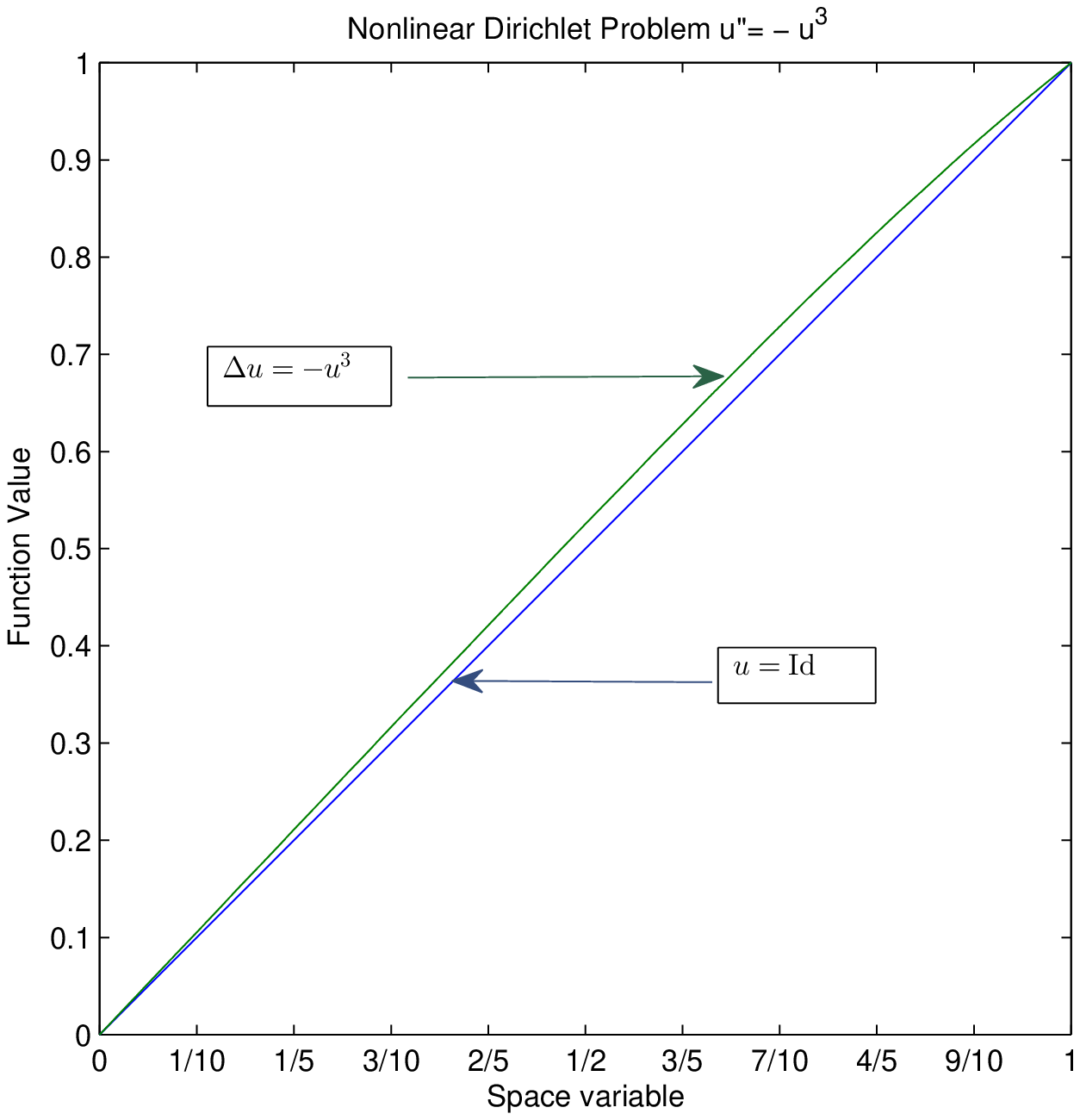}
   \end{minipage} \hfill
   \begin{minipage}[c]{.46\linewidth}
   \includegraphics[bb=88 264 505 577,width=5.25cm,height=5.25cm,keepaspectratio]{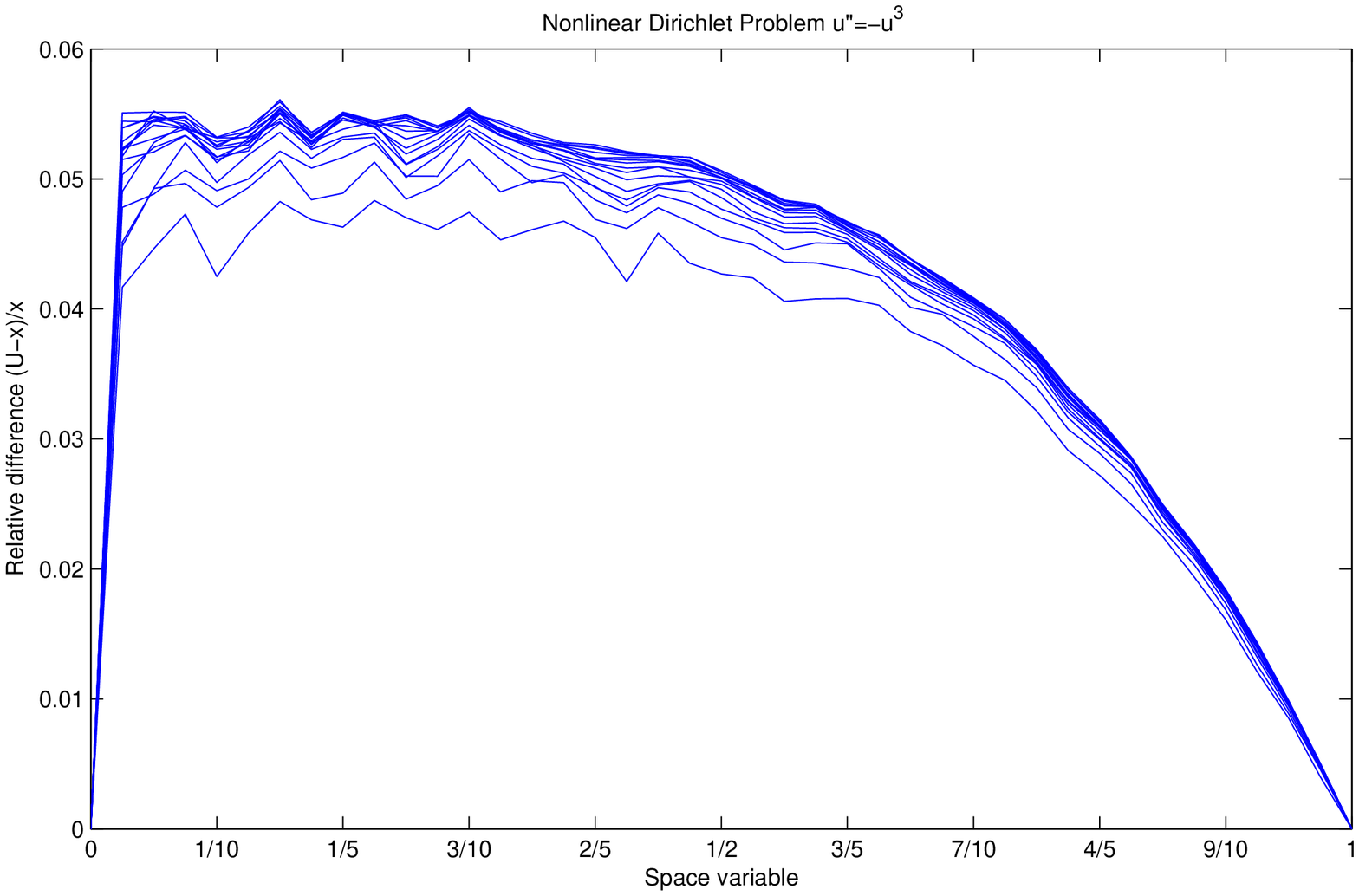}
\end{minipage}
  \caption{Solution pour $a=-1$ - Courbes de convergence numérique}
  \label{Fig_a_vaut_moins_1}
\end{figure}

Compte tenu des difficultés de convergence numérique et du temps de calcul, le programme de calcul a été amélioré en premier lieu par le choix aléatoire du point de départ des marches aléatoires, ce qui évite le biais du choix croissant des points, et par l'introduction d'un calcul de moyennes progressives à l'aide du poids \verb!bary! qui permet de cumuler l'effet des marches aléatoires successives. Par rapport au programme précédent, le programme ci-après a donné une convergence numérique plus rapide, en particulier pour $a=-1$ (voir la figure \ref{Fig_a_vaut_moins_1} où l'évolution de l'erreur relative empirique $(u-x)/x$ est reportée).
\begin{verbatim}
% Itération du problème non linéaire
miter1 =  miter * maxpt ;
for iter = 1:miter1
    m   = randi( [ 2 , maxpt ] ) ; % Choix aléatoire    
    ys  = 0.0 ; dir = 0.0 ;
    for I = 1 : nt % Nb de marches aléatoires par point
        indice = m  ; ya = 0.0;
        while ( (indice > 1) && (indice < maxpt1 ) )
            ya = ya + U3demi( indice ) ; % + au lieu de -
            indice = indice + (-1)^randi([0, 1]) ;
        end
        dir = dir + (indice>1);
        ys  = ys  + ya ;
    end
    U(m) = ( (h2*ys + dir )/nt + bary(m)*U(m) )/(bary(m)+1) ; 
    U3demi(m) = 0.5 .* U(m).* U(m).* U(m) ;
    bary(m) = bary(m) + 1 ;
end
\end{verbatim}
Ce programme a donné les résultats représentés dans les figures \ref{Fig_a_vaut_moins_1}.

\noindent Pour faciliter l'écriture d'un programme de calcul parallèle avec MATLAB, nous effectuons le changement de variable suivant\,:
$$
v = k\, u \qquad \mbox{avec $k = \mbox{sgn}(a)\sqrt{ |a| }$ }
$$
Alors le système devient\,:
\begin{eqnarray}
\left\{
       \begin{array}{rcll}
\Delta v    & = & \mbox{sgn}(a)\,v^3  & \mbox{dans } G       \\[2mm]
       v(0) & = & 0 & \\[2mm]
       v(1) & = & \mbox{sgn}(a)\sqrt{  |a| }   
       \end{array}
\right.  & \mbox{} & \mbox{}
\label{Variable_v}
\end{eqnarray}
Soit $M > 0$ assez grand. On cherche une solution $v$ bornée par $M$, dérivable à l'ordre deux, à dérivées continues et bornées sur l'intervalle fermé $G^*=[0,1]$. Alors on a:
$$
v \in {\cal C}_b^{2}\left(\,G^*\,\right) \Longrightarrow f\left(v\right) \in {\cal C}_b^{2}\left(\,G^*\,\right)
$$
où $f(v) =  \mbox{sgn}(a)\,v^3 $.

La fonction $f$ est lipschitzienne en $v$, uniformément sur $G$, ce qui satisfait les conditions du Théorème \ref{Theoreme:proba} du chapitre \ref{proba:fonctionnel}.

Ces propriétés sont présentées dans l'annexe \ref{Annexe1}.

\begin{figure}[tbp] 
  \includegraphics[bb=86 231 683 609,width=18cm,height=14cm,keepaspectratio]{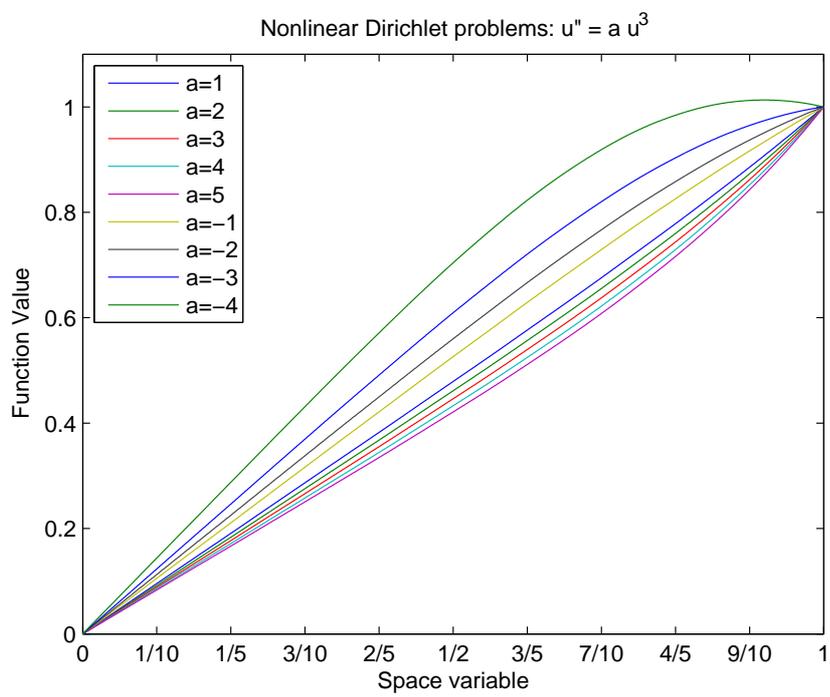}
  \caption{Solutions pour $a\in{-4,-3,-2,-1,1,2,3,4,5}$}
  \label{Fig_fuseau}
\end{figure}

La méthode a conduit aux résultats de la figure \ref{Fig_fuseau} avec une subdivision en 20 intervalles et un nombre de marches aléatoires égal à $20\,000$ en moyenne par point.

Les valeurs de $a$ correspondant à $u'(1)=0$ et à $u'(1)=-1$ ont été obtenues avec la méthode stochastique et les fonctions bvp4c de Matlab et NDSolve de Mathematica (voir le tableau \ref{Tableau_yprime}).

\begin{table}[hbt]{\small 
\caption{Valeurs de $a$ pour lesquelles $u'(1)=0$ et $u'(1)=-1$}
\label{Tableau_yprime}
$$\BeginTable    \OpenUp11    \def\C{\JustCenter}
\def\H#1{\JustCenter \Lower{\it #1}}
\BeginFormat |4  r | l | l | l |4 \EndFormat
\_4
|  |\JustLeft{Méthode}  | \JustLeft{MATLAB} | \JustLeft{MATHEMATICA}|\\+10
| \JustLeft{$u'(1)$}|\JustLeft{stochastique}  | \JustLeft{bvp4c} | \JustLeft{NDSolve}|\\+12
\_4
| $ 0 $ | $a = -3.43755656$  | $a=-3.43761588$ | $a=-3.43759287$ |\\
| $ -1$ | $a = -4.39$        | $a=-4.33564372$ | $a=-4.33553126$ |\\
\_4
\EndTable$$ }
\end{table}

\begin{figure}[tbp] 
  \centering
  \includegraphics[bb=-6 245 602 596,width=14cm,height=8.05cm,keepaspectratio]{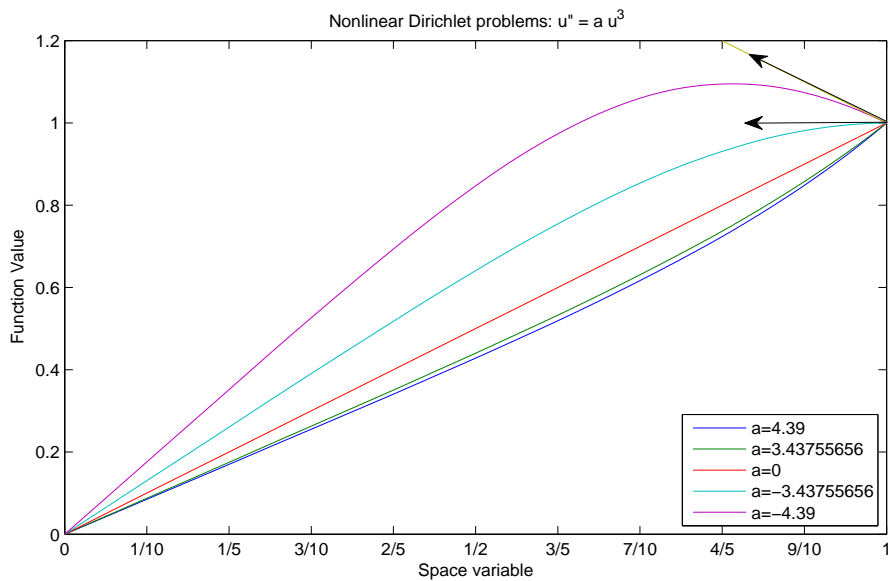}
  \caption{Solutions pour $u'(1)=0$ et $u'(1)=-1$}
  \label{Fig_yprime_nul}
\end{figure}

Pour ces deux cas, la méthode stochastique a conduit aux résultats de la figure \ref{Fig_yprime_nul} avec une subdivision en 20 intervalles et un nombre de marches aléatoires égal à $10\,000$ en moyenne par point.

\section{Commentaires}

Pour $a>0$, on trouve des solutions très proches (moins de 1\%) de celles obtenues avec les logiciels Mathematica (commande NDSolve, en particulier avec les méthodes "Shooting" et "ExplicitRungeKutta") et MATLAB (en particulier avec la commande bvp5c). Les courbes ont donc la même forme convexe que celle de la figure \ref{fig:Fig_a_positif} de l'annexe \ref{Annexe1}.

Pour $a<0$, on trouve des solutions concaves pour $a>-4.5$. Par contre, $a=-4.5$ et $a=-5$ conduisent à une convergence numérique "lente" dans un premier temps (par accumulation de marches aléatoires), puis à une divergence numérique de manière "explosive" dans un second temps (où la solution atteinte numériquement est prise comme valeur initiale). Par conséquent, il ne faut pas se limiter à une convergence numérique, mais il est nécessaire de réitérer la méthode en prenant la valeur obtenue dans un premier temps comme valeur initiale ; le "bruit blanc" du mouvement brownien permet de tester la stabilité de cette valeur obtenue lors de la première étape ; si cette valeur est stable, alors une solution est atteinte ; sinon, la méthode ne permet pas d'affirmer que la valeur obtenue dans un premier temps soit une solution (même instable).

Pour $a = -3.43755656$  et $a = -4.39$, les courbes obtenues ont les formes des figures \ref{Fig_a_negative1} de l'annexe \ref{Annexe2}. La forme  \ref{Fig_a_negative2} de l'annexe \ref{Annexe2} n'a jamais été obtenue lors des essais numériques de la méthode stochastique, sauf de manière transitoire. La méthode stochastique privilégie les solutions les "plus stables" en un sens qui reste à définir. Les autres solutions sont soit instables, soit métastables quand elles possèdent un domaine d'attraction réduit, si ce domaine existe (les calculs n'ont pas permis de trouver de tels domaines). 

 \chapter{Approche fonctionnelle}\label{proba:fonctionnel}

 \section{Introduction}\label{Intro:proba:fonctionnel}

Le problème de Dirichlet non linéaire (\ref{Dirichlet-non lineaire}) a été résolu numériquement à l'aide de la représentation (\ref{representation potentiel non lineaire de u}--\ref{representation potentiel non lineaire de Y}) par reports successifs des valeurs obtenues. Ce chapitre présente une approche fonctionnelle de ce type de représentation sous des hypothèses supplémentaires d'existence et de régula\-rité\index{hypo:régularité} de la solution. Cette approche par une méthode de type Picard\index{méthode Picard} est utilisée dans \cite{freidlin}, exclusivement à l'aide de la formule de Feynman-Kac\index{Feynman}, pour établir certains résultats concernant la solution de problèmes de Cauchy\index{pbCauchy} dans $\mathbb{R}^d$ avec des conditions particulières sur les données et sur les opérateurs. L'objectif principal de ce chapitre est de justifier, autant que possible, l'utilisation des  représentations qui permettent de résoudre numériquement le problème non linéaire (\ref{Dirichlet-non lineaire}).

Nous considérons le domaine $G = \left] 0 , L \right[$ où $L > 0$ est une constante donnée\,; nous notons $x$ la variable d'espace appartenant à $G$\,; $n$ désigne la normale\index{normale int} unitaire intérieure définie sur  $\partial G$.

Considérons le problème\index{pb:Dir:non lin}\,:
\begin{eqnarray}
\left\{
       \begin{array}{rcll}
-\frac12 \, \Delta u  & = & f(u)   & \mbox{dans } G       \\[2mm]
                    u & = & g & \mbox{sur } \partial G       \end{array}
\right.  & \mbox{} & \mbox{}
\label{sys:ptfixe}
\end{eqnarray}
où l'inconnue $u$ est définie de $G$ dans $\mathbb{R}$ et les données sont\,:
\begin{eqnarray*}
f\,: & &\mathbb{R} \longrightarrow \mathbb{R}  \\
     & &  u \longmapsto f(u) \\
g\,: &  & \partial G \longrightarrow \mathbb{R} \\
     &  & x \longmapsto g(x) 
\end{eqnarray*}

Une représentation stochastique permet l'application d'un théorème de point fixe\index{th:pt:fixe}. Ceci conduit à montrer dans le paragraphe~\ref{Reso:proba:fonctionnel}, l'existence et l'unicité de la solution $u$ de (\ref{sys:ptfixe}) avec un second membre non linéaire $f(u)$ sous des conditions sur $f$.
\section{Une méthode de point fixe}
\label{Reso:proba:fonctionnel}
Sur le plan fonctionnel, il s'agit de se ramener à traiter des représentations stochastiques de la solution de (\ref{sys:ptfixe}) par une méthode de point fixe.

La suite de ce paragraphe concerne principalement la résolution du sys\-tème (\ref{sys:ptfixe}) avec $f$ suffisamment régulière, par utilisation d'un théorème de point fixe de type Picard\index{th:pt:fixe:Picard}.

On introduit les espaces fonctionnels suivants\,:
\begin{Defliste}{(iii)}
  \item[(i)]{\boldmath ${\cal C}_b^k\left(\,G\,\right)$}, respectivement {\boldmath ${\cal C}_b^k\left(\,\partial G\,\right)$}, l'espace des fonctions bornées, $k$ fois dérivables en $x\in G$, respectivement $x\in \partial G$, à dérivées bornées\,;
  \item[(ii)]{\boldmath ${\cal C}_b\left(\,G\,\right)$},  l'espace  des fonctions continues et bornées, muni de la norme sup\,:
\begin{eqnarray}
u & \mapsto & \vert u \vert_\infty = \max \{\, \vert u(x) \vert \,:\, x \in G \,\}.
\label{norme:sup}
\end{eqnarray}
\end{Defliste}

\noindent Les étapes de la méthode peuvent être définies ainsi :
\begin{enumerate}
\item on prend $u^{(1)}$ dans  ${\cal C}_b^{2}\left(\,G\,\right)$ et on vérifie que l'on a\,: $f\left(\,u^{(1)}\,\right)$ dans ${\cal C}_b^{2}\left(\,G\,\right)$\,;
\item en supposant que, pour $f\left(\,u^{(1)}\,\right)$  dans ${\cal C}_b\left(\,G\,\right)$, il existe une solution $u^{(2)}$ du problème\,:
\begin{eqnarray}
\left\{
       \begin{array}{rcll}
-\frac12 \, \Delta u^{(2)}  & = & f\bigl( u^{(1)} \bigr)  & \mbox{dans } G       \\[2mm]
                   u^{(2)}  & = & g & \mbox{sur } \partial G       \end{array}
\right.  & \mbox{} & \mbox{}
\label{sys:proba:ordre:deux}
\end{eqnarray}
vérifiant\,: $u^{(2)} \in {\cal C}_b\left(\,G\,\right)$,\\
on obtient une représentation stochastique associée\,;
\item on définit l'application $u^{(2)} = {\cal G}\left(u^{(1)}\right)$ et on montre que ${\cal G}$ a un point fixe unique dans un espace adapté.
\end{enumerate}

On doit vérifier préalablement que, pour tout $u \in {\cal C}_b^{2}\left(\,G\,\right)$, on a $f\left(u\right) \in {\cal C}_b^{2}\left(\,G\,\right)$, pour pouvoir définir la représentation stochastique associée au pro\-blème. 

\noindent{\bf Exemple.}

La fonction $f(u) = a \, u^3$ est de classe ${\cal C}^\infty$ en $u$ et vérifie\,:
\label{prolongement:C:infini}\begin{equation}
u \in {\cal C}_b^{2}\left(\,G\,\right) \Longrightarrow f\left(u\right) \in {\cal C}_b^{2}\left(\,G\,\right)\,.
\label{prop.2.3:proba}\label{f verifie regularite}
\end{equation}

La seconde étape nécessite la mise en évidence d'une formulation stochastique associée à (\ref{sys:proba:ordre:deux}), noté plus simplement\,:
\begin{eqnarray}
\left\{
       \begin{array}{rcll}
-\frac12 \, \Delta u  & = & f  & \mbox{dans } G      \\[2mm]
                   u  & = & g & \mbox{sur } \partial G      \end{array}
\right.  & \mbox{} & \mbox{}
\label{sys:proba:3}
\end{eqnarray}
où les fonctions $f$ et $g$ vérifient, par hypothèse, les propriétés\,:
\begin{eqnarray}
 f &\in& {\cal C}_b^{2}\left(\,G\,\right)
   \label{f continument derivable} \\[2mm]
 g &\in& {\cal C}_b^{2}\left(\,\partial G \,\right)  
\label{proba:f,g}
\end{eqnarray}

\noindent {\bf Proposition:} On suppose $f$, $g$ donnés vérifiant (\ref{f continument derivable}--\ref{proba:f,g}).

Alors, le système linéaire (\ref{sys:proba:3}) admet une solution unique $u$ dans ${\cal C}_b^{2}\left(\,G\,\right)$.
\label{proba:ladyzenkaja}
{\it Démonstration: } voir \cite{lady}.
                                                     
\noindent Il reste à établir le résultat de cette approche  fonctionnelle.  Nous adaptons au cadre stochastique, à l'aide des hypothèses d'existence et de régularité, les arguments utilisés dans le cadre déterministe.

\begin{thm}
On suppose que\,:
\begin{enumerate}
\item le domaine en espace est défini par:
$$
G = \left]0,L\right[, \quad \mbox{avec~~$ 0 < L < +\infty$}\,;
$$
\item le second membre $f$  vérifie\,:
$$
u \in {\cal C}_b^{2}\left(\,G\,\right) \Longrightarrow f\left(u\right) \in {\cal C}_b^{2}\left(\,G\,\right)
$$
et la fonction $f$ est lipschitzienne\index{lipschitzienne} en $u$, uniformément sur $G$\,;
\item la fonction $g$ vérifie (\ref{proba:f,g})\,;
\item les données du système linéaire  (\ref{sys:proba:3}) sont telles que, pour tout $f$ dans ${\cal C}_b\left(\,G \,\right)$, il existe une solution unique $u$ telle que\,:
\begin{equation}
u \in {\cal C}_b\left(\,G \,\right)\,.
\label{hypothese}
\end{equation}
\end{enumerate}
Alors (\ref{sys:ptfixe}) admet une solution unique $u$ dans ${\cal C}_b\left(\, G \,\right)$, 
qui est limite de fonctions dans ${\cal C}_b^{2}\left(\, G \,\right)$\,. 

Cette solution\index{repre:stochas} $u$ est représentée, pour tout $x \in G\,$, par\,:

\begin{eqnarray*}
u(x) & = & E\left[ \, \int_0^{\tau} f\biggl( u \left( X_s^{x}\right), X_s^{x} \biggr) \right] + E\left[ \, g\left( X_{\tau}^{x} \right) \,\right] 
\end{eqnarray*}

\vspace*{-5mm}
\begin{equation} \label{formulation stochastique}
\end{equation}
\label{Theoreme:proba}
\end{thm}
{\it Démonstration\,: } Pour obtenir une représentation stochastique de la solution, nous restons dans un cadre continu et borné.

Soit $u^{(1)}$ dans ${\cal C}_b\left(\,G \,\right)$. Considérons la fonction $f$ définie sur $G$ par\,:
$f\equiv f\left(\,u^{(1)}\,\right)$. Les hypothèses de la proposition \ref{f verifie regularite} entraînent\,:
$$
f\equiv f\left(\,u^{(1)}\,\right) \in {\cal C}_b\left(\,G \,\right)\,,
$$
puis, d'après l'hypothèse (\ref{hypothese}), on obtient l'existence et l'unicité de la solution, notée $u^{(2)}$, de (\ref{sys:proba:ordre:deux}) où l'on a posé\,:
 $f\equiv f\left(\,u^{(1)}\,\right)$.

Lorsque l'on prend $u^{(1)}$ dans ${\cal C}_b^{2}\left(\,G \,\right)$, alors d'après la proposition \ref{proba:ladyzenkaja}, on a\,: $u^{(2)}\in {\cal C}_b^{2}\left(\, G \,\right)$, et, en appliquant la formule de Itô avec une démonstration similaire à celle du théorème 5.1 de \cite[pages 167--168]{freidlin}, on obtient la représen\-tation suivante\,:

\begin{equation}
u^{(2)}(x) = E\left[ \, \int_0^\tau f\biggl( u^{(1)} \left( X_s^{x} \right) ,  X_s^{x} \biggr) \right] + E\left[ \,  g\left( X_\tau^{x}\right) \,\right] \label{formul:stochastiqueBis}
\end{equation}

De (\ref{hypothese}), l'application\,:
\begin{eqnarray*}
{\cal G}\,:~ {\cal C}_b\left(\,G \,\right) & \longrightarrow & {\cal C}_b\left(\,G \,\right) \\
u^{(1)} & \longmapsto & u^{(2)} = {\cal G} \left(\,u^{(1)}\,\right)
\end{eqnarray*}
où $u^{(2)}$ est la solution unique de (\ref{sys:proba:ordre:deux}) avec  $f\equiv f\left(\,u^{(1)}\,\right)$, est bien définie.

Montrons que ${\cal G}$ est une contraction\index{contraction} sur ${\cal C}_b\left(\,G \,\right)$. Soit $u^{(1)}$, $\overline{u}^{(1)}$ dans ${\cal C}_b\left(\,G \,\right)$ et les solutions correspondantes dans ${\cal C}_b\left(\,G \,\right)$\,: 
$u^{(2)} = {\cal G} \left(\,u^{(1)}\,\right)$, $\overline{u}^{(2)}={\cal G} \left(\,\overline{u}^{(1)}\,\right)$.\\
Nous utilisons la notation générique\index{notation:générique}\,: $\delta q = q - \overline{q}$\,. En particulier, notons\,: $\delta u^{(1)}=u^{(1)}-\overline{u}^{(1)}$, $\delta u^{(2)}=u^{(2)}-\overline{u}^{(2)}$, $\delta f\left(\,u^{(1)}\,\right)=f\left(\,u^{(1)}\,\right) -f\left(\,\overline{u}^{(1)}\,\right)$, $\delta u^{(2)}= \delta {\cal G}\left(\,u^{(1)}\,\right) = {\cal G}\left(\,u^{(1)}\,\right) - {\cal G}\left(\,\overline{u}^{(1)}\,\right)$.

Cherchons un réel $\nu$ dans $\left]0,1\right[$ tel que\,:
$$
\left|\, \delta {\cal G}\left(\,u^{(1)}\,\right) \,\right|_\infty \leq \nu \, \left|\, \delta u^{(1)}  \,\right|_\infty .
$$
Les solutions $ u^{(2)}$ et $\overline{u}^{(2)}$ de (\ref{sys:proba:ordre:deux}) sont telles que leur différence $\delta u^{(2)}$ vérifie\,:
\begin{eqnarray*}
\left\{
\begin{array}{rclc}
- \frac12\,\Delta \left(\, \delta u^{(2)} \,\right)  & = & \left(\,\delta f \,\right) \bigl( u^{(1)} \bigr)  & \left(\,G\,\right)\\[2mm]
  \delta u^{(2)} & = & 0  & \left(\,\partial G \,\right)\,.
\end{array}
\right.
\end{eqnarray*}
D'après l'hypothèse (\ref{hypothese}), ce système admet une  solution unique\,:
$$
\delta u^{(2)} \in {\cal C}_b\left(\,G \,\right)\,.
$$
Alors, pour tout $x\in G$, on a\,:
\begin{eqnarray*}
\left(\, \delta u^{(2)} \,\right)(x) & = & 
E\left[ \, \int_0^\tau \left(\delta f\right)\biggl( u^{(1)} \left( X_s^{x}\right) , X_s^{x} \biggr) \,\right]
\end{eqnarray*}
Cette représentation\index{repres:stoch} est donnée, par exemple, dans \cite{freidlin} pour une solution dans ${\cal C}_b^{2}\left(\,G\,\right)$ d'un problème linéaire. Comme on a seulement $\delta u^{(2)} \in {\cal C}_b\left(\,G \,\right)$, cette représentation est déduite de la formule de Itô\index{Ito}, comme dans \cite{dautray}, en approchant  $\delta u^{(2)}$ par une suite de fonctions régulières dans ${\cal C}_b^{2}\left(\,G\,\right)$.

La représentation précédente entraîne pour tout $x\in G$\,:
\begin{eqnarray*}
\left|\,\delta u^{(2)}(x)\,\right| & \leq & E\left[ \, \int_0^\tau \left|\left(\delta f\right)\left( u^{(1)}  \left(  X_s^{x} \right), X_s^{x} \right)  \right|\,ds \right] \,.
\end{eqnarray*}

Dans la suite de cette démonstration, $M$ désigne une constante générique\index{cte:générique} telle que\,: $0 < M < +\infty\,.$

Comme la fonction $f$ est lipschitzienne\index{lipschitzienne} en $u$, uniformément sur $G$, on a\,:
\begin{eqnarray*}
\left|\,\delta u^{(2)}(x)\,\right| & \leq & M \, E\left[ \, \int_0^\tau \left|\left(\delta u^{(1)} \right) \left( X_s^{x} \right) \right|\,ds \, \right] \,, \quad \forall\,x\in G.
\end{eqnarray*}
On en déduit\,:
\begin{eqnarray*}
\left|\,\delta u^{(2)}(x)\,\right| & \leq & M \, E\left[ \, \int_0^\tau \left|\delta u^{(1)} \right|_\infty\,ds \, \right] \,, \quad \forall\,x\in G\,, 
\end{eqnarray*}
et\,:
\begin{eqnarray*}
\left|\,\delta u^{(2)}(x)\,\right| & \leq & M \,E[\tau]\, \left| \delta u^{(1)} \right|_\infty \,, \quad \forall\,x\in G
\end{eqnarray*}
où $E[\tau] < +\infty$.

Pour terminer la preuve de l'existence d'un réel $\nu$ dans $\left]0,1\right[$ rendant l'application ${\cal G}$ contractante\index{appli:contractante} sur ${\cal C}_b\left(\,G\right)$, nous raisonnons par récurrence\index{récurrence}.\\[4mm]
\noindent {\it \`A l'ordre $m=2$}, on prend l'inégalité   précédente.\\[2mm]

\noindent {\it \`A l'ordre $m$}, on fait l'hypothèse de récurrence suivante\,:
\begin{eqnarray*}
\left|\,\delta u^{(m)}(x)\,\right| & \leq & \frac{\textstyle \left(\,M\,E[\tau] \,\right)^{m-1}}{\textstyle (m-1)!}\, \left| \delta u^{(1)} \right|_\infty \,, \quad \forall\,x\in G.
\end{eqnarray*}
\noindent {\it \`A l'ordre $m+1$}, on a successivement\,:
\begin{enumerate}
\item de la même façon qu'à l'ordre $m=2$\,:
\begin{eqnarray*}
\left|\,\delta u^{(m+1)}(x)\,\right| & \leq & M \, E\left[ \, \int_0^\tau \left|\left(\delta u^{(m)} \right) \left( X_s^{x} \right) \right|\,ds \, \right] \,,
\end{eqnarray*}
\item d'après l'hypothèse de récurrence à l'ordre $m$\,:
\begin{eqnarray*}
\left|\,\delta u^{(m+1)}(x)\,\right| & \leq & M \, E \left[ \, \int_0^\tau \frac{\textstyle \left(\,M\,E[\tau]\,\right)^{m-1}}{\textstyle (m-1)!}\, \left| \delta u^{(1)} \right|_\infty \, ds \,\right] \,,
\end{eqnarray*}
\item par calcul de l'intégrale\,:
\begin{eqnarray*}
\left|\,\delta u^{(m+1)}(x)\,\right| & \leq & \frac{\textstyle M^m}{\textstyle (m-1)!}\,
\left| \delta u^{(1)} \right|_\infty \, E\left[ \, \int_0^\tau E[\tau]^{m-1} \, ds \,\right]  \,,\\[2mm]
\left|\,\delta u^{(m+1)}(x)\,\right| & \leq & \frac{\textstyle M^m}{\textstyle (m-1)!}\,
\left| \delta u^{(1)} \right|_\infty \, E[\tau]^{m}  \,,\\[2mm]
\left|\,\delta u^{(m+1)}(x)\,\right| & \leq & \frac{\textstyle \left(\,M\,E[\tau]\,\right)^m}{\textstyle m!}\,
\left| \delta u^{(1)} \right|_\infty\, \quad \forall\,x\in G.
\end{eqnarray*}
\end{enumerate}
La récurrence étant établie, on en déduit\,:
\begin{eqnarray}
\left|\,\delta u^{(m+1)}\,\right|_\infty  & \leq &  
\frac{\left(\,M\,E[\tau]\right)^{m}}{m!} \, \left|\,\delta u^{(1)} \,\right|_\infty \,.
\label{recurrence:proba}
\end{eqnarray}
Or la série $\displaystyle\sum_{m=1}^\infty \frac{\left(\,M\,E[\tau]\,\right)^m}{m!}$ converge, donc il résulte de (\ref{recurrence:proba}) que $\left(\,u^{(m)}\,\right)_{m\geq 1}$ forme une suite de Cauchy qui converge dans l'espace ${\cal C}_b\left(\,G\,\right)$ vers une limite unique $u$.

Par récurrence, nous définissons l'itérée\index{itérée} de ${\cal G}$ à l'ordre $m$, notée ${\cal G}^{(m)}$. Pour $m$ suffisamment grand, ${\cal G}^{(m)}$ est contractante de coefficient $\nu$ tel que\,:
$$
\left|\, \delta {\cal G}^{(m)} \left(\,u^{(1)}\,\right) \,\right|_\infty \leq \nu \, \left|\, \delta u^{(1)}  \,\right|_\infty\, .
$$
L'application contractante ${\cal G}^{(m)}$ admet donc un point fixe unique\index{pt fixe} dans ${\cal C}_b\left(\,G \,\right)$. D'après un corollaire du théorème du point fixe de Banach\index{Banach}, on sait que ${\cal G}$ admet aussi un point fixe unique $u$ dans ${\cal C}_b\left(\,G \,\right)$ tel que\,:
$ u = {\cal G}(u) $\,.

Finalement, on en déduit que $u$ est la solution unique de (\ref{sys:ptfixe}), et que $u$ est limite de fonctions de ${\cal C}_b^{2}\left(\,G\,\right)$, représentées sous la forme (\ref{formul:stochastiqueBis}).

La représentation (\ref{formulation stochastique}) résulte de (\ref{formul:stochastiqueBis}) et de $ u = {\cal G}(u) $\,.
\qed
%
  \section{Commentaires}\label{comment:proba}
Les résultats généraux d'Analyse Fonctionnelle appliqués dans ce chapitre peuvent être trouvés, par exemple, dans \cite{lady,dautraylions}. 

Le résultat du théorème \ref{Theoreme:proba} peut s'étendre à des conditions aux limites de Neumann ou de Robin avec un ouvert borné $G$ de $\mathbb{R}^d$, dont la frontière $\partial G$ et la normale intérieure $n$ sont suffisamment régulières\index{données à la frontière régulières} (voir les conditions nécessaires à l'existence de ce genre de représentations stochastiques dans \cite[pages 166 et 254]{freidlin} et \cite[page 117]{dautray} où $\partial G$, respectivement $n$, est de classe ${\cal C}^3$, \mbox{resp.} ${\cal C}^2$).

\chapter{Conclusion}
\label{Conclusion}
Des méthodes stochastiques de calcul de solutions de problèmes détermi\-nistes de \mbox{Dirichlet} linéaires et non linéaires ont été présentées. Les représenta\-tions résultent d'une application rigoureuse de la formule de Itô pour les problèmes réguliers et de son application formelle pour les problèmes à domaines non réguliers. Les cas linéaires et quasi-linéaires peuvent ainsi être traités. Ces méthodes sont basées sur des représentations stochastiques qui donnent directement des algorithmes aisément program\-mables. 

Les programmes sont courts, faciles à construire et à vérifier pas à pas. On évite d'entrer en mémoire un maillage de discrétisation du domaine et de gérer les tableaux de numérota\-tion qui l'accompagnent.

Suivant le problème considéré (géométrie du domaine, expressions de la fonction source $f$ et des conditions aux limites $g$\,), le temps de calcul est plus ou moins long avec le processeur séquentiel que nous utilisons.

Comme les méthodes de Monte-Carlo classiques, cette méthode stochastique admet une vitesse de convergence en $1/\sqrt{NT}$. Les problèmes  quasi-linéaires nécessitent davantage de tirages.
Des propriétés essentielles des méthodes de Monte-Carlo sont conservées : d'une part, le calcul de la solution en un point choisi indépendamment des autres dans le cas linéaire, d'autre part, l'adaptation au calcul parallèle dans les cas linéaires et non-linéaires. Mais notons également que les marches aléatoires sont simulées simplement à partir d'épreuves répétées de type Bernouilli.


La résolution numérique basée sur des représentations stochastiques peut s'appliquer en dimension deux ou trois sans difficultés autres que celles liées au temps de calcul. Les représentations stochastiques peuvent aussi être associées à la méthode itérative de Picard\index{méthode de Picard} (voir le chapitre \ref{proba:fonctionnel} pour des résultats partiels d'existence et d'unicité de la solution).

Sur le plan fonctionnel, ces méthodes stochastiques ont nécessité des propriétés supplémentaires de régularité, en particulier la dérivabilité\index{dérivabilité} de la solution au sens classique, et des propriétés sur les données telles que le système linéaire (\ref{sys:proba:3}) admette une solution dans ${\cal C}_b\left(\, G \,\right)$. Concernant la régularité des solutions de problèmes aux limites non linéaires ou de problèmes posés dans des ouverts à frontière non régulière, signalons que dans \cite{lady} et \cite{freidlin} pour des problèmes paraboliques non linéaires et linéaires, des estimations et des représentations stochastiques sont établies pour des solutions suffisamment régulières, de classe ${\cal C}^2$ en espace.

\paragraph{Remerciements}Ce travail a bénéficié du soutien de l'Université de La Réunion.

\appendix
\chapter{Système non linéaire $\Delta u = a\,u^3$, $u(0)=0$, $u(1)=1$}
\label{Annexe1}

Considérons le système non linéaire (\ref{Sys_u_trois}) :
\begin{eqnarray}
\left\{
       \begin{array}{rcll}
\Delta u  & = & a\, u^3   & \mbox{dans } G=]0;1[       \\[2mm]
     u(0) & = & 0 & \\[2mm]
     u(1) & = & 1 & \end{array}
\right.  & \mbox{} & \mbox{}
\label{Sys_u_trois_Annexe}
\end{eqnarray}
avec $a\in \mathbb{R}$.

\noindent Lorsque $a= 0$, ce problème admet la solution évidente $u(x)=x$.

\noindent Lorsque $a=  1$, ce problème de Dirichlet admet une solution qui vérifie \cite{Zwi98}\,:
$$
\forall\,x\in G, ~ x^{\left(1+\sqrt{5}\right)/2} \leq u(x) \leq x.
$$

\section{Changement de variable et propriétés}
Le changement de variable :
$$
v = k\, u \qquad \mbox{avec $k = \mbox{sgn}(a)\sqrt{ |a| }$ }
$$
donne le système suivant :
\begin{eqnarray}
\left\{
       \begin{array}{rcll}
\Delta v    & = & \mbox{sgn}(a)\,v^3  & \mbox{dans } G       \\[2mm]
       v(0) & = & 0 & \\[2mm]
       v(1) & = & \mbox{sgn}(a)\sqrt{  |a| }   
       \end{array}
\right.  & \mbox{} & \mbox{}
\label{Variable_v2}
\end{eqnarray}
On suppose que $v$ est bornée par $M > 0$ assez grand. On s'intéresse à une solution $v$ bornée par $M$, dérivable à l'ordre deux, à dérivées continues et bornées sur $[0,1]$. 
\begin{proposition}
Soit ${\cal C}_b^{2}$ l'ensemble des fonctions définies sur $[0,1]$, dérivables à l'ordre deux, à dérivées continues et bornées.
Pour tout $a\in\mathbb{R}$, on a :
$$
v \in {\cal C}_b^{2}\left(\,[0,1]\,\right) \Longrightarrow f\left(v\right) \in {\cal C}_b^{2}\left(\,[0,1]\,\right)
$$
où $f(v) =  \mbox{sgn}(a)\,v^3 $.
\end{proposition}
\noindent{Preuve :} évidente.

Pour satisfaire les conditions du Théorème \ref{Theoreme:proba} du chapitre \ref{proba:fonctionnel}, montrons la proposition suivante. 
\begin{proposition}
La fonction $f$ définie par $f(v) =  \mbox{sgn}(a)\,v^3 $ est lipschit\-zienne en $v$, uniformément sur $G$.
\end{proposition}
\noindent{Preuve :} Soient $v_1$, $v_2 \in {\cal C}_b^{2}\left(\,[0,1]\,\right)$.
\begin{eqnarray*}
\left|\, f(v_1) - f(v_2) \,\right| & \le & |\mbox{sign}(a)|\,\left|\, v_1^3 - v_2^3 \,\right| \\[2mm]
  & \le & \max_{[0,1]} \left|\, v_1^2 + v_1v_2 + v_2^2 \,\right| \,\times \, \left|\, v_1 - v_2 \,\right| \\[2mm]
 & \le & 3\,M^2 \left|\, v_1 - v_2 \,\right|
\end{eqnarray*}\qed

\begin{proposition}
Soit $v$ une fonction continue bornée solution de (\ref{Variable_v2}).

On suppose que $v$ n'est pas identiquement nulle sur un intervalle $[0,\lambda]$, $0<\lambda<1$.
\begin{enumerate}
\item $v"$ est continue et bornée sur $[0,1]$.
\item Si $a>0$, $v'$ ne s'annule pas sur $[0,1]$ et v est strictement croissante sur $[0,1]$
\item Si $a<0$ et s'il existe $x_0\in]0,1]$ tel que $v'(x_0)=0$, alors on a : $|v(x_0)|=\sqrt{\sqrt{2}|v'(x_0)|}$.
\end{enumerate}
\end{proposition}
\noindent{Preuve :} L'équation $v"=\mbox{sgn}(a)\,v^3$ entraîne que $v"$ appartient à ${\cal C}_b^{2}$.

Intégrons $v"=\mbox{sgn}(a)\,v^3$ sur $[0,1]$. On a :
\begin{eqnarray*}
v'v" &=& \mbox{sgn}(a)\,v'v^3 \\
\frac 12 (v')^2 &=& \frac 14 \mbox{sgn}(a)\,v^4 + \mbox{C}
\end{eqnarray*}
où $C$ est une constante réelle.

Les conditions aux limites entraînent :
\begin{eqnarray*}
v(0) = 0  & \Longrightarrow & \mbox{C} = \frac 12 (v'(0))^2 \\
v(1) = \mbox{sgn}(a)\sqrt{  |a| } & \Longrightarrow & \mbox{C} = \frac 12 (v'(1))^2 - \frac 14 \mbox{sgn}(a)\,a^2
\end{eqnarray*}
Par conséquent, on obtient :
\begin{equation}
\left( v' \right)^2 = \frac 12 \mbox{sgn}(a)\,v^4 + \left( v'(0) \right)^2
\end{equation}

S'il existe $x_0\in]0,1]$ tel que $v'(x_0)=0$, alors on a :
\begin{equation}
\frac 12 \mbox{sgn}(a)\,v(x_0)^4 + \left( v'(0) \right)^2 = 0
\label{formule}
\end{equation}

Si $a>0$ et $v$ non identiquement nulle à droite de zéro, alors il n'existe pas $x_0$ annulant $v'$, qui conserve ainsi le même signe sur $[0,1]$ ; or $v(0)=0$ et $v(1)=1$, donc $v'>0$ sur $[0,1]$ et $v$ est strictement croissante sur $[0,1]$. De plus, $v$ est convexe sur $[0,1]$.

Si $a<0$ et s'il existe $x_0\in]0,1]$ tel que $v'(x_0)=0$, alors $x_0$ vérifie :
\begin{eqnarray*}
  v(x_0)^4 & = & 2\,v'(0)^2 \\
  v(x_0)^2 & = & \sqrt{2}|v'(0)| \\
|v(x_0)| & = & \sqrt{\sqrt{2}|v'(0)| }
\end{eqnarray*}
\qed

\begin{remark} Comportement de la solution au voisinage de zéro.
\begin{enumerate}
  \item Le calcul approché de
$$
v'(0) = \lim_{h\rightarrow 0} \frac{v(h)-v(0)}h = \lim_{h\rightarrow 0} \frac{v(h)}h
$$
permet d'estimer : $v(x_0) = \pm\sqrt{\sqrt{2}|v'(0)| }$.
  \item Quelques conséquences de (\ref{formule}) :
\begin{enumerate}
  \item $\left( v'(1) \right)^2 - \left( v'(0) \right)^2 = \frac 12 \mbox{sgn}(a)\,a^2$
  \item Si $a>0$, alors $|v'| > |v'(0)|$ sur $]0,1]$ et $v'(1)^2 \geq \frac 12 a^2$.
  \item Si $a<0$ et s'il existe $x_0\in]0,1]$ tel que $v'(x_0)=0$, alors $|v(x_0)|>1$ si et seulement si $|v'(0)|>\frac{\sqrt{2}}2$.
  \item Si $a<0$, alors $|v'| < |v'(0)|$ sur $]0,1]$, $v'(1)^2 = v'(0)^2 - \frac 12 a^2$ et $v'(1)^2 \leq v'(0)^2$. De plus, on a : $a^2 = 2\left( v'(0)^2 - v'(1)^2 \right)$. Le calcul approché de $a$ tel que $v'(1)=0$, permet de calculer une valeur approchée de $v'(0)^2$, soit :
$$
v'(1) = 0 \Longrightarrow v'(0)^2 = \frac{a^2}2.
$$
\end{enumerate}
  \item Formule de Taylor de $v$ en zéro.
$$
v(h) = hv'(0) + \frac{h^5}{20}\mbox{sign}(a) v'(0)^3 + o(h^5).
$$
\end{enumerate}
\end{remark}

\section{Formes des courbes représentatives}

Les solutions $u$ du problème (\ref{Sys_u_trois_Annexe}) obtenues par la méthode stochastique ont été comparées aux solutions obtenues par différences finies. Pour $a<0$, la méthode stochastique fournit et privilégie une solution strictement positive.

Indiquons les formes des courbes représentatives de solutions $u$ du pro\-blème (\ref{Sys_u_trois_Annexe}) suivant le signe de $a$. Pour $a<0$, nous ne présentons que les formes de courbes où la solution ne s'annule au plus qu'une fois.

{\bf Cas : $a>0$}

$u$ est positive, croissante, continue, bornée et convexe sur $[0,1]$.

\begin{figure}[h] 
  \centering
  \includegraphics[bb=-223 145 819 696,width=5cm,height=3.7cm,keepaspectratio]{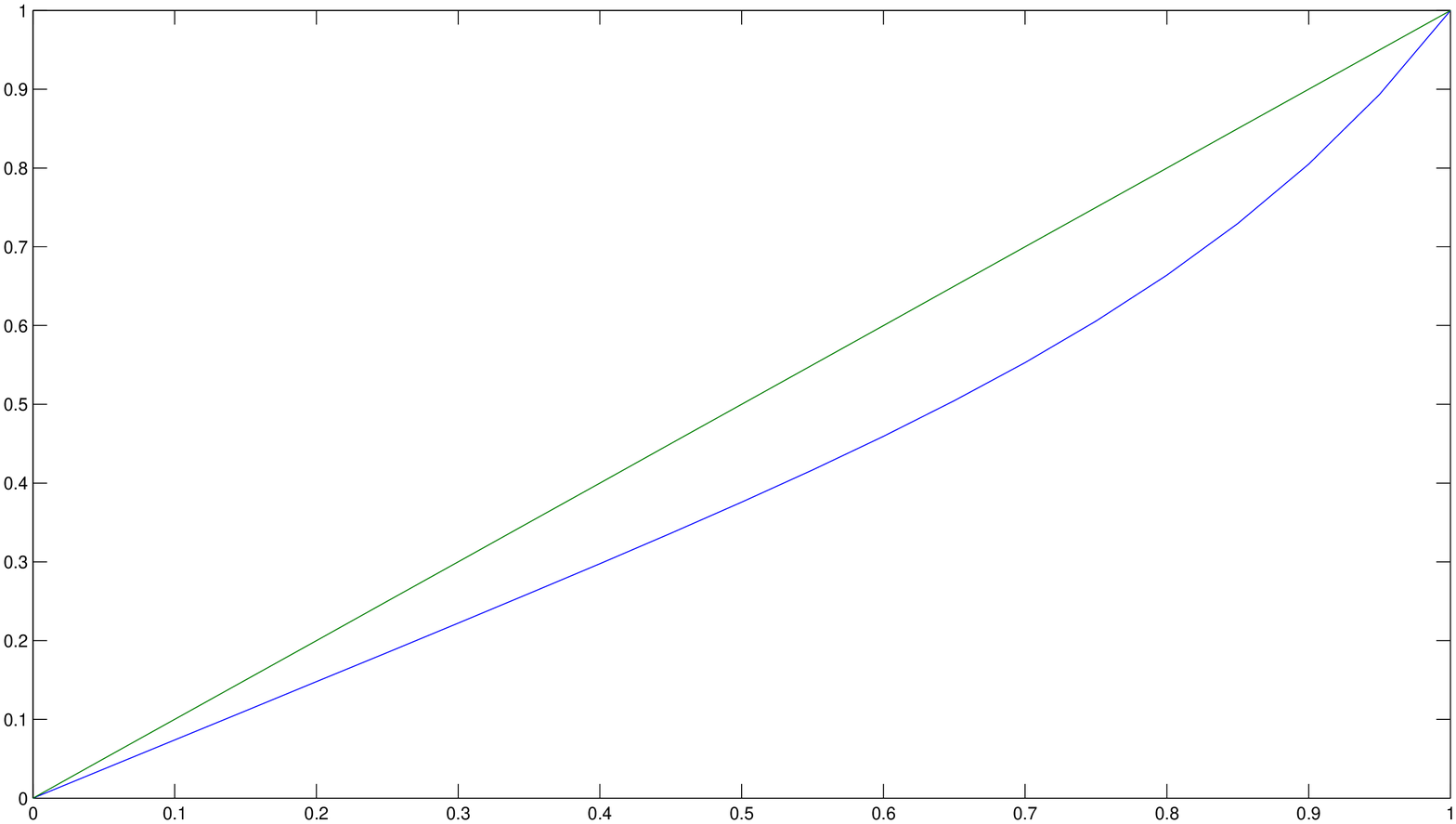}
  \caption{$a=10$}
  \label{fig:Fig_a_positif}
\end{figure}

{\bf Cas : $a<0$}

\begin{figure}
  \includegraphics[bb=-223 145 819 696,width=7cm,height=5.7cm,keepaspectratio]{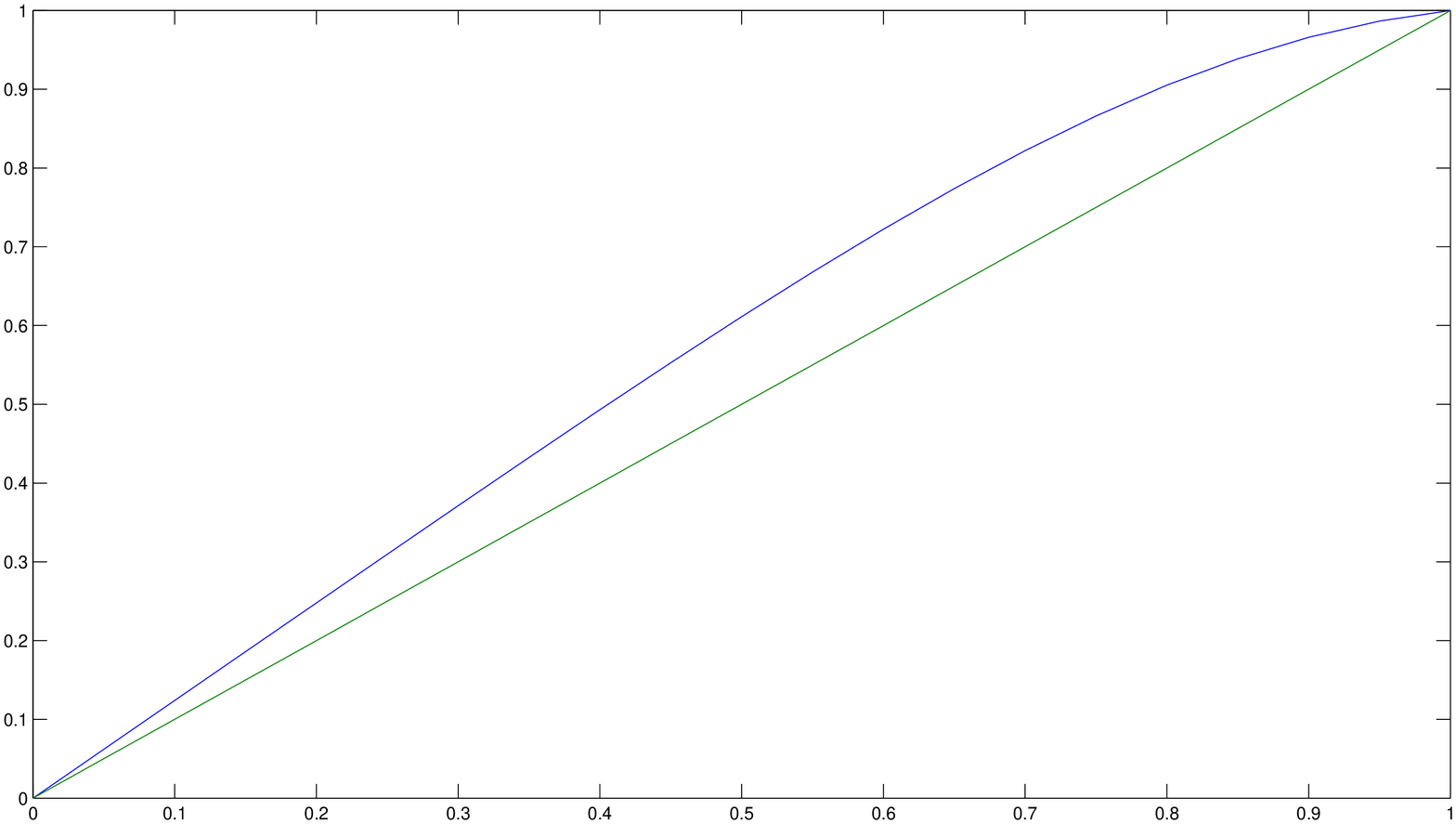}
 \hfill
  \includegraphics[bb=-223 145 819 696,width=7cm,height=5.7cm,keepaspectratio]{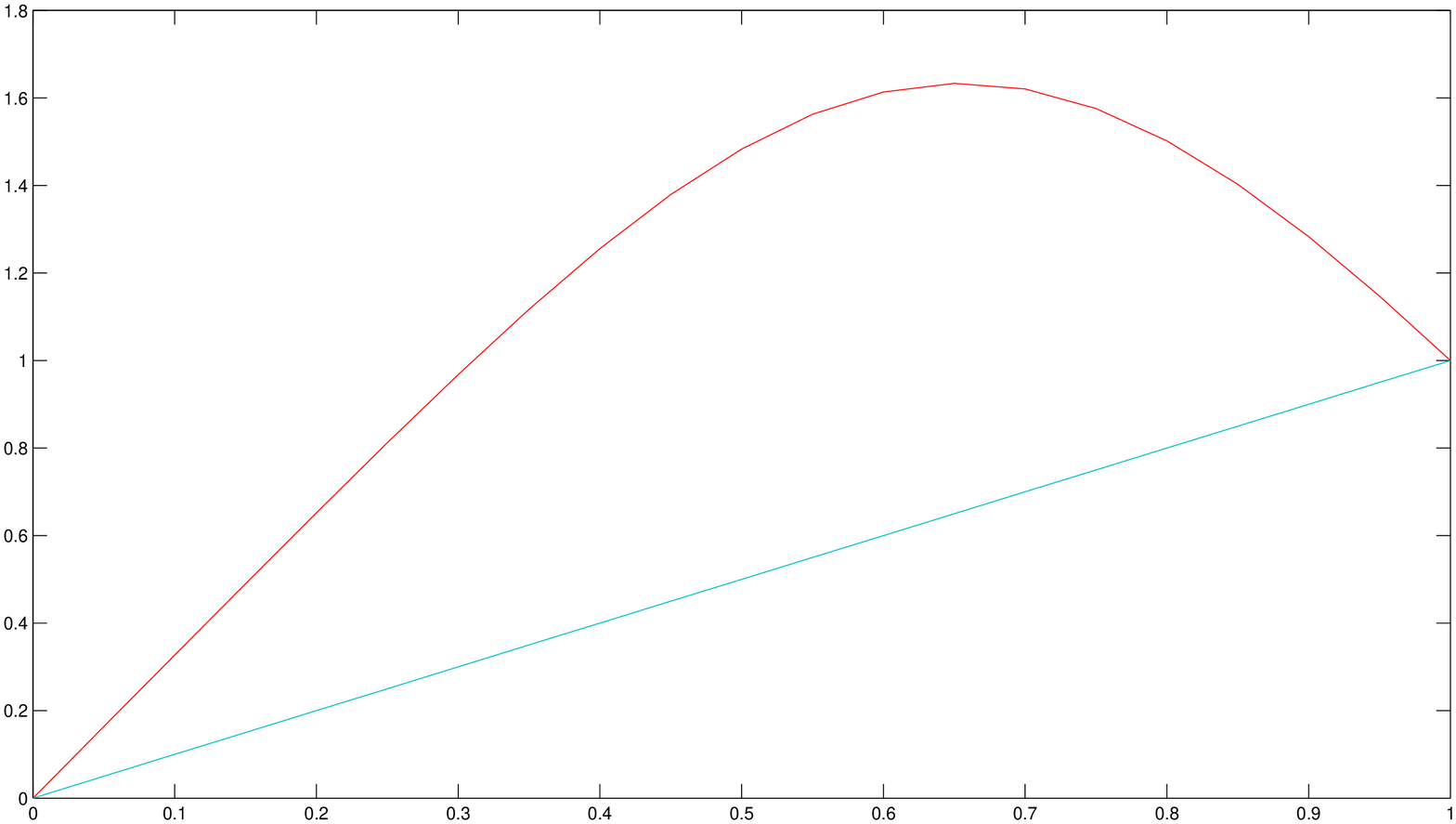}
\caption{Formes $F_1$ et $F_2$ obtenues avec $a<0$}\label{Fig_a_negative1}
\end{figure}

\begin{figure}
\centering
  \includegraphics[bb=-223 145 819 696,width=7cm,height=5.7cm,keepaspectratio]{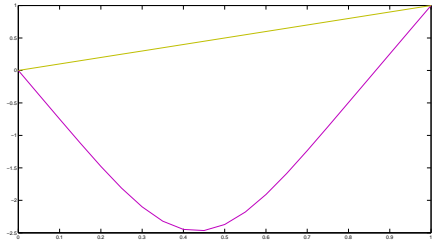}
\caption{Forme $F_3$ obtenue avec $a<0$}\label{Fig_a_negative2}
\end{figure}

Les différentes figures "régulières" obtenues avec le logiciel MATLAB ont la forme des figures \ref{Fig_a_negative1} et \ref{Fig_a_negative2} où le segment représente la fonction identité, pour les cas suivants :
\begin{enumerate}
  \item Forme $F_1$ : $a=-3$ avec \begin{verbatim}bvpinit(linspace(0,1,10),[  0.5 0.5  ])\end{verbatim}   $u$ est positive, strictement croissante, continue, bornée et concave sur $[0,1]$. $u$ atteint son maximum en $1$. $u$, $u'$ et $u"$ ne s'annulent pas sur $]0,1[$. 
  \item Forme $F_2$ : $a=-3$ avec \begin{verbatim}bvpinit(linspace(0,1,10),[  3   3    ])\end{verbatim}   $u$ est positive, continue, bornée et concave sur $[0,1]$, croissante, puis décroissante. $u$ admet un maximum sur $[0,1]$.  $u$ et $u"$ ne s'annulent pas sur $]0,1[$. $u'$ s'annule une seule fois sur $]0,1[$.
  \item Forme $F_3$ : $a=-3$ avec \begin{verbatim}bvpinit(linspace(0,1,10),[ -3  -3    ])\end{verbatim}   $u$ est continue, bornée et convexe sur $[0,1]$, décroissante, puis croissante. $u$ admet un minimum sur $[0,1]$ et s'annule une seule fois sur $]0,1[$. $u"$ s'annule sur $]0,1[$ en même temps que $u$, qui est d'abord convexe, puis concave. $u'$ s'annule une seule fois sur $]0,1[$. 
\end{enumerate}

\begin{figure}[tbp] 
  \centering
  \includegraphics[bb=-6 245 602 596,width=14cm,height=8.05cm,keepaspectratio]{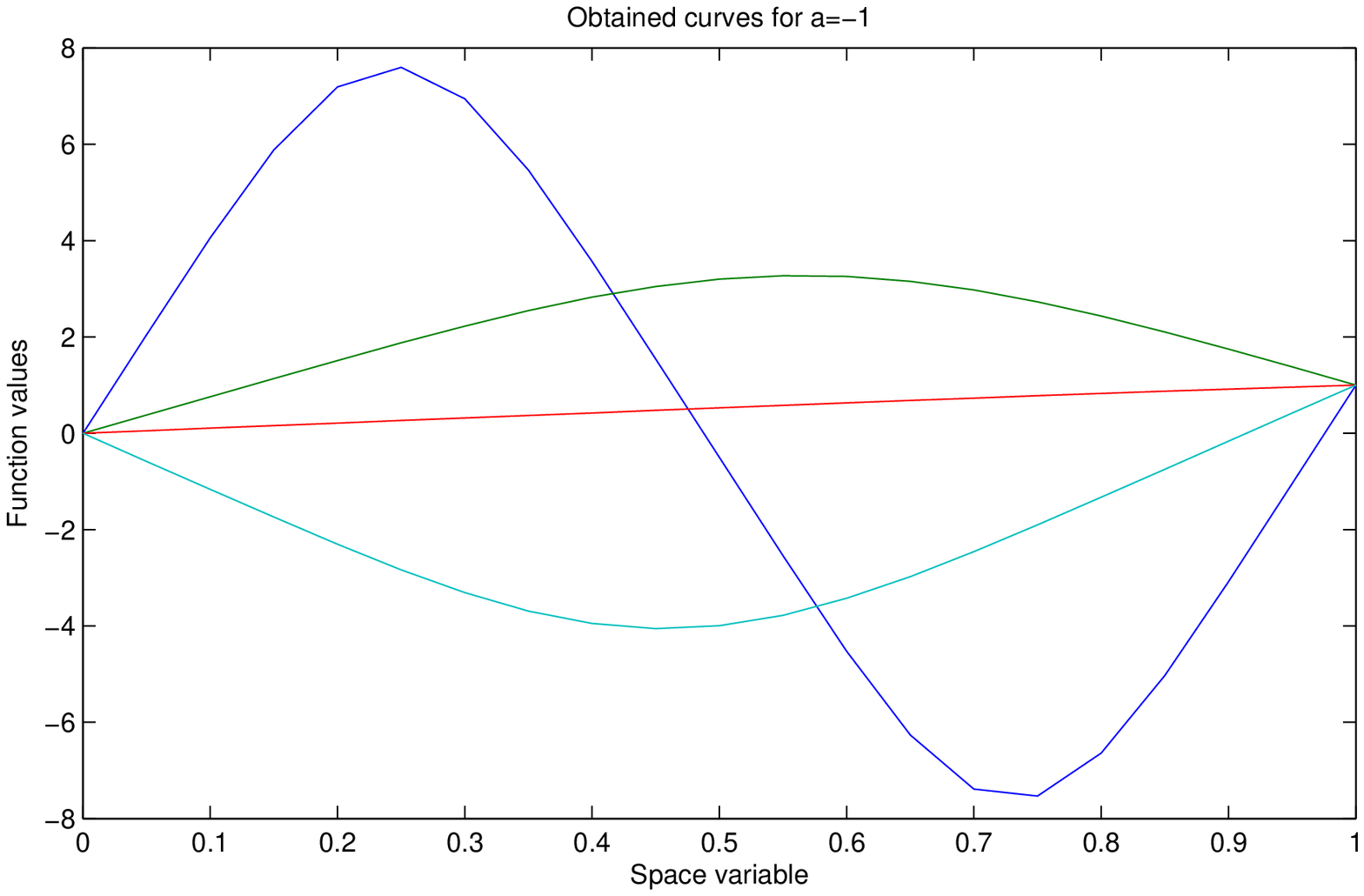}
  \caption{Courbes obtenues avec $a=-1$}
  \label{fig:Fig_Curves_a_moins_1}
\end{figure}

Formes de courbes obtenues à l'aide de la commande \verb!bvpc5! de MATLAB pour $a=-1$ selon les choix de valeurs {\em a priori} de la solution et de sa dérivée en dix points :
\begin{enumerate}
  \item Les valeurs \verb!bvpinit(linspace(0,1,10),[ 1  1  ])! donnent la même courbe que celle obtenue avec la méthode stochastique, à moins de 0.3\% près. Cette courbe est la plus proche de celle de la courbe identité.
  \item Les valeurs \verb!bvpinit(linspace(0,1,10),[ 1.79  1  ])! donnent une courbe convexe avec un minimum proche de $-4$ atteint en $x\simeq 0.45 $.
  \item Les valeurs \verb!bvpinit(linspace(0,1,10),[ 1.8  1  ])! donnent une courbe oscillante de type sinusoïdal où la fonction s'annule 2 fois. Le maximum atteint en $x\simeq 0.275$ vaut environ $7.6$ et le minimum atteint en $x\simeq 0.75$ vaut $-7.6$ environ. La courbe devient polygonale aux voisinages des extrema.
 \item Les valeurs \verb!bvpinit(linspace(0,1,10),[ 1.801  1  ])! donnent une courbe convexe identique à celle obtenue avec \\ \verb!bvpinit(linspace(0,1,10),[ 1.79  1  ])!.
   \item Les valeurs \verb!bvpinit(linspace(0,1,10),[ 1.81  1  ])! donnent une courbe polygonale (irrégulière) où la fonction s'annule 7 fois. L'algorithme de la méthode des différences finies ne converge pas dans ce cas. Cette courbe à oscillations de très grandes amplitudes n'est pas représentée sur la figure \ref{fig:Fig_Curves_a_moins_1}.
\item Les valeurs \verb!bvpinit(linspace(0,1,10),[ 2  1 ])! donnent une cour\-be concave. Le minimum atteint en $x\simeq 0.625$ vaut environ $3.25$.
 \end{enumerate}

\begin{remark}
La méthode stochastique ne fournit que la première courbe ; numériquement, les valeurs des autres solutions calculées par une méthode de différences finies ont été prises comme valeurs initiales de la méthode stochastique, ce qui a conduit soit à une divergence (dernier cas), soit à une convergence vers la première courbe (autre cas) ; en un sens qui est à préciser, la méthode stochastique donne une solution stable (les effets de "bruit blanc" ne font pas diverger l'algorithme pour certaines valeurs {\em a priori}).
\end{remark}

Pour $a=-4$, les quatre courbes suivantes ont été obtenues à l'aide de MATLAB :
\begin{figure}[h] 
  \centering
  \includegraphics[bb=-6 245 602 596,width=14cm,height=8.05cm,keepaspectratio]{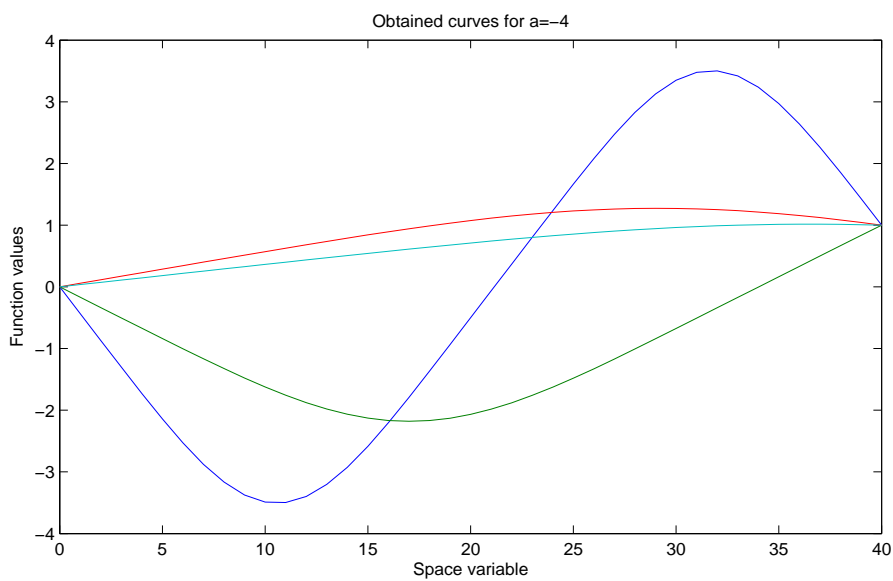}
  \caption{Courbes obtenues avec $a=-4$}
  \label{fig:Fig_a_moins_4}
\end{figure}

\noindent avec successivement :\\
\verb!bvpinit(linspace(0,1,10),[ 0.5  1 ])!,\\ \verb!bvpinit(linspace(0,1,10),[ 1 1 ])!, \\
\verb!bvpinit(linspace(0,1,10),[ -1  -1 ])! et\\ \verb!bvpinit(linspace(0,1,10),[ -1.5  -1 ])!.

Dans ce dernier cas, la solution oscillante de type sinusoïdal est différente de celle obtenue avec $a=-1$.

Seule la solution du premier cas a été obtenue avec la méthode stochastique.

Pour $a\le -5$, aucune solution n'est obtenue à l'aide de la méthode stochastique. Pour $a=-5$ et $a=-6$, les solutions calculées avec MATLAB ont été utilisées comme valeurs initiales de l'algorithme stochastique et ont conduit à la divergence numérique de la méthode stochastique (qui donnent des résultats de plus en plus grands jusqu'à la limite du logiciel). 

\begin{remark}
Pour une même valeur de paramètre $a<0$, les différentes solutions calculées avec MATLAB dépendent des valeurs {\em a priori} que nous avons choisies tandis que l'éventuelle solution calculée avec la méthode stochastique est unique et ne dépend pas des valeurs initiales que nous avons utilisées.
\end{remark}

\chapter{Temps d'atteinte et leurs courbes de régression}
\label{Annexe2}

\section{Introduction}
Etant donnée une partie $V$ de $\mathbb{R}^n$, le temps d'atteinte de $V$ à partir du point $x$ par le processus $x+W$ est la variable aléatoire réelle $\tau_V^x$ définie par\,:
\begin{eqnarray*}
\tau_V^x\,:~\Omega_x    & \longrightarrow &   \mathbb{R}^+ \\
~ ~ ~ ~ ~\omega      & \longmapsto     & \tau_V^x(\omega) = \inf \left\{ s>0 \left| \, X_s(\omega) \in V \right. \right\}
\end{eqnarray*}
où $\Omega_x$ est l'ensemble des trajectoires issues de $x$.

Par la suite, le temps d'atteinte de $V$ est la fonction $\tau_V$ :
\begin{eqnarray*}
\tau_V = E\left[\tau_V^.\right]\,:~G    & \longrightarrow &   \mathbb{R}^+ \\
~ ~ ~ ~ ~x     & \longmapsto     & E\left[\tau_V^x\right].
\end{eqnarray*}
Lorsqu'il n'y a pas d'ambiguité, la partie $V$ n'est pas précisée et $\tau_V$ s'écrit $\tau$.

D'une façon générale, considérons le temps d'atteinte $\tau$ et montrons que $\tau$ dépend du coefficient de diffusion $D$.

\noindent Soit $B(O,r)$ une boule ouverte de $\mathbb{R}^d$, de centre $O$ et de rayon $r$, de frontière ${\cal C}(O,r)$.\\
Considérons $u\,: B(O,r) \longrightarrow \mathbb{R}$ la solution unique du problème\,:
\begin{eqnarray}
\left\{
       \begin{array}{rcll}
 -\frac12 \, D \,\Delta u & = & 1 & \mbox{dans } B(O,r)       \\[2mm]
                        u & = & 0 & \mbox{sur  } {\cal C}(O,r)
       \end{array}
\right.  & \mbox{} & \mbox{}
\label{D tend vers zero}
\end{eqnarray}
où $D > 0$ est supposé constant.

Alors $u$ admet la représentation\index{représentation} (\ref{first representation for u}--\ref{first representation for Y}) avec le processus\index{processus stochas}\,:
\begin{equation}
X_t^x = x + \sqrt{D}\,W_t\,.
\label{processus_diffusion2}
\end{equation}

D'après le raisonnement appliqué par \cite[page 253, avec $D \equiv 1$]{karatzas}, on en déduit que le temps d'atteinte\index{tps:atteinte} $\tau^x$ du complémentaire de $B(O,r)$ par le processus $X_t^x$, à partir d'un point intérieur $x$, vérifie\,:
$$
u(x) = E\left[\, \tau^x \, \right] = \frac{r^2 - |x|^2}{d\,D} \,, \quad x \in B(O,r).
$$
Par conséquent, lorsque $D$ tend vers zéro par valeurs positives, $E\left[ \,\tau^x \, \right]$ tend vers $+\infty$.

De plus, lorsque $D$ décroît, les accroissements\index{accroissements du processus} du processus (\ref{processus_diffusion2}) décrois\-sent avec $\sqrt{D}$ car $dX_t^x=\sqrt{D}~dW_t$\,; ceci limite notre champ d'investigation lors des essais numériques. C'est pourquoi, avant toute résolution numérique, on s'intéresse au comportement de la solution $u$ lorsque $D$ tend vers zéro sur $G$ pour éviter les éventuelles difficultés numé\-ri\-ques.

L'étude du système (\ref{Sys_u_trois}) écrit sous la forme :
\begin{eqnarray}
\left\{
       \begin{array}{rcll}
\frac 12 \frac 1 a \Delta u  & = &  u^3   & \mbox{dans } G=]0;1[       \\[2mm]
       u(0) & = & 0 &  \\[2mm]
       u(1) & = & 1 &  \end{array}
\right.  & \mbox{} & \mbox{}
\end{eqnarray}
conduit à considérer un coefficient de diffusion $D=\frac 1a$ qui tend vers zéro lorsque $a$ tend vers $+\infty$. \`A la limite, la solution du système est la fonction indicatrice $\chi_1$ du point 1 ; la solution du système (\ref{Sys_u_trois}) pour $a\gg 1$ est régularisante de $\chi_1$.

La vitesse de convergence des algorithmes dépend du pas $h$ de la marche aléatoire\,; le processus brownien $X$ discrétisé et le temps d'atteinte $\tau$ de la frontière $\partial G$ en dépendent également. Lors des essais numériques, notre champ d'investigation est limité à des valeurs de $h$ qui permettent au processus d'atteindre la frontière en des temps raisonnables, donc aux vitesses de calcul des ordinateurs utilisés ; pour éviter d'éventuelles difficultés numériques, on s'intéresse au comportement des temps d'atteinte avant toute résolution numérique.

\section{Résultats des essais numériques : temps d'atteinte et courbes de régression}

Les courbes de régression quadratique sont déduites des équations données par l'outil Curve Fitting Tool du logiciel MATLAB.

\subsection{Temps d'atteinte $\tau$ de la frontière de la couronne cir\-cu\-lai\-re}

On s'intéresse au temps d'atteinte restreint à un segment porté par un rayon, en particulier le temps d'atteinte restreint à l'intervalle $[1,3]$ sur l'axe des abscisses.

\begin{figure}[tbp] 
  \centering
  \includegraphics[bb=88 264 505 577,width=12cm,height=10cm,keepaspectratio]{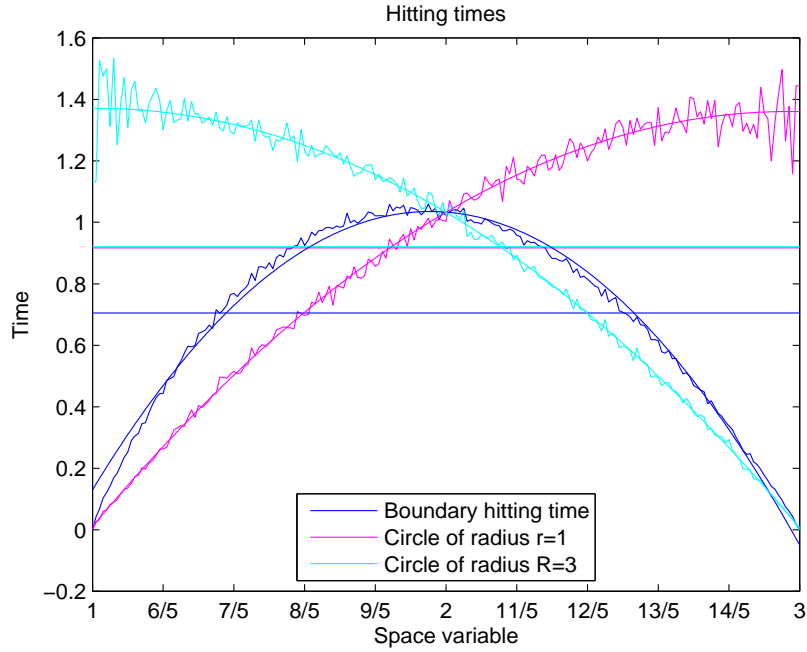}
  \caption{Couronne circulaire -- Temps d'atteinte}
  \label{fig:FigCircularHittingTimes}
\end{figure}

Equation de régression quadratique du temps d'atteinte $\tau$ de la frontière ${\cal C}_1(O,1)\cup{\cal C}_2(O,3) $ :
$$
f(x)=-0.9932  x^2 + 3.883 x -2.76 
$$
Les coefficients sont donnés avec un intervalle de confiance à 95\% :
$$\begin{array}{rclr}
       p1 &=&     -0.9932  &(-1.007, -0.9795)\\
       p2 &=&       3.883  &(3.828, 3.938) \\
       p3 &=&       -2.76  &(-2.812, -2.707)
\end{array}$$
Le temps moyen en seconde pour atteindre la frontière est $E = 0.7052$ où $E$ est la moyenne des $E[\tau^{x_i}]$, $(x_i)$, $i=1$,\ldots,$p$, est une subdivision du segment $[0,1]$ et $\tau^{x_i}$ est le temps d'atteinte d'un mouvement brownien issu de $x_i$.

Equation de régression quadratique du temps d'atteinte du disque interne fermé $D_1={\cal D}(O,1)\cup{\cal C}_1(O,1)$ :
$$
f(x)=-0.3473  x^2 + 2.064  x  -1.706
$$
Les coefficients sont donnés avec un intervalle de confiance à 95\% :
$$\begin{array}{rclr}
       p1 &=&      -0.3473 & (-0.3599, -0.3347)\\
       p2 &=&       2.064  & (2.013, 2.115) \\
       p3 &=&      -1.706  & (-1.754, -1.658)
\end{array}$$

Le pourcentage de marches aléatoires atteignant le disque interne $D_1$ est : $f_1= 40.8521$\%$=0.408521$.

Le temps moyen en seconde pour atteindre le disque interne par les marches aléatoires absorbées par le disque interne (le temps n'est pas pondéré par le nombre de fois où les marches aléatoires atteignent le disque interne $D_1$) : $E_1=0.9169$ où $E_1$ est la moyenne des $E[\tau^{x_i}_1]$, $\tau^{x_i}_1$ est le temps d'atteinte d'un mouvement brownien issu de $x_i$ et atteignant le disque interne $D_1$.

La pondération donne le temps moyen : $T_1=\frac 1p\sum_{i=1}^p E[\tau^{x_i}_1]\times n_i =0.6491$ où $n_i$ est le nombre de fois où la marche aléatoire atteint le disque interne $D_1$ à partir du point $x_i$.

Equation de régression quadratique du temps d'atteinte du domaine com\-plé\-mentaire $\mathbb{R}^2\setminus {\cal D}(O,3)$ de ${\cal D}(O,3)$ :
$$
f(x)=-0.3437  x^2 +  0.6926   x + 1.022 
$$
Les coefficients sont donnés avec un intervalle de confiance à 95\% :
$$\begin{array}{rclr}
       p1 &=&       -0.3437  & (-0.3539, -0.3336) \\
       p2 &=&       0.6926   & (0.6517, 0.7334)\\
       p3 &=&       1.022    & (0.9831, 1.06)
\end{array}$$
Le temps moyen en seconde pour atteindre la partie externe au disque de rayon 3 par les marches aléatoires absorbées par cette partie (le temps n'est pas pondéré par le nombre de fois où les marches aléatoires atteignent la partie externe) : $E_2=0.9210$  où $E_2$ est la moyenne des $E[\tau^{x_i}_2]$, $\tau^{x_i}_2$ est le temps d'atteinte d'un mouvement brownien issu de $x_i$ et atteignant la partie externe au disque de rayon 3.

La pondération donne le temps moyen : $T_2=\frac 1p\sum_{i=1}^p E[\tau^{x_i}_2]\times n_i =0.7438$ où $n_i$ est le nombre de fois où la marche aléatoire atteint la partie externe à partir du point $x_i$.

Compte tenu du pourcentage de marches aléatoires atteignant le disque interne $D_1$, on retrouve le temps moyen d'atteinte de la frontière :
$$
E = f_1 T_1 + (1-f_1) T_2 =  0.7052
$$

\subsection{Temps d'atteinte de la frontière de la couronne rectangulaire 1}

Considérons la couronne rectangulaire : $G = G_1 \setminus G_2$ où $G_1$ et $G_2$ sont les rectangles :
\begin{eqnarray*}
G_1 &=& \left\{\, (x,y) \in \mathbb{R}^2 \mid \,\max(|x|,|y|) < 3 \,\right\} \\
G_2 &=& \left\{\, (x,y) \in \mathbb{R}^2 \mid \,\max(|x|,|y|) \le 1 \,\right\}.
\end{eqnarray*}
On s'intéresse au temps d'atteinte restreint à l'intervalle $[1,3]$ porté par l'axe des abscisses.

\begin{figure}[tbp] 
  \centering
  \includegraphics[bb=88 265 505 575,width=12cm,height=10cm,keepaspectratio]{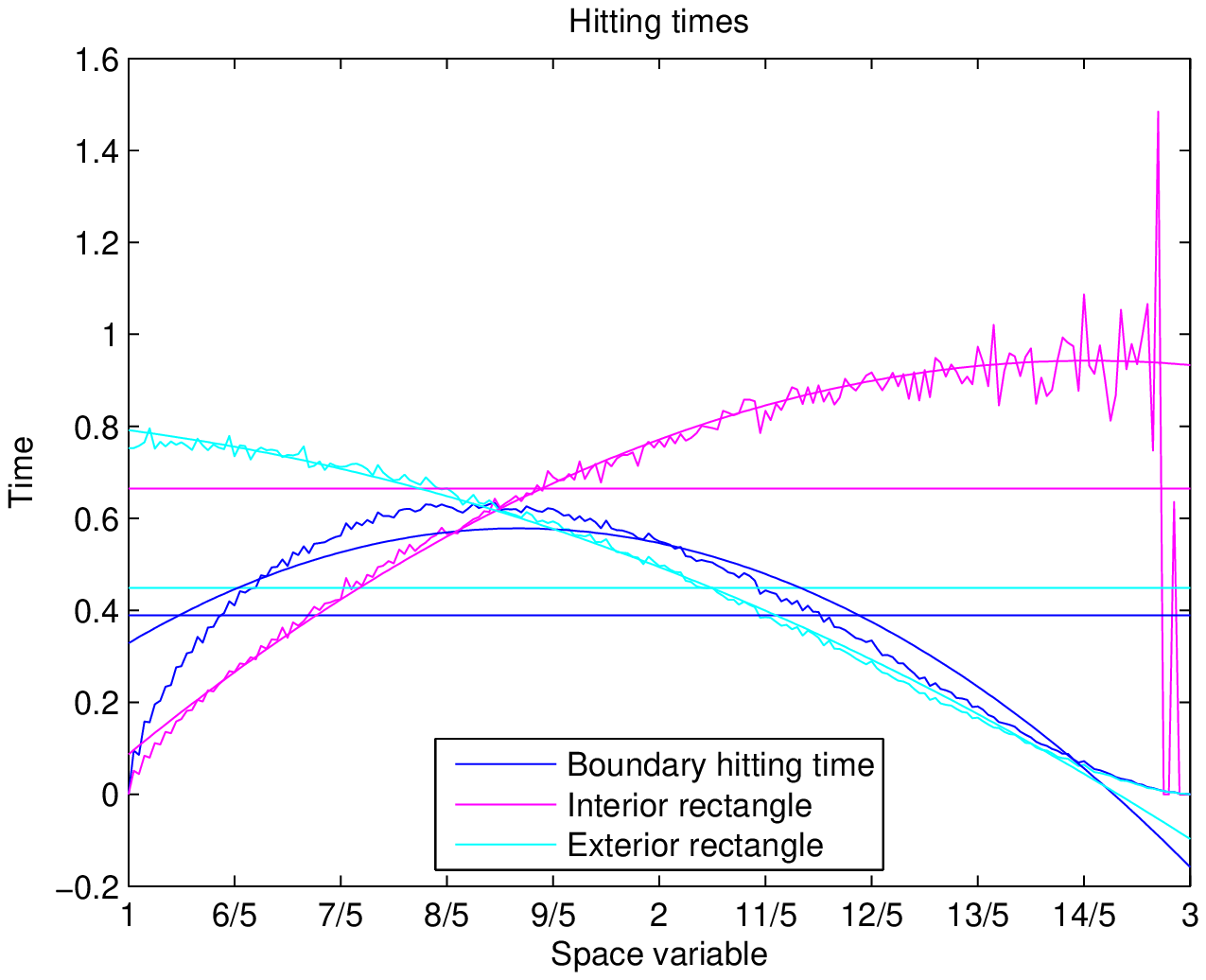}
  \caption{Couronne rectangulaire 1 -- Temps d'atteinte}
  \label{fig:FigR1HittingTimes}
\end{figure}

Equation de régression quadratique du temps d'atteinte de la frontière de $G$ :
$$
f(x)=-0.4609 x^2 + 1.6  x -0.811
$$
Les coefficients sont donnés avec un intervalle de confiance à 95\% :
$$\begin{array}{rclr}
       p1 &=&     -0.4609  &  (-0.487, -0.4349) \\
       p2 &=&      1.6     &  (1.495, 1.705) \\
       p3 &=&     -0.811   &  (-0.9106, -0.7114)
\end{array}$$
Le temps moyen en seconde pour atteindre la frontière est $E = 0.3892$.

Equation de régression quadratique du temps d'atteinte du rectangle fermé $[-1,+1]^2$ :
$$
f(x)=-0.2805 x^2 + 1.535  x -1.175 
$$
Les coefficients sont donnés avec un intervalle de confiance à 95\% :
$$\begin{array}{rclr}
       p1 &=&     -0.2805 & (-0.2978, -0.2632)\\
       p2 &=&      1.535  & (1.466, 1.605) \\
       p3 &=&      -1.706  & (-1.241, -1.109)
\end{array}$$

Le pourcentage de marches aléatoires atteignant le rectangle interne est : $f_1= 27.5465$\%$=0.275465$.

Le temps moyen en seconde pour atteindre le rectangle interne par les marches aléatoires absorbées (le temps n'est pas pondéré par le nombre de fois où les marches aléatoires atteignent le rectangle interne) : $E_1=0.6648$ où $E_1$ est la moyenne des $ E[\tau^{x_i}_1]$, $\tau^{x_i}_1$ est le temps d'atteinte d'un mouvement brownien issu de $x_i$ et atteignant le rectangle interne.

La pondération donne le temps moyen : $T_1=\frac 1p\sum_{i=1}^p E[\tau^{x_i}_1]\times n_i =0.4447$ où $p$ est le nombre de points de la subdivision $(x_i)$, $n_i$ est le nombre de fois où la marche aléatoire atteint le disque interne $D_1$ à partir du point $x_i$.

Equation de régression quadratique du temps d'atteinte du domaine com\-plé\-mentaire $\mathbb{R}^2\setminus ]-3,+3[^2$ de $G_1$ :
$$
f(x)=-0.1469   x^2 +  0.1429   x + 0.7954
$$
Les coefficients sont donnés avec un intervalle de confiance à 95\% :
$$\begin{array}{rclr}
       p1 &=&       -0.1469   & (-0.1542, -0.1396) \\
       p2 &=&       0.1429    & (0.1135, 0.1724)  \\
       p3 &=&       0.7954    & (0.7675, 0.8233)
\end{array}$$
Le temps moyen en seconde pour atteindre la partie externe par les marches aléatoires absorbées par cette partie (le temps n'est pas pondéré par le nombre de fois où les marches aléatoires atteignent la partie externe) : $E_2=0.4488$  où $E_2$ est la moyenne des  $ E[\tau^{x_i}_1]$, $\tau^{x_i}_1$, $\tau^{x_i}_2$ est le temps d'atteinte d'un mouvement brownien issu de $x_i$ et atteignant la partie externe.

La pondération donne le temps moyen : $T_2=\frac 1p\sum_{i=1}^p E[\tau^{x_i}_2]\times n_i =0.3681$ où $p$ le nombre de points $x_i$ d'une subdivision du segment $[0,1]$, $n_i$ est le nombre de fois où la marche aléatoire atteint la partie externe à partir du point $x_i$.

Compte tenu du pourcentage de marches aléatoires atteignant le rectangle interne $[-1,+1]^2$, on retrouve le temps moyen d'atteinte de la frontière :
$$
E = f_1 T_1 + (1-f_1) T_2 =  0.3892
$$

\subsection{Temps d'atteinte de la frontière de la couronne rectangulaire 2}

Considérons la couronne rectangulaire : $G = G_1 \setminus G_2$ où $G_1$ et $G_2$ sont les rectangles :
\begin{eqnarray*}
G_1 &=& \left\{\, (x,y) \in \mathbb{R}^2 \mid \,\sup(|x|,|y|) < 3 \,\right\} \\
G_2 &=& \left\{\, (x,y) \in \mathbb{R}^2 \mid \,\sup(|x|,|y|) \le 1 \,\right\}.
\end{eqnarray*}
On s'intéresse au temps d'atteinte restreint à l'intervalle $[1,3]$ porté par l'axe des abscisses.

\begin{figure}[tbp] 
  \centering
  \includegraphics[bb=88 264 505 577,width=12cm,height=10cm,keepaspectratio]{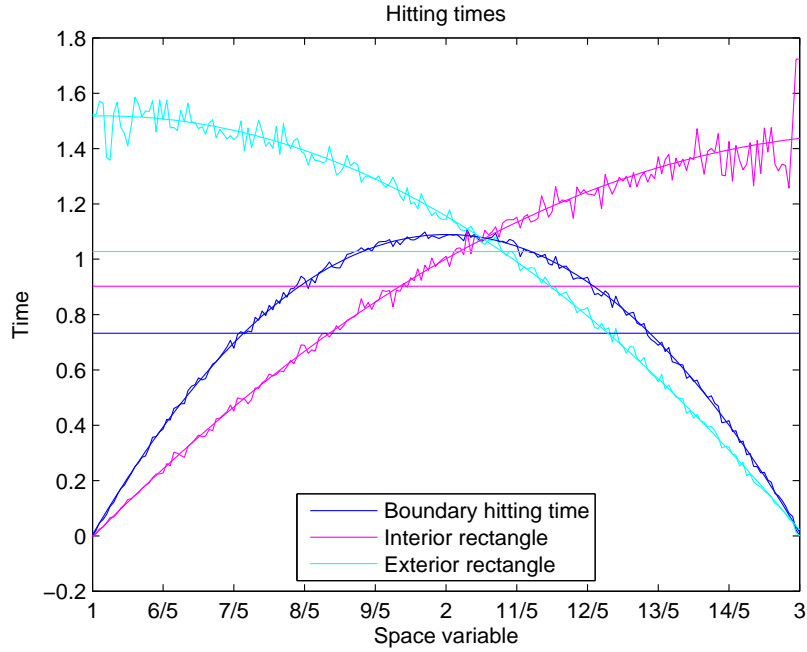}
  \caption{Couronne rectangulaire 2 -- Temps d'atteinte}
  \label{fig:FigR2HittingTimes}
\end{figure}

Equation de régression quadratique du temps d'atteinte de la frontière de $G$ :
$$
f(x)=-1.081 x^2 + 4.329   x -3.244
$$
Les coefficients sont donnés avec un intervalle de confiance à 95\% :
$$\begin{array}{rclr}
       p1 &=&     -1.081  &  (-1.087, -1.075) \\
       p2 &=&      4.329  &  (4.304, 4.354) \\
       p3 &=&     -3.244  &  (-3.267, -3.22)
\end{array}$$
Le temps moyen en seconde pour atteindre la frontière est $E = 0.7324$.

Equation de régression quadratique du temps d'atteinte du rectangle fermé $G_2$ :
$$
f(x)=-0.2849 x^2 +1.86  x -1.58
$$
Les coefficients sont donnés avec un intervalle de confiance à 95\% :
$$\begin{array}{rclr}
       p1 &=&     -0.2849  &  (-0.2961, -0.2738)\\
       p2 &=&      1.86    & (1.815, 1.905) \\
       p3 &=&     -1.58    & (-1.623, -1.537)
\end{array}$$

Le pourcentage de marches aléatoires atteignant le rectangle interne $G_2$ est : $f_1= 46.5136$\%$=0.465136$.

Le temps moyen en seconde pour atteindre le rectangle interne $G_2$ par les marches aléatoires absorbées par $G_2$ (le temps n'est pas pondéré par le nombre de fois où les marches aléatoires atteignent le rectangle interne $G_2$) : $E_1=0.9026$ où $E_1=\sum_{i=1}^p E[\tau^{x_i}_1]$, $p$ le nombre de points $x_i$ d'une subdivision du segment $[0,1]$,  $\tau^{x_i}_1$ est le temps d'atteinte d'un mouvement brownien issu de $x_i$ et atteignant le rectangle interne $G_2$.

La pondération donne le temps moyen : $T_1=\sum_{i=1}^p E[\tau^{x_i}_1]\times n_i =0.6501$ où $n_i$ est le nombre de fois où la marche aléatoire atteint le disque interne $G_2$ à partir du point $x_i$.

Equation de régression quadratique du temps d'atteinte du domaine  com\-plé\-mentaire ${\|x\|_1\ge 3}$ de $G_1$ :
$$
f(x)=-0.3865    x^2 +   0.7979  x + 1.106
$$
Les coefficients sont donnés avec un intervalle de confiance à 95\% :
$$\begin{array}{rclr}
       p1 &=&     -0.3865    & (-0.3981, -0.375) \\
       p2 &=&      0.7979   & (0.7514, 0.8445)  \\
       p3 &=&       1.106   & (1.062, 1.15)
\end{array}$$
Le temps moyen en seconde pour atteindre la partie externe au rectangle $G_1$ par les marches aléatoires absorbées par cette partie (le temps n'est pas pondéré par le nombre de fois où les marches aléatoires atteignent la partie externe) : $E_2=1.0275$  où $E_2$ est la moyenne des  $ E[\tau^{x_i}_2]$ et  $\tau^{x_i}_2$ est le temps d'atteinte d'un mouvement brownien issu de $x_i$ et atteignant la partie externe au rectangle $G_1$.

La pondération donne le temps moyen : $T_2=\frac 1p\sum_{i=1}^p E[\tau^{x_i}_2]\times n_i =0.8040$ où  $p$ est le nombre de points $x_i$ d'une subdivision du segment $[0,1]$ et $n_i$ est le nombre de fois où la marche aléatoire atteint la partie externe à partir du point $x_i$.

Compte tenu du pourcentage de marches aléatoires atteignant le rectangle interne $G_2$, on retrouve le temps moyen d'atteinte de la frontière :
$$
E = f_1 T_1 + (1-f_1) T_2 =  0.7324
$$


\addcontentsline{toc}{chapter}{Bibliographie}
\label{Bibliographie}

\end{document}